\documentclass{article}
\usepackage[utf8]{inputenc}
\usepackage{amssymb,amsmath,amsthm}
\usepackage{graphicx} 
\usepackage{hyperref,url,float}

\newcommand{\E}{\mathbb{E}}
\renewcommand{\P}{\mathbb{P}}
\newcommand{\R}{\mathbb{R}}

\renewcommand{\d}{\mathrm{d}}
\newcommand{\grad}{\nabla}

\newcommand{\g}{>}
\renewcommand{\l}{<}
\newcommand{\vol}{\text{vol}}

\newcommand{\manifolds}{\mathcal{M}}

\newcommand{\G}{\mathcal{G}}

\newcommand{\F}{\mathcal{F}}

\newcommand{\gromovspace}{\mathcal{GS}}
\newcommand{\finitegromovspace}{\mathcal{FGS}}
\newcommand{\gromovdistance}{\mathrm{d}_{\mathcal{GS}}}

\newcommand{\shift}{\mathrm{shift}}
\newcommand{\inv}{\mathrm{inv}}

\newtheorem{theorem}{Theorem}[section]
\newtheorem{example}[theorem]{Example}
\newtheorem{lemma}[theorem]{Lemma}
\newtheorem{corollary}[theorem]{Corollary}
\newtheorem{definition}[theorem]{Definition}

\title{Brownian motion on stationary random manifolds}
\author{Pablo Lessa\thanks{CMAT, Facultad de Ciencias, UdelaR, Iguá 4225, 11400 Montevideo, Uruguay. Email: lessa@cmat.edu.uy}}
\date{}

\begin{document}
\maketitle
\begin{abstract}
We introduce the notion of a stationary random manifold and develop the basic entropy theory for it.  Examples include manifolds admitting a compact quotient under isometries and generic leaves of a compact foliation.  We prove that the entropy of an ergodic stationary random manifold is zero if and only if the manifold satisfies the Liouville property almost surely, and is positive if and only if it admits an infinite dimensional space of bounded harmonic functions almost surely.  Upper and lower bounds for the entropy are provided in terms of the linear drift of Brownian motion and average volume growth of the manifold.  Other almost sure properties of these random manifolds are also studied.
\end{abstract}

\section*{Introduction}

In two articles published in the late 80s (see \cite{kaimanovich1986} and \cite{kaimanovich1988}) Vadim Kaimanovich laid out a plan for the study of asymptotic entropy of Brownian motion on Riemannian manifolds.  By analogy with the, by then, well established theory of Avez entropy for random walks on discrete groups (see \cite{kaimanovich-vershik1983}), he defined an entropy for Riemannian manifolds and outlined how it would relate to the Liouville property, algebraic properties of the fundamental group of the manifold, and volume growth among other things.  The main idea is that the Riemannian metric on the manifold must have some sort of recurrence in order for this statistical approach to work.  The three main cases where one suspects that this recurrence condition is satisfied are: manifolds with transitive isometry group (e.g. the entropy theory for continuous groups was treated by Avez in \cite{avez1976} and Derrienic in \cite{derriennic1985}), manifolds with a compact quotient (see \cite{kaimanovich1986}, \cite{varopoulos1986}, and \cite{ledrappier1996}), and generic leaves of compact foliations in the sense of Lucy Garnett's harmonic measures (see \cite{garnett1983} and \cite{kaimanovich1988}).

The entropy theory of discrete groups has seen its results successively generalized to more general types of graphs than Cayley graphs.  So for example in \cite{kaimanovich-woess2002} the theory is worked out for graphs whose isometry group is transitive while Benjamini and Curien have introduced in \cite{benjamini-curien2012} the concept of a stationary random graph which simultaneously generalizes and includes the case of Cayley graphs and graphs with transitive isomorphism group.   

Following Benjamini and Curien's lead we introduce the concept of a stationary random manifold, which is a rooted random manifold whose distribution is invariant under re-rooting by Brownian motion, and develop the basic theory of entropy for them.  This concept includes manifolds with transitive isometry group, compact quotient, and generic leaves of compact foliations and is also invariant under weak limits (as long as the manifolds involved have uniformly bounded geometry).

Using this concept we carry out Kaimanovich's plan.  That is, we define a non-negative quantity $h(M)$ (the Kaimanovich entropy) associated to each stationary random manifold $M$ which we prove exists and equals zero if and only if $M$ satisfies the Liouville property (i.e. all bounded harmonic functions are constant) almost surely and is positive if and only if almost surely the space of bounded harmonic functions on $M$ is infinite dimensional (see Theorem \ref{entropytheorem}).  We also consider the linear drift $\ell(M)$ of Brownian motion on the random manifold $M$ which exists and is non-negative by virtue of Kingman's subadditive ergodic theorem.  The main result, Theorem \ref{inequalities}, is the chain of inequalities
\[\frac{1}{2}\ell(M)^2 \le h(M) \le \ell(M)v(M)\]
where $v(M)$ is the mean volume growth of the ergodic stationary random manifold $M$.  This result has several applications which we discuss in Section \ref{sectioninequalities}.

In Section \ref{sectionbrownianmotion} we construct a path space and corresponding a Brownian motion on any stationary Riemannian manifold.  Using this concept we prove some results about linear drift which where used in the proof of Theorem \ref{inequalities} (see Corollary \ref{citethisfromchapter2}), obtain a Furstenberg type formula for linear drift (Theorem \ref{furstenbergformula}), and a weak result (Theorem \ref{ghystypetheorem} and Corollary \ref{ghyscorollary}) in the spirit of Ghys work on surface laminations (see \cite{ghys1995}).

In the case of stationary random Hadamard manifolds with pinched negative curvature under a reversibility hypothesis (which is automatically satisfied in the case of a single manifold with compact quotient), we improve the lower bound for entropy to $2\ell(M)^2 \le h(M)$  (see Theorem \ref{lsquarelessthanh}), a result first proved by Kaimanovich in the case with compact quotient (see \cite[Theorem 10]{kaimanovich1986}).  The same inequality is valid for all manifolds with compact quotient without any curvature restrictions as shown by Ledrappier in \cite[Theorem A]{ledrappier2010}, however we have no analogous result for stationary random manifolds.

We have organized the article as follows.  In the first section we treat results about a single fixed Riemannian manifold.  In the second, we treat result on stationary random manifolds without the need for Brownian motion (except indirectly through the heat kernel).  Finally, in the third section we define Brownian motion on a stationary random manifold and explore its main properties.

The results in article are included in the author's Phd. Thesis.

\subsection*{Acknowledgements}

I would like to thank François Ledrappier and Matilde Martinez for advising me during the past few years.  Special thanks are also due to Vadim Kaimanovich and  Jesús A. Álvarez Lopez for helping to improve this work (also, the example of a manifold whose Brownian motion is not steady we give in Theorem \ref{vadimsexample} is due to Vadim).  I benefited greatly from conversations with Fernando Alcalde, Sébastién Álvarez, Christian Bonatti, Matías Carrasco, Yves Coudene, Gilles Courtois, Françoise D'albo, Bertrand Deroin,  Raphaël Krikorian, and Andrés Sambarino.  I would also like to thank Samuel Senti and the other organizers and participants of the conference Probability in Dynamics 2014 at URFJ, Río de Janeiro.  This work was partially financed by CSIC and SNI/ANII.
\section{Liouville properties and Zero-one laws on Riemannian manifolds}

In this section we establish the basic results linking the asymptotic behavior of Brownian motion with potential theory of the underlying Riemannian manifold (see Theorem \ref{zerooneliouvilletheorem} and Lemma \ref{tailbijection}).  We give an example, due to Vadim Kaimanovich, of a manifold on which the tail and invariant $\sigma$-algebras are not equivalent for Brownian motion (see Theorem \ref{vadimsexample}),   and show that any such example necessarily has unbounded geometry (see Theorem \ref{steadybrownian}).  To conclude we give proofs of the main facts about mutual information (i.e. Theorem \ref{informationtheorem} and Lemma \ref{finitedimensionallemma}) which will be used later on to prove results about Kaimanovich entropy (see Theorem \ref{entropytheorem}).

The results in this section are not new but they are not available (to the best of the author's knowledge) in the literature in the form needed for our applications.  Analogous results for discrete time Markov chains were treated in \cite{kaimanovich1992} and can be used to give alternative proofs of some of the results discussed here.   Mutual information and its relationship to entropy was discussed in  \cite{derriennic1985} in the context of random walks on continuous groups.

\subsection{Brownian motion and the backward heat equation}

\subsubsection{Laplacian, heat semigroup, and heat kernel}

Consider a connected complete $d$-dimensional Riemannian manifold $M$.  We denote by $\Delta$ the self-adjoint extension of the non-positive definite Laplacian on $L^2(M)$ and by $p(t,x,y)$ the heat kernel.  We recall that $p(t,x,y) = p(t,y,x)$ and $p$ is a smooth positive function on $(0,+\infty)\times M \times M$ such that $u(t,x) = P^tf(x) = \int p(t,x,y)f(y)\d y$ is the unique solution to the heat equation $\partial_t u = \Delta_x u$ with initial data $u(0,x) = f(x)$ for any given $f \in L^2(M)$.   Existence of the heat kernel and its basic properties are treated in detail in \cite{grigoryan2009}.

It can be shown that $\int p(t,x,y)\d y \le 1$ for all $t \g 0$.  If the last integral is always equal to $1$ then one says that $M$ is stochastically complete.

The Euclidean plane minus one point $\R^2 \setminus \lbrace 0\rbrace$ is an example of a stochastically complete manifold which is not geodesically complete.  An example of a complete but not stochastically complete manifold can be obtained by endowing the plane $\R^2$ with the Riemannian metric given in polar coordinates by
\[\d s^2 = \d r^2 + p(r)^2 \d \theta^2,\]
for some function $p$ satisfying $p(r) = e^{r^3}$ for all $r$ large enough.

\subsubsection{Brownian motion}

On $\R$ one has $p(t,x,y) = (4\pi t)^{-1/2} \exp\left(-\frac{(x-y)^2}{4t}\right)$.  We notice that the density of the time $t$ of a standard Brownian motion starting at $x \in \R$ can be written as $p(t/2,x,y)$.  When passing to a Riemannian manifold we have decided to keep this factor of $\frac{1}{2}$ which distinguishes the heat kernel from the transition density function of Brownian motion.  To avoid confusion we keep the notation $q(t,x,y) = p(t/2,x,y)$.

Given a stochastically complete manifold $M$ and a point $x \in M$ we define Weiner measure $\P_x$ starting at $x$ on the space $\Omega = C([0,+\infty),M)$ of continuous paths from $[0,+\infty)$ to $M$ as the unique Borel measure (the topology being that of uniform convergence on closed intervals) such that for all Borel sets $A_1,\ldots,A_n \subset M$ and all positive times $t_1 \l \cdots \l t_n$ the probability
\[\P_x(\omega_{t_1} \in A_1,\cdots, \omega_{t_n} \in A_n)\]
of the set of paths $\omega \in \Omega$ which visit each $A_i$ at the corresponding time $t_i$ is given by the integral
\[\int\limits_{A_1\times \cdots \times A_n}q(t_1,x,x_1)q(t_2-t_1,x_1,x_2)\cdots q(t_n-t_{n-1},x_{n-1},x_n)\d x_1\cdots \d x_n.\]

A Brownian motion with initial distribution $\mu$ (a Borel probability on $M$) is defined to be an $M$ valued stochastic process whose distribution is given by
\[\int \P_x \d \mu(x).\]

With the above definition one can prove the existence of manifold valued Brownian motion via Kolmogorov's continuity theorem using further properties of the heat kernel (i.e. upper bounds in terms of distance).

We will also use an alternative description of Brownian motion (usually attributed to Eells, Elsworthy, and Malliavin, e.g. see \cite[pg. 75]{hsu2002}) as a diffusion on the orthogonal frame bundle $O(M)$.

Consider the smooth vector fields $V_i, i = 1,\ldots,d$ on $O(M)$ such that the flow of $V_i$ applied to a frame $X = (x,v_1,\ldots,v_d) \in O(M)$ (here $x \in M$ and the $v_i$ form an orthonormal basis of the tangent space at $x$) moves the basepoint along the geodesic with initial condition $v_i$ and transports the frame horizontally.   Then any solution to the Stratonovich stochastic differential equation
\[\d X_t = \sum\limits_{i = 1}^d V_i(X_t)\circ \d W_t^i\]
driven by a standard Brownian motion $(W^1_t,\ldots,W^d_t)$ in $\R^d$, projects to a Brownian motion on $M$.

The equivalence of these two approaches is established in \cite[Propositions 3.2.2 and 4.1.6]{hsu2002}.

\subsubsection{Zero-one laws}

For each $t \ge 0$ define the $\sigma$-algebra $\F_t$ of events occurring before time $t$ as the Borel subsets of $\Omega = C([0,+\infty),M)$ generated by the open sets of the topology of uniform convergence on the interval $[0,t]$.  Similarly we let $\F^t$ be the $\sigma$-algebra of events occurring after time $t$ which is generated by the open sets of the topology of uniform convergence on compact subsets of the interval $[t,+\infty)$.  Events belonging to all $\F^t$ are called tail events and form the tail $\sigma$-algebra defined by
\[\F^\infty = \bigcap\limits_{t \ge 0}\F^t.\]

We recall that any Borel subset in $\Omega$ can be approximated (meaning the probability of the symmetric difference can be made arbitrarily small) by a finite disjoint unions of events of the form
\[\lbrace \omega \in \Omega:  \omega_{t_1} \in A_1,\ldots, \omega_{t_n} \in A_n\rbrace\]
where $t_1 \l \ldots \l t_n$ and the sets $A_i$ are Borel subsets of $M$.  Similarly each set in $\F_t$ can be approximated by finite disjoint unions of events of the above form with $t_n \le t$ and each set in $\F^t$ by events of the above form with the restriction $t_1 \ge t$.

The Markov property allows one to express the probability of a tail event with respect to the measure $\P_x$ as averages over $y$ of the probabilities with respect to $\P_y$ of a `shifted' event.  More concretely let $\shift^t:\Omega \to \Omega$ be defined for $t \ge 0$ by 
\[(\shift^t\omega)_s = \omega_{t+s},\]
one has the following property.

\begin{lemma}\label{backwardequation}
Let $M$ be a complete connected and stochastically complete Riemannian manifold.  For each tail event $A$ the function
\[u(t,x) = \P_x(\shift^t(A))\]
solves the backward heat equation
\[\partial_t u(t,x) = -\frac{1}{2}\Delta u(t,x).\]
\end{lemma}
\begin{proof}
Fix $T \g 0$ and set $v(t,x) = u(T-t,x)$ for each $t \in (0,T)$ and $x \in M$. By applying the Markov property one obtains
\[v(t,x) = \int q(t,x,y)\P_y\left(\shift^{T} \shift^{T-t}A\right)\d y = \int p(t/2,x,y)v(0,y)\d y\]
which implies that $\partial_t v(t,x) = \frac{1}{2}\Delta v(t,x)$ from which the desired result follows.
\end{proof}

We say that an event $A \subset \Omega$ is trivial if it has probability $0$ or $1$ with respect to all measures $\P_x$.  Brownian motion on $M$ is said to satisfy the zero-one law if all tail events are trivial.  Lemma \ref{backwardequation} allows one to show that triviality of a tail event for $\P_x$ is independent of the choice of $x \in M$ (in particular the zero-one law can be verified at a single $x \in M$).

\begin{corollary}\label{triviality}
Let $M$ be a complete connected and stochastically complete Riemannian manifold and $A \subset \Omega$ be a tail event.  Then $A$ is trivial if and only if $\shift^t(A)$ has probability $0$ or $1$ with respect to some $\P_x$ for some $t \ge 0$.
\end{corollary}
\begin{proof}
Apply the maximimum principle to $u(t,x)$ defined in Lemma \ref{backwardequation}.
\end{proof}

An event $A$ is said to be invariant if $(\shift^t)^{-1}(A) = A$ for all $t \ge 0$ (this implies $\shift^t(A) = A$ since the shift maps are surjective).  The $\sigma$-algebra of all invariant events is denoted by $\F^{\inv}$.  Since invariant events are also tail events one may apply Lemma \ref{backwardequation} to obtain the following.

\begin{corollary}
Let $M$ be a complete connected and stochastically complete Riemannian manifold.  For each invariant event $A$ the function
\[v(x) = \P_x(A)\]
is harmonic (i.e. $\Delta v(x) = 0$ for all $x$).
\end{corollary}

We say Brownian motion is ergodic on $M$ if all invariant events are trivial.  By Corollary \ref{triviality} ergodicity is equivalent to triviality of all invariant events with respect to a single probability $\P_x$.

\subsubsection{Liouville properties}

A manifold $M$ is said to satisfy the Liouville property (some times we just say $M$ is Liouville) if it admits no non-constant bounded harmonic functions.  Similarly we say $M$ is backward-heat Liouville if it admits no non-constant bounded solutions to the backward heat equation (defined for all $t \ge 0$).

\begin{theorem}\label{zerooneliouvilletheorem}
Let $M$ be a complete connected and stochastically complete Riemannian manifold.  Then $M$ is backward-heat Liouville if and only if its Brownian motion satisfies the zero-one law.  Similarly, $M$ is Liouville if and only if its Brownian motion is ergodic.
\end{theorem}
\begin{proof}
Suppose $M$ is backward-heat Liouville and $A$ is a tail event.   Then by Lemma \ref{backwardequation} the function
\[u(t,x) = \P_x(\shift^t A)\]
solves the backward equation and by hypothesis must be constant.

Given times $t_1 \l \cdots \l t_n$ and Borel sets $A_1,\ldots A_n \subset M$ we calculate using the Markov property (which is possible because $A \in \F^{t_n}$) to obtain that the probability
\[\P_x\left(\omega_{t_i} \in A_i\text{ for }i=1,\ldots,n\text{ and }\omega \in A\right)\]
of the trajectory belonging to $A$ while hitting each $A_i$ at the corresponding time $t_i$ is given by
\[\int\limits_{A_1\times \cdots\times A_n}q(t_1,x,x_1)\cdots q(t_n - t_{n-1},x_{n-1},x_n)\P_{x_n}(\shift^{t_n}A)\d x_1\cdots \d x_n\]
which since $u(t,x)$ is constant yields
\[\P_x\left(\omega_{t_i} \in A_i\text{ for }i=1,\ldots,n\right)\P_x\left(\omega \in A\right).\]

This implies that $A$ is independent from $\F_t$ for all $t$ so that $A$ is independent from itself and must have probability $0$ or $1$.  We conclude that if $M$ is backward-heat Liouville then its Brownian motion satisfies the zero-one law (notice that the proof mimics that of the classical zero-one law).

The same argument shows that if $M$ is Liouville then its Brownian motion is ergodic.

On the other hand if there is a bounded backward solution $u(t,x)$ defined for all $t \ge 0$ then 
\[u(t,\omega_t)\]
is a bounded martingale with respect to any $\P_x$.  Since $u(t,\cdot)$ is not constant (otherwise $u$ would be constant) the random variable $u(t,\omega_t)$ is not almost-surely constant with respect to $\P_x$.  On the other hand the martingale convergence theorem implies that the limit
\[f(\omega) = \lim\limits_{t \to +\infty}u(t,\omega_t)\]
exists almost surely with respect to $\P_x$ and that its conditional expectation to $\F_t$ is $u(t,\omega_t)$.  This shows that $f$ is not almost-surely constant with respect to $\P_x$ and, since $L$ is tail measurable, there are non-trivial tail events.

In the case where one assumes that there is a non-constant bounded harmonic function $v(x)$ one has that $u(t,x) = v(x)$ is a bounded backward solution independent of $t$.  The same argument above works with the additional fact that the limit $f$ is shift invariant and hence yields non-trivial invariant events.
\end{proof}

We conclude this subsection reexamining the last part of the previous proof (i.e. the construction of bounded tail measurable function $f:\Omega \to \R$ starting from a bounded backward solution $u(t,x)$).  In view of Corollary \ref{triviality} all the measures $\P_x$ are mutually absolutely continuous when restricted to the tail $\sigma$-algebra $\F^{\infty}$.  We call the measure class of any and all $\P_x$ the harmonic measure class on $\F^{\infty}$.  We say a tail measurable function $f:\Omega \to \R$ is invariant if $f\circ \shift^t = f$ for all $t \ge 0$.

\begin{lemma}\label{tailbijection}
Let $M$ be a complete connected and stochastically complete Riemannian manifold.  There is a one to one correspondence associating to each bounded solution $u(t,x)$ to the backward equation $\partial_tu(t,x) = -\frac{1}{2}\Delta u(t,x)$ the bounded tail measurable function
\[f_u(\omega) = \lim\limits_{t \to +\infty} u(t,\omega_t)\]
considered up to modifications on zero-measure sets with respect to the harmonic measure class.  Furthermore $f_u$ can be modified on a null set with respect to the harmonic measure class so that it is shift invariant if and only if $u(t,x) = v(x)$ for some bounded harmonic function $v:M\to \R$.
\end{lemma}
\begin{proof}
First of all we fix $x \in M$ and notice that $u(t,\omega_t)$ is a bounded martingale with respect to $\P_x$ so that the limit $f_u(\omega)$ exists $\P_x$-almost surely.  Since the existence of the limit $f_u$ is a tail event this implies that $f_u$ is well defined almost surely with respect to the harmonic measure class on $\F^{\infty}$.

We will now show that $u \mapsto f_u$ is injective.  

For this purpose suppose $f_u = f_v$ almost surely with respect to $\P_x$.  By the martingale convergence theorem the conditional expectation of $f_u$ to $\F_t$ with respect to $\P_x$ is given by
\[\E_x\left(f_u|\F_t\right) = u(t,\omega_t)\]
and similarly for $f_v$ so that one has for each $t \ge 0$ that
\[u(t,\omega_t) = v(t,\omega_t)\]
for $\P_x$ almost every $\omega \in \Omega$.   Since $\omega_t$ has a strictly positive density $q(t,x,\cdot)$ under $\P_x$ and the functions $u(t,\cdot)$ and $v(t,\cdot)$ are continuous this implies that $u(t,\cdot) = v(t,\cdot)$ for each $t$ so that $u = v$ as claimed.

If $u(t,x) = v(x)$ for some harmonic function $v$ then 
\[f_u(\omega) = \lim\limits_{t \to +\infty} v(\omega_t) = \lim\limits_{t \to +\infty}v(\omega_{t+s}) = f_u(\shift^s\omega)\]
almost surely with respect to the harmonic measure class so $f_u$ can be modified on a zero measure set to be invariant.

Reciprocally assume that $f_u$ is shift invariant.  One has
\[\lim\limits_{t \to +\infty}u(t,\omega_t) = f_u(\omega) = f_u(\shift^s\omega) = \lim\limits_{t \to +\infty}u(t,\omega_{t+s}) = \lim\limits_{t \to +\infty}u(t-s,\omega_t).\]
Setting $u_s(t,x) = u(t-s,x)$\begin{footnote}{One can extend $u_s$ to $t \le s$ uniquely using the heat equation.}\end{footnote} one obtains that $f_u = f_{u_s}$ so that by the previously established injectivity $u = u_s$.  Since this works for all $s$ we obtain that $u(t,x) = v(x)$ for some harmonic function $v$.

It remains only to show that the map $u \mapsto f_u$ is surjective.

By Lemma \ref{conditionalexpectation} below for each $t$ and $x$ there is a probability $\P_{(t,x)}$ on $\F^t$ which satisfies
\[\P_{(t,x)}\left(A\right) = \P_x\left(\shift^t(A)\right).\]

Denoting by $\E_{(t,x)}$ the expectation with respect to $\P_{(t,x)}$ and setting 
\[u(t,x) = \E_{(t,x)}\left(f(\omega)\right)\]
one has by the martingale convergence theorem and Lemma \ref{conditionalexpectation} that $f = f_u$.  Hence $u \mapsto f_u$ is surjective as claimed.
\end{proof}

\begin{lemma}\label{conditionalexpectation}
 Let $M$ be a complete connected and stochastically complete Riemannian manifold.  For each $t \ge 0$ the map $\shift^t$ is a bijection between the $\sigma$-algebras $\F^T$ and $\F^{T-t}$ on $\Omega$ for all $T \ge t$.  In particular each $\shift^t$ is a bijection on $\F^\infty$.  
 
Furthermore, denoting by $\P_{(t,x)}$ the unique probability on $\F^t$ which satisfies 
\[\P_{(t,x)}\left(A\right) = \P_x\left(\shift^t(A)\right)\]
for all $A \in \F^t$ one has that the conditional expectation of any bounded and tail measurable function $f: \Omega \to \R$ to the $\sigma$-algebra $\F_t$ relative to the probability $\P_{x_0}$ ($x_0$ being any chosen point in $M$) is given by
\[\E_{x_0}\left(f(\omega)|\F_t\right) = u(t,\omega_t)\]
where $u(t,x) = \E_{(t,x)}(f(\omega))$ is the expectation of $f$ relative to $\P_{(t,x)}$ for all $t \ge 0$ and $x \in M$.
\end{lemma}
\begin{proof}
We had glossed over this point earlier (e.g. in Lemma \ref{backwardequation}) but the continuity of $\shift^t$ does not imply that if $A \in \F^t$ then $\shift^t(A)$ is Borel.  

However, if $\omega \in A$ for some $A \in \F^t$ then all continuous paths which coincide with $\omega$ after time $t$ also belong to $A$.  This property implies (valid for all $t \ge 0$) that $\shift^t$ is a bijection between $\F^T$ and $\F^{T-t}$ (even though $\shift^t$ certainly is not injective as a function on $\Omega$) for all $T \ge t$.

The second claim amounts to establishing the fact that 
\begin{equation}\label{conditionalexpectationequation}
\E_{x_0}\left(f(\omega)1_A(\omega)\right) = \E_{x_0}\left(u(t,\omega_t)1_A(\omega)\right) 
\end{equation}
for all $A \in \F_t$.

Suppose first that $f = 1_B$ for some $B \in \F^\infty$ and 
\[A = \lbrace \omega \in \Omega: \omega_{s_i} \in A_i, i =1,\ldots,m\rbrace\]
where the $A_i$ are Borel subsets of $M$ and $t_1 \l \cdots \l t_n \le t$.

Then one has
\begin{align*}
\E_{x_0}&\left(f(\omega)1_A(\omega)\right) = \P_{x_0}\left(A \cap B\right)
\\ &= \int\limits_{A_1\times\cdots \times A_m}q(s_1,x_0,x_1)\cdots q(t - s_{m},x_m,y)\P_y\left(\shift^t B\right)\d x_0\cdots \d x_m \d y
\\ &= \E_{x_0}\left(u(t,\omega_t)1_A(\omega)\right).
\end{align*}

Since any $A \in \F_t$ can be approximated (with respect to $\P_{x_0}$) by finite disjoint unions of events of the above form we have established the claim for bounded tail measurable functions that are indicators of a tail set.  

For the general case notice that given two functions for which Equation \ref{conditionalexpectationequation} holds one has that the equation holds for any linear combination of them.  Furthermore, if $f$ is the monotone limit of a sequence of non-negative functions for which Equation \ref{conditionalexpectationequation} is known to hold then by the monotone convergence theorem the equation holds for $f$ as well.  This proves that the claim holds for all bounded tail measurable functions.
\end{proof}

\subsection{A bounded backward heat solution}

If $v$ is a bounded harmonic function on a manifold $M$ then $u(t,x) = v(x)$ is a bounded solution to the backward heat equation.  The following construction, due to Vadim Kaimanovich (see also \cite[pg. 23]{kaimanovich1992}), yields an example of a manifold admitting bounded solutions to the backward heat equation which are not of the above form.

\begin{figure}[H]
\centering
\begin{minipage}{\textwidth}
\includegraphics[width=\textwidth]{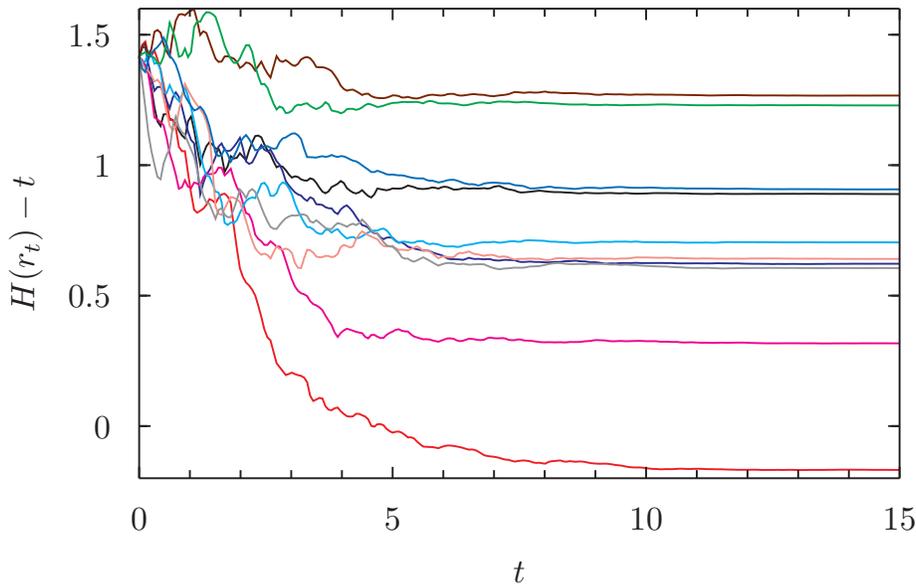}
\end{minipage}
\caption{Ten trajectories of the process $H(r_t)-t$.}
\label{radialprocessfigure}
\end{figure}

\begin{theorem}\label{vadimsexample}
Consider the smooth Riemannian metric $g$ on the plane $\R^2$ which in polar coordinates has the form
\[\d s^2 = d r^2 + p(r)^2 \d \theta^2\]
with $p(r) = re^{\frac{1}{2}r^2}$.   There exists a smooth bounded function $u(t,x)$ which solves the backward heat equation with respect to this metric and such that $u(t,\cdot)$ is not harmonic for any $t \in \R$.
\end{theorem}
\begin{proof}
To see that such an expression in polar coordinates yields a smooth metric at the origin of $\R^2$ we calculate explicitly the coefficients of the metric (letting $e_1,e_2$ be the canonical basis of $\R^2$ and $(x,y) = (r\cos(\theta),r\sin(\theta))$) and obtain
\[g_{11} = g(e_1,e_1) = 1 + y^2 (p(r)^2/r^2 - 1)/r^2\]
\[g_{12} = g(e_1,e_2) = -xy(p(r)^2/r^2 - 1)/r^2\]
\[g_{22} = g(e_2,e_2) = 1 + x^2 (p(r)^2/r^2 - 1)/r^2,\]
so the claim follows because $(p(r)^2 /r^2 - 1)/r^2$ can be extended analytically to $r = 0$ (just consider the power series of $p(r)$).

Consider a solution $r_t$ to the Ito differential equation
\[\left\lbrace\begin{array}{ll}r_0 = 1\\ \d r_t = \d X_t + f(r_t)\d t\end{array}\right.\]
where $f(r) = \frac{1}{2}p'(r)/p(r) = (r + 1/r)/2$ and $X_t$ is a standard Brownian motion on $\R$.  If one sets
\[\tau_T = \int_0^T \frac{1}{f(r_t)^2}\d t\]
and
\[\theta_t = Y_{\tau_t}\]
where $Y_t$ is an Euclidean Brownian motion independent from $X_t$, then $(r_t\cos(\theta_t),r_t\sin(\theta_t))$ is a Brownian motion for the metric $g$ (see \cite[Example 3.3.3]{hsu2002}).

We will show that there is a non-trivial tail event for the process $r_t$ which is not shift invariant.

For this purpose notice that the fact that $f(r) \ge 1$ implies that
\[r_T = 1 + X_T + \int_0^T f(r_t)\d t \ge \left(1 + X_T + T\right)^+\]
for all $T$ where $x^+ = x$ if $x \g 0$ and $0$ otherwise.

Next set $H(r) = \log(1 + r^2)$ and notice that $H'(r) = h(r) = 1/f(r)$ if $r \g 0$.  By the Ito formula one has
\[\d H(r_t) = h(r_t)\d X_t + (1 + \frac{1}{2}h'(r_t))\d t.\]

We will show that the limit
\[L = \lim\limits_{t \to +\infty} H(r_t) - t\]
exists almost surely.  Clearly $L$ is tail measurable with respect to the filtration associated to $r_t$ and is not shift invariant (replacing $r_t$ by $r_{t+s}$ changes the value of $L$ by $s$ as well).  If we show that $L$ is not almost surely constant then there are non-trivial tail events (of the form $\lbrace L \g a\rbrace$) which are not shift invariant.

Notice that
\[H(r_T) - T = \int_0^T h(r_t)\d X_t + \frac{1}{2}\int_0^Th'(r_t)\d t.\]

Using the inequality $r_t \ge \left(1 + X_t + t\right)^+$ one obtains that for almost all trajectories one eventually has $r_t \g t/2$.  Combined with the fact that $|h'(r)| = O(1/r^2)$ when $r \to +\infty$ one obtains that
\[\int_0^{+\infty} |h'(r_t)| \d t \l +\infty\]
almost surely.

To show that the martingale part of $H(r_T) - T$ converges it suffices to show that its variance is bounded.  By the Ito isometry one has
\[\E\left[\left(\int_0^T h(r_t)\d X_t\right)^2\right] = \int_0^T \E\left[h(r_t)^2\right]\d t.\]

To bound the integrand we separate into two cases according to whether $|X_t| \g t/2$ or not and obtain (for $t \g 2$ using that $h \le 1$ and that $h$ is decreasing on $r \g 1$)
\[\E\left[h(r_t)^2\right] \le \P\left[|X_t| \g t/2\right] + h(t/2)^2 = \P\left[ |X_1| \g \sqrt{t}/2\right] + h(t/2)^2.\]
The right hand side is integrable because the first term decreases exponentially while the second is of order $O(1/t^2)$.

Hence we have established that the limit $L$ of $H(r_t) - t$ exists almost surely when $t \to +\infty$.  To complete the proof it remains to show that the random variable $L$ is not almost surely constant (see Figure \ref{radialprocessfigure} below for evidence supporting this claim).

Suppose that $L$ were almost surely equal to a constant $C$.  Let the stopping time $\sigma$ for $r_t$ be minimal among those with the property that $r_\sigma = 1$ and $r_t = 2$ for some $t \l \sigma$.  One always has $\sigma \g 0$ and, by the Varadhan-Stroock support theorem, there is a positive probability that $\sigma$ is finite.  The Markov property implies that on the set with $\sigma \l \infty$ one has
\[C = \lim\limits_{t \to +\infty}H(r_t) - t = \lim\limits_{t \to +\infty}H(r_{\sigma+ t}) - t = C + \sigma\]
contradicting the fact that $\sigma$ is positive.
\end{proof}

In the above example the radial process $r_t$ grows super-linearly so that $\tau_t$ converges almost surely as $t \to +\infty$ and hence so does $\theta_t$.  Events of the form $\theta_\infty = \lim_{t \to +\infty}\theta_t \in [a,b]$ are invariant and therefore may be used to define non-constant bounded harmonic functions.   As far as the author is aware the question of whether there exists a manifold satisfying the Liouville property but admitting non-constant bounded solution to the backward heat equation is still open.

\subsection{Steadyness of Brownian motion}

Recall that a Riemannian manifold is said to have bounded geometry if its injectivity radius is positive and its sectional curvature is bounded in absolute value.  We show that examples such as the one given in the previous subsection must have unbounded geometry (in the example the curvature at distance $r$ from the origin is $-(3+r^2)$).

Following Kaimanovich we say Brownian motion on $M$ is steady if every tail event can be modified on a null set with respect to the harmonic measure class on $\F^\infty$ to be invariant.  This is equivalent (via Lemma \ref{tailbijection}) to the property that every bounded solution to the backward heat equation is of the form $u(t,x) = v(x)$ for some bounded harmonic function $v$.

The following result was proved in the case $M$ has a compact quotient under isometries by Varopoulos (see \cite[pg. 359]{varopoulos1986}).  A more general result with no assumption on the injectivity radius of $M$ was announced by Kaimanovich with a proof sketch (see \cite[Theorem 1]{kaimanovich1986}).

\begin{theorem}\label{steadybrownian}
Let $M$ be connected Riemannian manifold with bounded geometry.  Then $M$ is stochastically complete and Brownian motion on $M$ is steady.  In particular every bounded solution $u(t,x)$ to the backward heat equation defined for all $t \ge 0$ is of the form $u(t,x) = v(x)$ for some harmonic function $v$.
\end{theorem}

The so-called zero-two law is a sharp criteria for equivalence of the tail and invariant $\sigma$-algebras of Markov chains (see \cite{derriennic1976}).  In our situation it amounts to the statement that
\[\sup\limits_{x \in M}\left\lbrace \lim\limits_{t \to +\infty}\int | p(t+\tau,x,y) - p(t,x,y)|\d y\right\rbrace\]
is either equal to $0$ or to $2$ for all $x \in M$ and all $\tau \g 0$ and furthermore the limit is $0$ if and only if Brownian motion is steady.

We will verify that the above limit cannot be $2$ in Lemma \ref{boundedgeometrylemma} below.  From this, steadiness of Brownian motion follows from the zero-two law.  A proof which does not rely on the zero-two law will be given at the end of this subsection.

\begin{lemma}\label{boundedgeometrylemma}
Let $M$ be a connected Riemannian manifold with bounded geometry.  For each $\tau \g 0$ there exists $\epsilon_\tau \g 0$ such that
\[\int |q(t+\tau,x,y) - q(t,x,y)|\d y \le 2-\epsilon_\tau\]
for all $x \in M$ and $t \ge \tau$.  In particular, if $u(t,x)$ is a solution to the backward equation $\partial_t u = -\frac{1}{2}\Delta_x u$ bounded by $1$ in absolute value then 
\[|u(t+\tau,x) - u(t,x)| \le 2 - \epsilon_\tau\]
for all $t \ge 0$ and all $x \in M$.
\end{lemma}
\begin{proof}
Let $K \g 0$ be a finite bound for the absolute value of all the sectional curvatures of $M$ and $\rho \g 0$ be strictly less than the injectivity radius at all points of $M$ and the diameter of the $d$-dimensional sphere of constant curvature $K$.

Fix $x \in M$ and let $\psi:\R^d \to M$ be a normal parametrization at $x$, i.e. $\psi(v) = \exp_x\circ L (v)$ where $\exp_x:T_xM \to M$ is the Riemannian exponential map at $x$ and $L:\R^d \to T_xM$ is a linear isometry between $\R^d$ (endowed with the usual inner product) and the tangent space $T_xM$ (with the inner product given by the Riemannian metric on $M$).

Consider the metric of constant curvature $-K$ ball $B_\rho(0)$ of radius $\rho$ centered at $0$ in $\R^d$ of the form $\d s^2 = \d r^2 + p(r)\d \theta^2$ where $\d \theta^2$ is the standard Riemannian metric on the unit sphere $S^{d-1} \subset \R^d$ and one sets 
\[p(r) = \sinh(\sqrt{K}r).\]

We denote by $\varphi(t,\cdot)$ the probability density of the time $t$ of Brownian motion started at $0$ and killed upon first exit from $B_\rho$ with respect to the constant curvature metric above.  The only fact about $\varphi$ we need is that it is everywhere positive on $B_\rho$ for all $t$.

Let $q_\rho^K(t,x,y)$ be defined for $y$ in the open ball $B_\rho(x)$ of radius $\rho$ centered at $x$ by
\[q^K_\rho(t,x,y) = \varphi(t,\psi^{-1}(y))\]
where $\psi^{-1}(y)$ is the unique preimage of $y$ in $B_\rho(0)$.

Theorem 1 of \cite{debiard-gaveau-mazet1976} states that for all $y \in B_\rho(x)$ one has
\[q^K_\rho(t,x,y) \le q_\rho(t,x,y)\]
where $q_\rho(t,x,\cdot)$ is the probability density of the time $t$ of Brownian motion on $M$ started at $x$ and killed upon first exit from the ball of radius $\rho$ centered at $x$.

Also one has $q_\rho(t,x,y) \le q(t,x,y)$ since the probability of Brownian motion on $M$ going from $x$ to a small neighborhood of $y$ in time $t$ diminishes if one demands that it never exit the ball of radius $\rho$ centered at $x$ before that.  Therefore one has
\[q^K_\rho(t,x,y) \le q(t,x,y)\]
for all $y \in B_\rho(x)$.

Define $\epsilon_\tau(x)$ by the equation
\[\epsilon_\tau(x) = \int\limits_{B_\rho(x)} \min(q^K_\rho(\tau,x,y),q^K_\rho(2\tau,x,y)) \d y.\]

Let $\omega$ be the Euclidean volume form on $B_\rho$ and $\lambda(p)\omega$ be the pullback of the volume form of $M$ under $\psi$.  Since the sectional curvature of $M$ is bounded from above by $K$ by \cite[Theorem 27]{petersen2006} one has $\lambda(p) \ge \sin(\sqrt{K}r)^{d-1}$ for all $p$ at distance $r$ from $0$ in $B_\rho$.  Since $\epsilon_\tau(x)$ can be calculated by integrating a fixed positive function on $B_\rho$ with respect to the form $\lambda(p)\omega$ one obtains that $\epsilon_\tau = \inf\lbrace \epsilon_\tau(x), x \in M\rbrace$ is positive.

Since $q(\tau,x,\cdot) \ge q_\rho^K(\tau,x,\cdot)$ and $q(2\tau,x,\cdot) \ge q_\rho^K(2\tau,x,\cdot)$ one obtains the following
\[\int |q(\tau,x,y) - q(2\tau,x,y)|\d y \le \int q(\tau,x,y) + q(2\tau,x,y) - \min(q(\tau,x,y),q(2\tau,x,y))\d y \le 2 - \epsilon_\tau.\]

From this it follows for all $t \ge 0$ that
\[\int |q(t+2\tau,x,y) - q(t+\tau,x,y)|\d y \le \int |q(2\tau,x,z) - q(\tau,x,z)|q(t,z,y)\d z \d y \le 2 - \epsilon_\tau\]
as claimed.

To conclude we observe that if $u$ satisfies the backward equation and is bounded in absolute value by $1$ then one has
\[|u(t+\epsilon,x) - u(t,x)| = | \int (q(\tau,x,y)-q(2\tau,x,y))u(t+2\tau,y)\d y| \le 2 - \epsilon_\tau\]
which concludes the proof.
\end{proof}

As mentioned above one can prove Theorem \ref{steadybrownian} from the previous lemma using the zero-two law.  The proof below relies instead on the bijection between bounded tail measurable functions and solutions to the backward equation (see Lemma \ref{tailbijection}).

\begin{proof}[Proof of Theorem \ref{steadybrownian}]
Let $v_r(x)$ denote the volume of the ball of radius $r$ centered at a point $x \in M$.  The lower curvature bound implies that $v_r$ is less than or equal to the volume of a ball of radius $r$ in hyperbolic space of constant curvature $-K$ (see \cite[Lemma 35]{petersen2006}).  In particular one has
\[\int_1^{+\infty} \frac{r}{\log(v_r)}\d r = +\infty\]
so that $M$ is stochastically complete by \cite[Theorem 11.8]{grigoryan2009}.

Suppose that Brownian motion on $M$ is not steady.  Then one can find a non-trivial non-invariant (even up to modifications on null-sets with respect to the harmonic measure class) tail set $A$ and $\tau \g 0$ such that $B = \shift^\tau(A)$ is disjoint from $A$.  It follows from Lemma \ref{triviality} that $B$ is also non-trivial.

Consider the tail function defined by $f(\omega) = 1_A(\omega) - 1_B(\omega)$.  By Lemma \ref{tailbijection} there exists a bounded solution $u$ to the backward equation such that $f = f_u$.  By Lemma \ref{conditionalexpectation} one knows that $u$ is bounded by $1$ in absolute value almost everywhere and by continuity of $u$ this holds everywhere.

Notice that for almost every $\omega$ with respect to the harmonic measure class one has:
\[\lim\limits_{t \to +\infty}u(t,\omega_t) = f(\omega)\]
and
\[\lim\limits_{t \to +\infty}u(t-\tau,\omega_t) = \lim\limits_{t \to +\infty}u(t,\omega_{t+\tau}) = f(\shift^\tau(\omega)).\]

In particular by choosing such a generic path in $B$ one obtains that there exists $\omega \in \Omega$ such that
\[\lim\limits_{t \to +\infty}u(t,\omega_t) = -1\]
and
\[\lim\limits_{t \to +\infty}u(t-\tau,\omega_t) = 1.\]

This implies that there exist values of $t$ and $x$ such that $u(t-\tau,x) - u(t,x)$ is arbitrarily close to $2$, contradicting Lemma \ref{boundedgeometrylemma}.
\end{proof}

\subsection{Mutual information}

Suppose $M$ is a stochastically complete manifold whose Brownian motion satisfies the zero-one law.  Then given $x \in M$, $A \in \F_t$ and $B \in \F^\infty$ one has that $A$ and $B$ are independent under $\P_x$ i.e. $\P_x(A \cap B) = \P_x(A)\P_x(B)$.  The converse is also true, i.e. if each tail event $B$ is independent from the events in $\F_t$ for all $t$ then the Brownian motion on $M$ satisfies the zero-one law (Proof: As in the proof of the classical zero-one law, one approximates $B$ by events in $\F_t$ to show that it is independent from itself and hence trivial).

A, perhaps convoluted, but useful way of rephrasing this is the following:  Consider the function $\omega \mapsto (\omega,\omega)$ from $\Omega$ to $\Omega\times \Omega$.  Since this function is continuous one can push forward $\P_x$ to obtain a probability $\widehat{\P}_x$ on $\Omega\times \Omega$.  The measure $\widehat{\P}_x$ describes the joint distribution of two copies of the same Brownian motion on $M$.  On the other hand the probability $\P_x \times \P_x$ on $\Omega \times \Omega$ describes the joint distribution of two independent Brownian motions on $M$ starting at $x$.    The two probabilities $\widehat{\P}_x$ and $\P_x \times \P_x$ are very different (e.g. they are mutually singular).  However, assuming the zero-one law is satisfied, if one restricts them both to the $\sigma$-algebra $\sigma\left(\F_t\times \F^\infty\right)$ generated by sets of the form $A \times B$ with $A \in \F_t$ and $B \in \F^\infty$ then they coincide.  In fact, Brownian motion on $M$ satisfies the zero-one law if and only 
if $\widehat{\P}_x$ and $\P_x \times \P_x$ coincide when restricted to $\sigma\left(\F_t \times \F^\infty\right)$ for all $t \ge 0$.

The mutual information between two random variables is a non-negative number which is zero if and only if they are independent.  Given a $\sigma$-algebra $\F$ of Borel sets in $\Omega$ one may consider the identity map $\omega \mapsto \omega$ as a random variable from $\Omega$ endowed with the Borel $\sigma$-algebra to $\Omega$ endowed with $\F$, and hence one may define mutual information between $\sigma$-algebras.

Concretely, given $x \in M$ we define the mutual information between $\F_t$ and $\F^T$ (where $0 \le t \le T$ and possibly $T = \infty$) under $\P_x$ as
\[I_x\left(\F_t,\F^T\right) = \sup\left\lbrace \sum\limits_{i = 1}^{n}\log\left(\frac{\widehat{\P}_x(A_i)}{\P_x\times\P_x(A_i)}\right)\widehat{\P}_x(A_i) \right\rbrace\]
where the supremum is taken over all finite partitions $A_1,\ldots,A_n$ of $\Omega\times \Omega$ with each $A_i$ belonging to $\sigma\left(\F_t\times \F^T\right)$.  One may interpret the result as a measure of how much the behavior of Brownian motion after time $T$ (or the tail behavior if $T = \infty$) depends on what happened before time $t$.

The fact that $I_x\left(\F_t,\F^T\right)$ is always non-negative and is zero if and only if $\widehat{P}_x$ and $\P_x \times \P_x$ coincide on $\sigma\left(\F_t \times \F^T\right)$ follows from Jensen's inequality applied to the strictly convex function $-\log$ (see \cite[Lemma 3.1]{gray2011} for details).

Mutual information was used to unify results about random walks on discrete and continuous groups by Derriennic, in particular he established several results analogous to the Theorem below in that context (e.g. see \cite[Section III]{derriennic1985}).  In the case of a manifold with a compact quotient under isometries similar results to those below where established by Varopoulos (see \cite[Part I.5]{varopoulos1986}).  Results of this type where also announced by Kaimanovich both in the case when $M$ has a compact quotient and when $M$ is a generic leaf of a compact foliation (e.g. \cite[Theorem 2]{kaimanovich1986} and \cite[Lemma 1]{kaimanovich1988}).  In the context of discrete time Markov chain similar results are discussed in detail in \cite[Section 3]{kaimanovich1992}.
 
\begin{theorem}\label{informationtheorem}
Let $M$ be a complete connected and stochastically complete Riemannian manifold.  Then Brownian motion on $M$ satisfies the zero-one law if and only if $I_x\left(\F_t,\F^\infty\right) = 0$ for some $t \g 0$ and $x \in M$.  Furthermore, the following properties hold for all $x \in M$ and $0 \l t \le T \l \infty$:
\begin{enumerate}
\item $I_x(\F_t,\F^T) = \int \log\left(\frac{q(T-t,x_1,x_2)}{q(T,x,x_2)}\right)q(t,x,x_1)q(T-t,x_1,x_2)\d x_1 \d x_2.$
\item The function $T \mapsto I_x\left(\F_t,\F^T\right)$ is non-increasing and satisfies the inequality $I_x\left(\F_t,\F^\infty\right) \le \lim\limits_{T \to +\infty}I_x\left(\F_t,\F^T\right)$ with equality if some $I_x\left(\F_t,\F^T\right)$ is finite.
\end{enumerate}
\end{theorem}
\begin{proof}

If Brownian motion on $M$ satisfies the zero-one law then $\F_t$ is independent from $\F^\infty$ under $\P_x$ for all $x$ and therefore $I_x\left(\F_t,\F^\infty\right) = 0$.

Assume now that there is some $x \in M$ with $I_x\left(\F_t,\F^\infty\right) = 0$ and fix $B \in \F^\infty$.  We must show that $\P_x(B)$ is either $0$ or $1$.

For this purpose fix $s \l t$ and and an open subset $U$ of $M$ and notice that
\[\P_x(\omega_s \in U,  \omega \in B) = \int\limits_{U}q(s,x,y)\P_{(s,y)}(B) \d y.\]

On the other hand by hypothesis the above is also equal to
\[\P_x(\omega_s \in U)\P_x(B) = \int\limits_{U}q(s,x,y)\P_x(B)\d y\]
from which one obtains that
\[\P_{(s,y)}(B) = \P_x(B)\]
for almost all $y \in M$.  Since $u(s,y) = \P_{(s,y)}(B)$ is a solution to the backward equation it must be constant and equal to $\P_x(B)$ for all $s$ and $y$.

Consider now the set of paths where with $\omega_{t_i} \in A_i$ for all $i =1,\ldots,n$ where $t_1 \l \cdots \l t_n$ and the $A_i$ are Borel subsets of $M$.  One may calculate using the above to obtain
\begin{align*}\P_x\left(A \cap B\right) &= \int\limits_{A_1\times\cdots\times A_n}q(t_1,x,x_1)\cdots q(t_n-t_{n-1},x_{n-1},x_n)\P_{(t_n,x_n)}(B) \d y
\\ &= \P_x(A)\P_x(B)
\end{align*}
so that $B$ is independent from all events $A$ of this form.  Since $B$ may be approximated with respect to $\P_x$ by finite disjoint unions of events of the form $A$ above this shows that $B$ is independent from itself and hence has probability equal to $0$ or $1$ as claimed.

We will now establish the integral formula for $I_x\left(\F_t,\F^T\right)$ (property 1 above).  

The so-called Gelfand-Yaglom-Perez Theorem (see \cite[Theorem 2.1.2]{pinsker1964} and the translator's notes on page 23 or \cite[Lemma 7.4]{gray2011} for further detail) implies that
\[I_x\left(\F_t,\F^T\right) = \int f_T(\omega^1,\omega^2)\log(f_T(\omega^1,\omega^2)) \d(\P_x \times \P_x)(\omega^1,\omega^2)\]
where $f_T$ is the Radon-Nikodym derivative of $\widehat{P}_x$ restricted to $\sigma\left( \F_t,\F^T\right)$ with respect to $\P_x\times \P_x$ restricted to the same $\sigma$-algebra.  The formula then follows by substituting the explicit formula for $f_T$ that we will establish below in Lemma \ref{radonderivative}.  Notice that, because $x \mapsto x\log(x)$ is bounded from below, the integral formula always makes sense regardless of convergence considerations, but may assume the value $+\infty$.

We will now prove property 2 of the statement.

To begin notice that when $T$ increases the set of partitions used to define $I_{x}\left(\F_t,\F^T\right)$ decreases, hence the supremum taken over all such partitions decreases as well.  This implies that $T \mapsto I_{x}\left(\F_t,\F^T\right)$ is decreasing and also that
\[I_x\left(\F_t,\F^\infty\right) \le \lim\limits_{T \to +\infty}I_x\left(\F_t,\F^T\right).\]

Now assume that $I_x(\F_t,\F^{T_0})$ is finite and set $f = f_{T_0}$.  Notice that by definition of the Radon-Nikodym derivative one has
\[\widehat{\P}_x(A) = \int\limits_A f(\omega^1,\omega^2) \d (\P_x \times \P_x)(\omega^1,\omega^2)\]
for all $A \in \sigma(\F_t,\F^{T_0})$.  In particular the same equation is valid for all  $A$ in $\sigma(\F_t,\F^T)$ if $T \g T_0$.  This implies that whenever $T \g T_0$ the function $f_T$ coincides with the conditional expectation of $f$ to the $\sigma$-algebra $\sigma(\F_t,\F^T)$ with respect to $\P_x \times \P_x$.  Hence $f_T$ is a reverse martingale (all statements of this type are relative to the measure $\P_x\times \P_x$ from now on) when $T \to +\infty$ and converges almost surely to $f_\infty$ which is the Radon-Nikodym derivative of $\widehat{\P}_x$ with respect to $\P_x\times \P_x$ on $\sigma(\F_t,\F^\infty)$ (see \cite[pg. 483]{doob2001}).

It follows that $f_T\log(f_T)$ converges almost surely to $f_\infty\log(f_\infty)$ when $T$ goes to $+\infty$ and it remains to show only that these functions are uniformly integrable in order to obtain that
\[\lim\limits_{T \to +\infty} \int f_T \log(f_T ) \d(\P_x \times \P_x) = \int f_\infty \log(f_\infty) \d(\P_x \times \P_x)\]
and conclude (by the Gelfand-Yaglom-Perez Theorem as above) that
\[\lim\limits_{T \to +\infty} I_x(\F_t,\F^T) = I_x(\F_t,\F^\infty)\]
as claimed.

To simplify notation set $\varphi(x) = x\log(x)$ and $\G_T = \sigma(\F_t\times \F^T)$ (including the case $T = \infty$), and denote integration and conditional expectation with respect to $\P_x \times \P_x$ by $\E$.  We notice that $x \mapsto \varphi(x)$ is convex and always larger than or equal to $-e^{-1}$ on $x \ge 0$.

Setting $g = \varphi(f)$ and $g_T = \varphi(f_T)$ one has by Jensen's inequality
\[-e^{-1} \le g_T = \varphi(f_T) = \varphi\left(\E\left(f|\G_T\right)\right) \le \E\left(\varphi(f)|\G_T\right).\]

By the reverse martingale convergence theorem (see \cite[pg. 483]{doob2001}) the right hand side converges in $L^1$ to $\E\left(\varphi(f)|\G_\infty\right)$.  From this it follows that the functions $g_T$ are uniformly integrable which concludes the proof of claim 2.
\end{proof}

We now establish the result on the Radon-Nikodym derivative of $\widehat{\P}_x$ with respect to $\P_x \times \P_x$ which was used in the previous proof (see also \cite[pg. 354]{varopoulos1986}).
\begin{lemma}\label{radonderivative}
Let $M$ be a complete connected and stochastically complete Riemannian manifold.  Then for all $x$ and $0 \l t \l T \l +\infty$ the measure $\widehat{\P}_x$ restricted to $\sigma\left(\F_t\times \F^T\right)$ is absolutely continuous with respect to $\P_x \times \P_x$ restricted to the same $\sigma$-algebra and the corresponding Radon-Nikodym derivative is given by
\[f_T(\omega^1,\omega^2) = \frac{q(T-t,\omega^1_t,\omega^2_T)}{q(T,\omega^1_0,\omega^2_T)}.\]
\end{lemma}
\begin{proof}
Consider two subsets of $\Omega$ defined by
\[A = \lbrace\omega \in \Omega: \omega_{s_1} \in A_1,\ldots, \omega_{s_m} \in A_m\rbrace\]
\[B = \lbrace\omega \in \Omega: \omega_{t_1} \in B_1,\ldots, \omega_{t_n} \in B_n\rbrace\]
where $s_1 \l \cdots \l s_m = t$, $T = t_1 \l \cdots t_n$, and the sets $A_i$ and $B_j$ are Borel subsets of $M$.

By direct calculation using the definition of $f_T$ we obtain that
\[\int\limits_{A\times B} f_T(\omega^1,\omega^2)\d \P_x\times \P_x(\omega^1,\omega^2) = \int\limits_{A \times B} \frac{q(T-t,\omega^1_t,\omega^2_T)}{q(T,\omega^1_0,\omega^2_T)} \d \P_x(\omega^1) \d \P_x(\omega^2).\]

The right hand side coincides (via the definition of $\P_x$) with the integral over $A_1\times \cdots A_m \times B_1 \times \cdots \times B_n$ of
\[\frac{q(T-t,x_m,y_1)}{q(T,x,y_1)} q(s_1,x,x_1)\cdots q(s_m-s_{m-1},x_{m-1},x_m)q(t_1,x,y_1)\cdots q(t_n - t_{n-1},x_{n-1},x_n)\]
which after cancellation yields
\[\P_x(\omega_{s_1} \in A_1,\ldots,\omega_{s_m} \in A_m,\omega_{t_1} \in B_1,\ldots,\omega_{t_n} \in B_n).\]

This last probability is seen to be equal to $\widehat{P}_x(A \times B)$ by definition of $\widehat{\P}_x$.

Hence we have established that the integral of $f_T$ with respect to $\P_x\times \P_x$ over any set of the form $A \times B$ as above is $\widehat{\P}_x(A \times B)$.  Since any set in $\G = \sigma(\F_t \times \F^T)$ can be approximated by finite disjoint unions of such sets we have that the integral of $f_T$ on any set of this $\sigma$-algebra with respect to $\P_x\times \P_x$ is equal to the probability of the set with respect to $\widehat{\P}_x$.  As $f_T$ is $\G$-measurable this implies that $f_T$ is (a version of) the Radon-Nikodym derivative of $\P_x$ with respect to $\P_x \times \P_x$ on $\G$ as claimed.
\end{proof}

\subsection{Finite dimensional spaces of bounded harmonic functions}

Brownian motion on $\R^3$ is transient but satisfies the zero-one law.  Using this fact one may construct manifolds whose space of bounded harmonic functions is finite dimensional.

To see this, consider a Riemannian metric $g$ on $\R^3 \setminus \lbrace 0 \rbrace$ which coincides with the usual flat metric outside of the open ball $B_2$ of radius $2$ centered at $0$.  Let $f:\R^3 \setminus \lbrace 0\rbrace \to \R^3 \setminus \lbrace 0\rbrace$ be the inversion with respect to the sphere of radius $1$ centered at $0$ and suppose that $f$ is an isometry between the complement of $B_2$ and $B_{1/2}\setminus \lbrace 0 \rbrace$ endowed with the metric $g$.  In other words, endowed with this metric one can think of $\R^3 \setminus \lbrace 0 \rbrace$ as two copies of $\R^3$ (flat) which have been glued by removing an open ball from each and inserting tube connecting the two boundaries.

The resulting Riemannian manifold has a two dimensional space of bounded harmonic functions generated by the two functions 
\[x \mapsto \P_x(\lim\limits_{t \to +\infty}\omega_t = +\infty)\text{ and }x \mapsto \P_x(\lim\limits_{t \to +\infty}\omega_t = 0).\]
Notice that the sum of these two functions is constant.  A similar construction using more copies of $\R^3$ yields manifolds with spaces of bounded harmonic functions of dimension $3,4,5$, etc.

We will later on show that examples of this type cannot be `recurrent' in the sense that they can neither admit a compact quotient, nor be generic with respect to a harmonic measure on a compact foliation.  In short, a stationary random manifold almost surely either satisfies the Liouville property of has an infinite dimensional space of bounded harmonic functions (for leaves of foliations this was announced in \cite{kaimanovich1988}).   The result is a consequence of the following basic estimate (which is essentially \cite[Lemma 3.11]{gray2011}, in short the mutual information between two random variables is less than the entropy of either of them, so if one of them takes only finitely many, say $d$, values one gets an upper bound of $\log(d)$).

\begin{lemma}\label{finitedimensionallemma}
Let $M$ be a complete and stochastically complete Riemannian manifold and assume that the space of bounded harmonic functions on $M$ is finite dimensional of dimension $d$.  Then one has
\[I_x(\F_t,\F^\infty) \le \log(d)\]
for all $x \in M$ and all $t \g 0$.
\end{lemma}
\begin{proof}
By Lemma \ref{tailbijection} the space of bounded tail measurable functions modulo modification on $\P_x$-null sets has dimension $d$.  Hence there exists a partition $B_1,\ldots, B_d$ of $\Omega = C([0,+\infty),M)$ where each $B_j$ belongs to the tail $\sigma$-algebra $\F^\infty$ and is an atom for $\P_x$.

By definition $I_x(\F_t,\F^\infty)$ is the supremum over finite partitions $C_1,\ldots, C_n$ of $\Omega\times \Omega$ into $\F_t\times \F^\infty$-measurable sets of
\[\sum\limits_{j}\log\left(\frac{\widehat{\P}_x(C_j)}{\P_x \times \P_x(C_j)}\right)\widehat{P}_x(C_j).\]

However, by Dobrushin's theorem (see \cite[Lemma 7.3]{gray2011}) this coincides with the supremum over partitions consisting of sets of the form $A_i \times B_j$ where $A_i \in \F_t$ and $B_1,\ldots,B_d$ are as above.

For any such partition one obtains using Jensen's inequality and the fact that for any $p_1,\ldots,p_d \in [0,1]$ adding up to $1$ one has $-\sum p_j\log(p_j) \le \log(d)$ and $\sum p_j^2 \le 1$ the following chain of inequalities:
\begin{align*}
\sum\limits_{i,j}\log&\left(\frac{\widehat{\P}_x(A_i \times B_j)}{\P_x \times \P_x(A_i \times B_j)}\right)\widehat{P}_x(A_i \times B_j) =\sum\limits_{i,j}\log\left(\frac{\P_x(A_i \cap B_j)}{\P_x(A_i)\P_x(B_j)}\right)\P_x(A_i\cap B_j)
\\ &= -\sum\limits_{j}\log\left(\P_x(B_j)\right)\P_x(B_j) + \sum\limits_{i,j}\log\left(\frac{\P_x(A_i \cap B_j)}{\P_x(A_i)}\right)\frac{\P_x(A_i\cap B_j)}{\P_x(A_i)}\P_x(A_i)
\\ &\le \log(d) + \sum\limits_{i}\log\left(\sum\limits_{j = 1}^d\frac{\P_x(A_i \cap B_j)^2}{\P_x(A_i)^2}\right)\P_x(A_i) \le \log(d).
\end{align*}

From this the claim follows by taking supremum.
\end{proof}

\section{Entropy of stationary random manifolds}

In this section, after some preliminary work on Gromov-Hausdorff and smooth convergence of manifolds (see Theorems  \ref{gromovspacetheorem} and \ref{uniformlyboundedgeometry}),  we introduce the notion of a stationary random manifold (see Sections \ref{harmonicsection} and \ref{stationarysection}) and develop the basic entropy theory for it.

Examples of stationary random manifolds include: manifolds with transitive isometry group, manifolds admitting a compact quotient, and generic leaves of foliations (see Theorem \ref{harmoniccorrespondence}).  Further examples can be obtained by taking weak limits as discussed in Section \ref{stationarysection}.

Our main results are the existence of Kaimanovich entropy $h(M)$, equivalence of $h(M) = 0$ to the almost sure Liouville property, equivalence of $h(M) \g 0$ to the almost sure existence of an infinite dimensional space of bounded harmonic functions on $M$ (see Theorem \ref{entropytheorem}), and the inequalities relating linear drift of Brownian motion and volume growth to entropy (and hence to the Liouville property, see Theorem \ref{inequalities}).  We give several applications of these results at the end of the section.

\subsection{The Gromov space and harmonic measures}

\subsubsection{The Gromov space}\label{chaptertwowheregromovspace}

In this subsection we construct a model of `the Gromov space' which is a complete separable metric space whose points represent the isometry classes of all proper (i.e. closed balls are compact) pointed metric spaces.  The topology on the Gromov space is that of pointed Gromov-Hausdorff convergence (see \cite[Chapter 8]{burago-burago-ivanov2001}).  

Our main point is that one can construct the Gromov space using well defined sets (i.e. avoiding use of `the set of all metric spaces') and without using the axiom of choice (see \cite[Remark 7.2.5]{burago-burago-ivanov2001} and the paragraph preceding it).  We will later be interested in certain probability measures on the Gromov space.

A sequence of pointed metric spaces $(X_n,o_n)$ (here $o_n$ is the basepoint of the space which we will sometimes abuse notation by omitting; also, we use $d$ to denote the distance on different metric spaces simultaneously) is said to converge in the pointed Gromov-Hausdorff sense to a pointed metric space $(X,o)$ if for each $r \g 0$ and $\epsilon \g 0$ there exists $n_0$ and for all $n \g n_0$ a function $f_n:B_r(o_n) \to X$ (we use $B_r(x)$ to denote the ball of radius $r$ centered at $x$ in a metric space) satisfying the following three properties:
\begin{enumerate}
 \item $f_n(o_n) = o$
 \item $\sup\left\lbrace |d(f_n(x),f_n(y)) - d(x,y)|: x,y \in B_r(o_n)\right\rbrace \l \epsilon$
 \item $B_{r-\epsilon}(o) \subset \bigcup\limits_{x \in B_r(o_n)} B_\epsilon(f_n(x))$.
\end{enumerate}

Given two metric spaces $X$ and $Y$ we say a distance on the disjoint union $X \sqcup Y$ is admissible if it coincides with the given distance on $X$ when restricted to $X \times X$ and similarly for $Y$.

Following Gromov (see \cite[Section 6]{gromov1981}) we metricize pointed Gromov-Hausdorff convergence by defining the distance $\gromovdistance(X_1,X_2)$ between two pointed metric spaces $(X_1,o_1)$ and $(X_2,o_2)$ as the infimum of all $\epsilon \in (0,\frac{1}{2})$ such that there exists an admissible distance $d$ on the disjoint union $X_1 \sqcup X_2$ which satisfies the three inequalities $d(o_1,o_2) \l \epsilon$,$\ d(B_{1/\epsilon}(o_1),X_2) \l \epsilon$, and $d(X_1,B_{1/\epsilon}(o_2)) \l \epsilon$; or $\frac{1}{2}$ if no such admissible distance exists (this truncation is necessary in order for $\gromovdistance$ to satisfy the triangle inequality as noted by Gromov in the above-mentioned reference).  For a proof of the following lemma see for example \cite{cristina2008}.

\begin{lemma}\label{cristinaslemma}
The distance $\gromovdistance$ metricizes pointed Gromov-Hausdorff convergence.
\end{lemma}

We now consider the set of finite pointed metric spaces of the form $(X,o)$ where $X = \lbrace 0,\ldots,n\rbrace$ for some non-negative integer $n$ and $o = 0$.  Consider two such pointed metric spaces to be equivalent if they are isometric via a basepoint preserving isometry (each equivalence class has finitely many elements) and let $\finitegromovspace$ (read `finite Gromov space') be the set of all equivalence classes.  One verifies that $(\finitegromovspace,\gromovdistance)$ is a separable metric space.

\begin{definition}
We define the Gromov space $(\gromovspace,\gromovdistance)$ as the metric completion of $(\finitegromovspace,\gromovdistance)$.
\end{definition}

With the above definition it follows immediately that the Gromov space is a complete separable metric space with $\finitegromovspace$ as  a dense subset.  It remains to show that each of its points `represents' an isometry class of proper pointed metric spaces and that all such classes are represented by some point.

\begin{theorem}\label{gromovspacetheorem}
For each point $p$ in $\gromovspace$ there exists a unique (up to pointed isometry) proper pointed metric space $(X,o)$ with the property that all sequences $(X_n,o_n)$ of representatives of Cauchy sequences in $\finitegromovspace$ converging to $p$ converge in the pointed Gromov-Hausdorff sense to $(X,o)$.

Furthermore, the thus defined correspondence between isometry classes of pointed proper metric spaces and points in the Gromov space is bijective.
\end{theorem}

\begin{proof}
Consider a sequence $(X_n,o_n)$ of finite metric spaces representing some Cauchy sequence in $\finitegromovspace$.   By taking a subsequence we assume that the distance between $X_n$ and $X_{n+1}$ is less than $2^{-n}$ for all $n$.

By definition there exists an adapted metric $d_n$ on $X_n \sqcup X_{n+1}$ with the property that $d_n(o_n,o_{n+1}) \l 2^{-n}$ and the ball of radius $2^n$ centered at the basepoint of either half is at distance less than $2^{-n}$ from the other half.

Let $Y$ be the countable disjoint union of all $X_n$.   We define a distance on $Y$ by letting $d(x,y)$ in the case $x \in X_n$ and $y \in X_{n+k}$ be the infimum of
\[d_n(x_0,x_1)+\cdots + d_{n+k-1}(x_{k-1},x_k)\]
over all sequences of $k$ elements with $x = x_0,\ldots, x_k = y$ and $x_i \in X_{n+i}$ for all $i$.  The other case is determined by symmetry.

Set $X = \widehat{Y} \setminus Y$ and $o = \lim_{n \to +\infty} o_n$ where $\widehat{Y}$ is the metric completion of $Y$.  We claim $(X,o)$ is proper and is the limit of $(X_n,o_n)$ in the pointed Gromov-Hausdorff sense.  Once this claim is established uniqueness of $(X,o)$ up to pointed isometry is given by \cite[Theorem 8.1.7]{burago-burago-ivanov2001}.  And, since pointed Gromov-Hausdorff convergence is characterized by $\gromovdistance$ (see Lemma \ref{cristinaslemma}), the triangle inequality implies that $(X,o)$ is also the limit of any Cauchy sequence equivalent to the one determined by $(X_n,o_n)$.

We will now establish the claim.

Fix $r \g 0$ and let $B$ be the closed ball of radius $r$ centered at $o$ in $X$.  We must show that $B$ is compact.

For this purpose notice that for all $n$ and all $k \ge 0$ one has that ball of radius $2^{n-1}$ centered at $o_{n+k}$ is at distance less than $2^{-(n-1)}$ from $X_n$.  If $2^{n-1} \g r$ then one can approximate any $x$ by a sequence $x_k$ in $Y$ with the property that eventually $d(x_k,o_{n_k}) \l 2^{n-1}$ (where one chooses $n_k$ so that $x_k \in X_{n_k}$).  It follows that $n_k \to +\infty$ (otherwise infinitely many $x_k$ would belong to the same $X_N$ which is finite, and ultimately one obtains that $x \in X_N$) and one obtains that the distance between $B$ and $X_n$ is less than or equal to $2^{-(n-1)}$ as well.  In particular since $X_n$ is finite this shows that $B$ can be covered by a finite number of balls of radius $2^{-(n-2)}$.  This establishes that $B$ is compact as claimed.

We have shown that each equivalence class of Cauchy sequences in $Y$ determines a unique isometry class of pointed proper metric spaces.  Now let $(X,o)$ be  a pointed proper metric space and for each $n$ let $X_n = \lbrace o_n = x_{n,0},x_{n,1},\ldots, x_{n,k_n}\rbrace$ be a finite subset of $B_{2^n}(o)$ which is $2^{-n}$-dense.  There is a unique point $p_n \in \finitegromovspace$ such that all of its representatives are isometric to the pointed metric space  $(X_n,o_n)$.  Since $(X_n,o_n)$ converges in the pointed Gromov-Hausdorff sense to $(X,o)$ it follows that any sequence of such representatives converges to $(X,o)$ as well.  From this one obtains that $p_n$ converges to some point $p$ in $\gromovspace$ which represents the isometry class of $(X,o)$.  Hence the correspondence between points in $\gromovspace$ and isometry classes of pointed proper metric spaces is bijective, which concludes the proof.
\end{proof}

In view of the above theorem we will no longer distinguish between a point in $\gromovspace$ and a pointed proper metric spaces $(X,o)$ in the isometry class represented by it.
\subsubsection{Spaces of manifolds with uniformly bounded geometry}

We say a complete connected Riemannian manifold $M$ has geometry bounded by $(r,\lbrace C_k\rbrace)$ (where $r$ is a positive radius and $C_k$ is a sequence of positive constants indexed on $k = 0,1,\ldots$) if its injectivity radius is larger than or equal to $r$ and one has
\[|\nabla^k R| \le C_k\]
at all points, where $R$ is the curvature tensor, $\nabla^k R$ is its $k$-th covariant derivate, and the tensor norm induced by the Riemannian metric is used on the left hand side.

We denote by $\manifolds\left(d,r,\lbrace c_k\rbrace\right)$ the subset of the Gromov space representing isometry classes of $d$-dimensional pointed Riemannian manifolds with geometry bounded by $(r,\lbrace c_k\rbrace)$.

Following \cite[10.3.2]{petersen2006} we say a sequence $(M_n,o_n,g_n)$ of pointed connected complete Riemannian manifolds (here $g_n$ is the Riemannian metric an $o_n$ the basepoint) converges smoothly to a pointed connected complete Riemannian manifold $(M,o,g)$ if for each $r \g 0$ there exists an open set $\Omega \supset B_r(o)$ and for $n$ large enough a smooth pointed (i.e. $f_n(o) = o_n$) embedding $f_n:\Omega \to M$ such that the pullback metric $f_n^*g_n$ converges smoothly to $g$ on compact subsets of $\Omega$.

The following result is a consequence of $C^k$ compactness theorems from Riemannian geometry (see for example \cite[Theorem 72]{petersen2006}), a proof is provided in full detail in \cite{lessa2013}.

\begin{theorem}\label{uniformlyboundedgeometry}
Let $\manifolds = \manifolds\left(d,r,\lbrace C_k\rbrace\right)$ for some choice of dimension $d$, radius $r$, and sequence $C_k$.  Then $\manifolds$ is a compact subset of the Gromov space on which smooth convergence is equivalent to convergence in the pointed Gromov-Hausdorff sense.
\end{theorem}

We will say a subset of the Gromov space `consists of manifolds with uniformly bounded geometry' if it is contained in some set of the form $\manifolds\left(d,r,\lbrace C_k\rbrace\right)$.  Elements of such subsets are represented by triplets $(M,o_M,g_M)$ ($o_M$ being the basepoint and $g_M$ the Riemannian metric).  We usually write just $M$ leaving the other two elements implicit and will refer to them as $o_M$ and $g_M$ when needed.

Recall that a complete Riemannian manifold $M$ is said to be stochastically complete if the integral of its heat kernel $p(t,x,y)$ with respect to $y$ equals $1$ for all $t \g 0$ and $x \in M$.  We denote by $q(t,x,y) = p(t/2,x,y)$ the transition probability density of Brownian motion on such a manifold $M$.   With this convention one has that $q(t,x,y)$ is $(2\pi t)^{-1/2}e^{-(x-y)^2/2t}$ on $\R$ and $(2\pi t)^{-3/2}e^{-t/2-d(x,y)^2/2t}d(x,y)/\sinh(d(x,y))$ on three dimensional hyperbolic space (see \cite[pg. 185]{davies-mandouvalos1988}).

We will need the following uniform upper bound on the transition density $q$ for manifolds with uniformly bounded geometry.
\begin{theorem}\label{heatkernelupperbound}
Let $\manifolds$ be a subset of the Gromov space consisting of $n$-dimensional manifolds with uniformly bounded geometry.  Then for each $t_0 \g 0$ and $D \g 2$ there exist a positive constant $C$ such that the inequality
\[q(t,x,y) \le C\exp\left(-\frac{d(x,y)^2}{Dt}\right)\]
hold for all $t \ge t_0$ and all pairs of points $x,y$ belonging to any manifold $M$ of $\manifolds$.
\end{theorem}
\begin{proof}
The on diagonal bound given by Theorem 8 of \cite[pg. 198]{chavel1984} (setting $r = \sqrt{t}$) yields a constant $c_1$ depending only on $n$ such that any complete manifold of dimension $n$ satisfies
\[q(t,x,x) \le \frac{c_1}{\vol(B_{\sqrt{t/2}}(x))}t^{-\frac{n}{2}}.\]

Because of the uniform bounds on curvature and injectivity radius one may bound the volume of the ball of radius $\sqrt{t/2}$ from below uniformly on $\manifolds$ by some multiple of $t^{n/2}$ for small $t$ and by a constant for large $t$.

Using this one obtains that
\[q(t,x,x) \le \frac{1}{\gamma(t)}\]
for all $t \g 0$ and all $x$ in a manifold of $\manifolds$ where $\gamma(t)$ is of the form
\[\gamma(t) = \max(c_2 t^n, c_3).\]

One verifies that there exists $c_4$ such that $\gamma(t) \le \gamma(2t) \le c_4\gamma(t)$ for all $t \g 0$ after which by Corollary 16.4 of \cite{grigoryan2009} one obtains that for each $D \g 2$ there exist $c_5,c_6 \g 0$ (depending on $\manifolds$) such that
\[q(t,x,y) \le \frac{c_5}{\gamma(c_6 t)}\exp\left( -\frac{d(x,y)^2}{Dt}\right)\]
for all $x,y \in M$, $t \g 0$ and $M \in \manifolds$.

Restricting to $t \ge t_0$ on obtains for each $D \g 2$ a constant $C$ such that
\[q(t,x,y) \le C\exp\left( -\frac{d(x,y)^2}{Dt}\right)\]
for all $x,y \in M$, $t \ge t_0$ and $M \in \manifolds$, as claimed.
\end{proof}
\subsubsection{Harmonic measures}\label{harmonicsection}

We say a probability measure $\mu$ on the Gromov space is harmonic if it gives full measure to some set $\manifolds$ of manifolds with uniformly bounded geometry and is invariant under re-rooting by Brownian motion, i.e. one has
\[\int f(M,o,g)\d\mu(M,o,g) = \int f(M,x,g)q(t,o,x)\d x \d\mu(M,o,g)\]
for all $t \g 0$ and all bounded measurable functions $f:\manifolds \to \R$.

In order for the above equation to make sense one needs to know that the inner integral on the right hand side is Borel measurable on the Gromov space.  We prove this in the following lemma together with further regularity properties which will be useful to construct harmonic measures.  The key point is that the heat kernel depends continuously on the manifold in the smooth topology.  This intuitive fact was used (for time-dependent metrics) by Perelman in his proof of the geometrization conjecture after which it has received careful treatment by several authors (see \cite{lu2012} and the references therein).  It had also been previously used by Lucy Garnett to prove the existence of harmonic measures on foliated spaces which we will consider in the next subsection (see \cite[Fact 1]{garnett1983} and \cite{candel2003}).

\begin{lemma}\label{difussionlemma}
Let $\manifolds$ be a compact subset of the Gromov space consisting of manifolds of uniformly bounded geometry and for each $t \g 0$,$r \g 0$ and each function $f:\manifolds \to \R$ define $P^tf$ and $P^t_rf$ on $\manifolds$ by
\[P^t_rf(M,o,g) = \int\limits_{B_r(o)} f(M,x,g)q(t,o,x)\d x\]
\[P^tf(M,o,g) = \int\limits_{M} f(M,x,g)q(t,o,x)\d x.\]
Then the following properties hold:
\begin{enumerate}
 \item If $f$ is continuous then $P^t_rf$ is continuous.
 \item If $f$ is continuous then $P^t_rf$ converges uniformly to $P^tf$ when $r \to +\infty$. In particular, $P^tf$ is continuous.
 \item If $f$ is bounded and Borel measurable then $P^tf$ is also Borel measurable and satisfies $\sup |P^tf| \le \sup |f|$.
\end{enumerate}
\end{lemma}
\begin{proof}
We begin by showing that if $f$ is continuous then $P^t_rf$ is as well.

Consider a manifold $(M,o,g) \in \manifolds$ and a sequence $(M_n,o_n,g_n) \in \manifolds$ converging to it.   By Theorem \ref{uniformlyboundedgeometry} the convergence is smooth so there exists an exhaustion $U_n$ of $M$ by precompact open sets and a sequence of diffeomorphisms $\varphi_n:U_n \to V_n \subset M_n$ such that $\varphi_n(o) = o_n$ and the pullback metric $\varphi_n^*g_n$ converges smoothly on compact subsets to $g$.

In this situation Theorem 2.1 of \cite{lu2012} (applied in the case where the fields $X_n$ and the potentials $Q_n$ are equal to $0$ and the metrics $g_n(\tau)$ are constant with respect to $\tau$) guarantees that the sequence of pullbacks $q_n(t,o,x) = q^{M_n}(t,o_n,\varphi(x))$ of the transition probability densities of Brownian motion on each $M_n$ converges uniformly on compact subsets of $[0,+\infty) \times M$ to a fundamental solution $\tilde{q}(t,o,\cdot)$ of the heat equation (that fact that we use $\Delta/2$ instead of $\Delta$ is clearly inessential) which satisfies
\[\int\limits_{M} \tilde{q}(t,o,x)\d x \le \sup\left\lbrace \int_{M_n}q^{M_n}(t,o_n,x)\d x \right\rbrace = 1.\]

By Theorem 4.1.5 \cite{hsu2002} the transition density $q(t,o,x)$ is the minimal fundamental solution so one has $q(t,o,x) \le \tilde{q}(t,o,x)$.  Combined with the fact that the integral of both kernels with respect to $x$ is at most $1$ and the $\int q(t,o,x) \d x = 1$ one obtains $q(t,o,x) = \tilde{q}(t,o,x)$ so that $q_n(t,o,\cdot)$ converges uniformly on compact sets to $q(t,o,\cdot)$.

Setting $F(x) = f(M,x,g)$ and $F_n(x) = f(M_n,\varphi_n(x),g_n)$ and using the fact that $f$ is uniformly continuous (because $\manifolds$ is compact) one obtains that $F_n \to F$ uniformly on compact subsets of $M$.

Finally because the pullback metrics $\varphi_n^*g_n$ converge smoothly to $g$, the Jacobian $J_n$ of $\varphi_n$ converges uniformly to $1$ on compact subsets and also the open sets $\Omega_n = \varphi_n^{-1}(B_r(o_n))$ converge in the Hausdorff distance to $B_r(o)$.

Combining these four facts (local uniform convergence of $q_n(t,o,\cdot)$ to $q(t,o,\cdot)$, $F_n$ to $F$, $J_n$ to $1$, and Hausdorff convergence of $\Omega_n$ to $B_r(o)$) with the fact that $f$ is bounded one obtains
\begin{align*}
\lim\limits_{n \to +\infty}P^t_rf(M_n,o_n,g_n) &= \lim\limits_{n \to +\infty}\int\limits_{B_r(o_n)}q(t,o_n,x)f(M_n,x,g_n)\d x 
\\ &= \lim\limits_{n \to +\infty}\int\limits_{\Omega_n}q_n(t,o,x)F_n(x)J_n(x)\d x
\\ &= \lim\limits_{n \to +\infty}\int\limits_{B_r(o)}q_n(t,o,x)F_n(x)J_n(x)\d x
\\ &= \lim\limits_{n \to +\infty}\int\limits_{B_r(o)}q(t,o,x)F(x)\d x
\\ &= \lim\limits_{n \to +\infty}\int\limits_{B_r(o)}q(t,o,x)f(M,x,g)\d x
\\ &= P^t_rf(M,o,g)
\end{align*}
which implies that $P^t_rf$ is continuous as claimed.

We will now show that $P^t_rf$ converges uniformly to $P^tf$ when $r \to +\infty$ if $f$ is continuous.

By Theorem \ref{heatkernelupperbound} for each fixed $t \g 0$ there exists a constant $C$ such that
\[q(t,o,x) \le C\exp\left(-d(o,x)^2/C\right)\]
for all manifolds $(M,o,g) \in \manifolds$ and all $x \in M$.  Furthermore by the Bishop comparison theorem one may increase $C$ above so that the volume of the ball of radius $r$ is bounded from above by $\exp(Cr)$.  

Combining these two facts one obtains that
\begin{align*}|(P^t_r-P^t)f(M,o,g)| &\le \int\limits_{d(o,x) \g r}|f(M,x,g)|q(t,o,x)\d x
\\ &\le C\sup|f|\int\limits_{d(o,x) \g r}\exp\left(-d(o,x)^2/C\right) \d x
\\ &\le C\sup|f| \sum\limits_{n = 1}^{+\infty} \exp(-(nr)^2/C + C(n+1)r) 
\\ &\le C\sup|f| \sum\limits_{n = 1}^{+\infty} \exp(-(nr)^2/C + C(n+1)r)
\\ &= C\sup|f|\varphi(r)
\end{align*}
where the last inequality is obtained by bounding the integral by the sum over anulii of the form $B_{(n+1)r}(o) \setminus B_{nr}(o)$ and the integral over each anulus by the maximum value time the volume of the ball of radius $(n+1)r$.

As soon as $r \g 2C^2$ one has that $C(n+1)r - n^2 r^2/C$ is decreasing with respect to $r$ for all $n \ge 1$.  Hence $\varphi(r) \to 0$ when $r \to +\infty$ and one obtains that $P^t_rf$ converges uniformly to $P^tf$ as claimed.

To conclude we will prove that $P^tf$ is Borel for all bounded Borel $f$ and that it is bounded in absolute value by $\sup|f|$.

For this purpose consider for some $C \g 0$ the family of functions $\F$ on $\manifolds$ bounded in absolute value by $C$ and such that $P^tf$ is Borel measurable.  We have shown that $\F$ contains the continuous functions bounded in absolute value by $C$.  By the dominated convergence theorem it is closed by pointwise limits (this is because Borel functions are the smallest class containing continuous functions and closed under pointwise limits, see for example \cite[Theorem 11.6]{kechris1995}).  Therefore it contains all Borel measurable functions bounded by $C$ in absolute value.  Since this works for all $C$ one has that $P^tf$ is Borel measurable for all bounded Borel measurable $f$.  The claim $\sup |P^t f| \le \sup |f|$ follows directly from the definition of $P^tf$ because one has $\int q(t,x,y)\d y = 1$ on all manifolds in $\manifolds$.
\end{proof}

The following theorem implies that one can associate at least one harmonic measure to each manifold of bounded geometry.  By this we mean that if $(M,g)$ has bounded geometry then the closure in the Gromov space of the set of pointed manifolds of the form $(M,x,g)$ supports at least one harmonic measure.  We will see other examples of harmonic measures in the next subsection.

\begin{theorem}\label{harmonicexistence}
If $\manifolds$ is a compact subset of the Gromov space consisting of manifolds with uniformly bounded geometry then there exists at least one harmonic measure supported on $\manifolds$.
\end{theorem}
\begin{proof}
For each probability $\mu$ on $\manifolds$ and $t \g 0$ define the measure $P^t\mu$ on $\manifolds$ using the Riesz representation theorem and Lemma \ref{difussionlemma} in such  way that for all continuous $f:\manifolds \to \R$ one has
\[\int_\manifolds f \d P^t\mu = \int\limits_\manifolds P^t f \d \mu.\]

The maps $\lbrace P^t: t \ge 0\rbrace$ form a commuting family of linear maps which leave the convex and weakly compact set of probability measures on $\manifolds$ invariant.   By the Markov-Kakutani fixed point theorem there is a common fixed point for all the $P^t$ which must be a harmonic measure.
\end{proof}

\subsubsection{Foliations and leaf functions}

Harmonic measures on foliations were introduced by Lucy Garnett in \cite{garnett1981} (see also \cite{garnett1983} and \cite{candel2003}).  In this subsection we will explore how they relate to harmonic measures on the Gromov space in the sense of our definition.

To begin we must fix a definition of foliation.  There are several definitions in the literature, the crucial feature for our purposes is that each leaf should be a Riemannian manifold.  An important example is given by the foliation defined by an integrable distribution of tangent subspaces on a Riemannian manifold, in this case each leaf inherits a Riemannian metric from the ambient space.

A $d$-dimensional foliation is a metric space $X$ partitioned into disjoint subsets called leaves.  Each leaf is a continuously and injectively immersed $d$-dimensional connected complete Riemannian manifold.  Furthermore, for each $x \in X$ there is an open neighborhood $U$, a Polish space $T$, and a homeomorphism $h:\R^d \times T \to U$ with the following properties:
\begin{enumerate}
 \item For each $t \in T$ the map $x \mapsto h(x,t)$ is a smooth injective immersion of $\R^d$ into a single leaf.
 \item For each $t\in T$ let $g_t$ be the metric on $\R^d$ obtained by pullback under $x \mapsto h(x,t)$ of the corresponding leaf's metric.  If a sequence $t_n$ converges to $t \in T$ then the Riemannian metrics $g_{t_n}$ converge smoothly on compact sets to $g_t$.
\end{enumerate}

As part of a program to study the geometry of topologically generic leaves Álvarez and Candel introduced the `leaf function' which is a natural function into the Gromov space associated to each foliation $X$ (see \cite{alvarez-candel2003}).   It is defined as the function  mapping each point $x$ in the foliation to the leaf $L(x)$ containing it considered as a pointed Riemannian manifold with basepoint $x$.

We will now establish that the leaf function is Borel measurable.  We do this by using a result of Solovay for which we must assume the existence of an inaccessible cardinal.  The author believes a more direct proof without this assumption is attainable in the same vein as \cite{lessa2013} where semicontinuity of the leaf function is established and related to Reeb type stability results.

\begin{lemma}\label{leaffunctionregularity}
Let $X$ be a compact foliation and $L$ its leaf function.  Then the following holds:
\begin{enumerate}
 \item $L$ takes values in a compact subset of the Gromov space consisting of manifolds with uniformly bounded geometry.
 \item $L$ is measurable with respect to the completion of the Borel $\sigma$-algebra with respect to any probability measure on $X$. 
\end{enumerate}
\end{lemma}
\begin{proof}
The first claim follows from Theorem \ref{uniformlyboundedgeometry} (see \cite{lessa2013} for details).

To establish the second claim suppose $X$ is a compact foliation, $\mu$ is a probability on $X$ and there exists an open set $U$ in the Gromov space such that $L^{-1}$ is not $\mu$-measurable (where $L$ is the leaf function of $X$).

By \cite[Théorème 4-3]{delarue1993} (see also \cite{rohlin1952}) there exists a full measure set $X' \subset X$ and a bi-measurable bijection $f:X' \to \R$ such that $f(X') = [0,m_0] \sqcup C$ where $m_0 \ge 0$ and $C$ is a countable subset of $\R$ disjoint form $[0,m_0]$, such that $f_*\mu$ equals the sum of Lebesgue measure on $[0,m_0]$ with a probability measure on $C$.

If follows that $f(L^{-1}(U))$ is not Lebesgue measurable.  And we have therefore constructed a non-Lebesgue measurable subset of $\R$ without using the axiom of choice.

Assuming the existence of an inaccessible cardinal this is not possible due to \cite[Theorem 1]{solovay1970}.
\end{proof}

A probability measure $m$ on a foliation $X$ is said to be harmonic (see \cite[Fact 4]{garnett1983}) if it satisfies
\[\int f(x)\d m(x) = \int q(t,x,y)f(y)\d y\d m(x)\]
for all bounded measurable functions $f:X \to \R$.

Every compact foliation admits at least one harmonic measure (see \cite{garnett1983} and \cite{candel2003}).  The following theorem implies that any result which establishes properties of generic manifolds for harmonic measures on the Gromov space immediately implies a similar result for generic leaves of compact foliations.

\begin{theorem}\label{harmoniccorrespondence}
Let $X$ be a compact foliation with leaf function $L$ and $m$ a harmonic measure on $X$.   Then the push-forward measure $L_*m$ is harmonic measure on the Gromov space.
\end{theorem}
\begin{proof}
Let $\manifolds$ be a compact set of the Gromov space containing the image of $L$ and consisting of manifolds with uniformly bounded geometry.  If $f:\manifolds \to \R$ is bounded and measurable then by definition of $L_*m$ and harmonicity of $m$ one has
\begin{align*}
\int\limits_\manifolds \int\limits_{L(x)}&q(t,x,y)f(M,y,g)\d y \d L_*m(M,x,g)
\\& = \int\limits_X \int\limits_{L(x)} q(t,x,y)f\circ L(y)\d y \d m(x)
\\& = \int\limits_X f\circ L(x) \d m(x) = \int\limits_\manifolds f(M,x,g)\d L_*m(M,x,g)
\end{align*}
so $L_*m$ is harmonic as claimed.
\end{proof}

\subsection{Asymptotics of random manifolds}
\subsubsection{Stationary random manifolds}\label{stationarysection}

We define a stationary random manifold to be a random element of the Gromov space whose distribution is a harmonic measure.  A typical example is obtained as follows: let $X$ be a compact foliation and $m$ a harmonic measure on $X$, the leaf function $L:X \to \gromovspace$ is a stationary random manifold defined on $(X,m)$.

As noted in the previous section if a bounded geometry manifold $(M,g)$ admits a finite volume quotient under isometries and one takes a random point $o$ in a fundamental domain of this action distributed according to the normalized volume measure then $(M,o,g)$ is a stationary random manifold.

Another way of obtaining stationary random manifolds is by taking weak limits.  For example, for each $n$ let $(M_n,g_n)$ be a compact hyperbolic surface with genus $n$ whose injectivity radius is larger than a prescribed constant $\epsilon \g 0$.  If $o_n$ is uniformly distributed on $(M_n,g_n)$ then the $(M_n,o_n,g_n)$ is stationary.  Since the sequence takes values in a set of manifolds with uniformly bounded geometry there is a weak limit $(M,o,g)$ which is also stationary.  

In the context of graphs a similar construction (where the sequence consists of finite binary trees) yields as a limit the so-called Canopy tree, which has been shown for example to encode information about the asymptotic spectrum of Schrödinger operators on the sequence (see \cite{aizenman-warzel2006}).  Hence, one might expect a stationary random manifold constructed from a sequence of compact manifolds as above to encode information about the asymptotic behavior of the spectrum of the Laplacian on the sequence.  However, no such result is known to the author at the time of writing.
\subsubsection{Entropy}

We introduce an asymptotic quantity `Kaimanovich entropy' associated to each stationary random manifold which measures the asymptotic behavior of the differential entropy between the time $t$ distribution of Brownian motion and the Riemannian volume measure.   Several alternate definitions for this quantity in different contexts as well as theorems and applications where announced in the interesting papers \cite{kaimanovich1986} and \cite{kaimanovich1988}.  For the case of manifolds with a compact quotient some of these properties (e.g. the so-called Shannon-McMillan-Breiman type theorem) were later proved in \cite{ledrappier1996}.

The main point is that Kaimanovich entropy relates directly to the mutual information between $\sigma$-algebras $\F_t$ and $\F^T$ for Brownian motion.  This allows one to show that Kaimanovich entropy is zero if and only if the random manifold satisfies the Liouville property almost surely.  Later on we will will relate entropy to other asymptotic quantities.

The following technical lemma sets the basis for our study.
\begin{lemma}\label{finitetimeentropy}
Let $\manifolds$ be a compact subset of the Gromov space consisting of manifolds with uniformly bounded geometry.  The following functions are finite and continuous with respect to $T \g t \g 0$ and $M \in \manifolds$:
\begin{align*}
h_t(M) &= -\int\limits_M q(t,o_M,x)\log(q(t,o_M,x))\d x,
\\ I_t^T(M) &= \int\limits_{M \times M}\log\left(\frac{q(T-t,x,y)}{q(T,o_M,y)}\right)q(t,o_M,x)q(T-t,x,y)\d x\d y. 
\end{align*}
Also, the following formula holds (where $P^t$ is defined by Lemma \ref{difussionlemma}):
\[I_t^T(M) = h_T(M) - (P^th_{T-t})(M).\]
\end{lemma}
\begin{proof}
The proof is similar to that of Lemma \ref{difussionlemma}.  We define 
\[h_t^r(M) = -\int_{B_r(o_M)}q(t,o_M,x)\log(q(t,o_M,x))\d x\]
with the purpose of showing that $h_t^r$ is continuous and converges uniformly to $h_t$ on $\manifolds$ when $r \to +\infty$.

Assume $\manifolds$ consists of $n$-dimensional manifolds with curvature greater than $-k^2$ and let $K(t,r)$ be the heat kernel at time $t$ between two points at distance $r$ in the $n$ dimensional hyperbolic plane with constant curvature $-k^2$.  Then by Theorem 2.2 of \cite{ichihara1988} one has $q(t,o_M,x) \ge K(t,d(o_M,x))$ for all $M \in \manifolds$.  From \cite[Theorem 3.1]{davies-mandouvalos1988} one obtains that $\log(K(t,r))$ is bounded from below by a polynomial $f$ in $r$ on any compact interval of positive times.

On the other hand by Theorem \ref{heatkernelupperbound} for each compact interval of times there exists $C \g 0$ such that one has $q(t,o_M,x) \le Ce^{-d(o_M,x)^2/3t}$ on all $M$ of $\manifolds$.

Combining these facts yields
\[-\log(C) \le -\log(q(t,o_M,x)) \le f(d(o_M,x))\]
for all $M \in \manifolds$.

By Bishop's volume comparison theorem the volume of the ball of radius $r$ in any manifold of $\manifolds$ is bounded by that in $n$-dimensional hyperbolic space with curvature $-k^2$.  Hence one has an upper bound for volume of the form  $\vol(B_r(o_M)) \le \exp(Cr)$ (notice that the previous inequality remains valid if one increases $C$ so there is no problem in using the same constant for both bounds).

Similarly to the proof of Lemma \ref{difussionlemma} one obtains for each $r$ a positive constant $\epsilon(r)$ which decreases to $0$ as $r \to +\infty$ such that
\[|h_t(M) - h_t^r(M)| \le \int\limits_{M \setminus B_r(o_M)} q(t,o_N,x)|\log(q(t,o_M,x))|\d y \le \epsilon(r)\]
for all manifolds in $\manifolds$.

Hence $h_t(M)$ is the uniform limit on $\manifolds$ of $h_t^r(M)$ when $r \to +\infty$ and it suffices to establish continuity of the later.  

Before doing that we establish continuity with respect to $t$ of $h_t$.  Given $M \in \manifolds$, $\epsilon \g 0$ and $t \g 0$, one can find $r \g 0$ such that $|h_s - h_s^r| \l \epsilon/3$ for all $s$ in a compact neighborhood of $t$ (notice that our bounds were obtained uniformly on such intervals).  Since $q(t,o_M,x)$ is continuous with respect to $t$ and $x$ one has that $q(s,o_M,\cdot)$ converges uniformly to $q(t,o_M,\cdot)$ on $B_r(o_M)$ when $s \to t$.  Hence $h_s^r(M) \to h_t^r(M)$.  Combining the two facts one obtains that there exists a neighborhood of $t$ on which
\[|h_s(M) -h_t(M)| \le |h_s(M) - h_s^r(M)|+ |h_s^r(M) - h_t^r(M)| + |h_t^r(M) - h_t(M)| \l \epsilon\]
which yields continuity of $h_t$ with respect to $t$ as claimed.

We now establish continuity of $h_t^r(M)$ with respect to $M$. Assume the sequence $(M_n,o_n,g)$ in $\manifolds$ converges to $(M,o,g)$.  By Theorem \ref{uniformlyboundedgeometry} convergence is smooth so that there exists an exhaustion $U_n$ of $M$ by precompact open sets and a sequence of diffeomorphisms $\varphi_n:U_n \to V_n \subset M_n$ such that $\varphi_n(o) = o_n$ and the pullback metric $\varphi_n^*g_n$ converges smoothly on compact subsets to $g$.

From this it follows that the Jacobian $J(x)$ of $\varphi_n$ at $x$ converges to $1$ uniformly on compact sets.  And by the results of \cite{lu2012} the functions $q_n(t,o,x) = q^{M_n}(t,o_n,\varphi_n(x))$ ($q^{M_n}$ being the transition density of Brownian motion on $M_n$) converge uniformly on compact sets to the transition density $q(t,o,x)$ of Brownian motion on $M$.  Using this the continuity of $h_t^r(M)$ follows as claimed.

We will establish the formula $I_t^T = h_T - P^th_{T-t}$, from this the continuity of $I_t^T$ follows from that of $h_T$ and $h_{T-t}$ by Lemma \ref{difussionlemma}.  The proof is the following computation using the property $\int q(t,x,y)q(s,y,z)\d y = q(t+s,x,z)$ of the heat kernel:
\begin{align*}
I_t^T(M) &= \int \log\left(\frac{q(T-t,x,y)}{q(T,o,y)}\right)q(t,o,x)q(T-t,x,y)\d x\d y 
\\ &= \int\log\left(q(T-t,x,y)\right)q(t,o,x)q(T-t,x,y)\d x \d y
\\ & -\int \log\left(q(T,o,y)\right)q(t,o,x)q(T-t,x,y)\d x \d y
\\ &= - \int q(t,o,x)h_{T-t}(M,x,g)\d x - \int \log\left(q(T,o,y)\right)q(T,o,y)\d y
\\ &= - P^t h_{T-t}(M) + h_T(M).
\end{align*}
\end{proof}

The following result implies in particular that a stationary random manifold almost surely either has an infinite dimensional space of bounded harmonic functions or satisfies the Liouville property (for generic leaves of a foliation this was announced in \cite{kaimanovich1988}).

\begin{theorem}\label{entropytheorem}
The following limit (Kaimanovich entropy) exists and is non-negative for any ergodic stationary random manifold $M$:
\[h(M) = \lim\limits_{t \to +\infty}\E\left(\frac{1}{t}h_t(M)\right).\]
Furthermore, $h(M) = 0$ if and only if $M$ is almost surely Liouville and, $h(M) \g 0$ if and only if the space of bounded harmonic functions on $M$ is infinite dimensional almost surely.
\end{theorem}
\begin{proof}
Let $H_t = \E(h_t(M))$ and notice that by dominated convergence it is continuous with respect to $t \g 0$, bounded by the maximum of $h_t$ on $\manifolds$, and by Lemma \ref{finitetimeentropy} one has
\[H_T - H_{T-t} = \E(h_T(M) - h_{T-t}(M)) = \E(h_T(M) - P^t h_{T-t}(M)) = \E(I_t^T(M)).\]

The mutual information $I_t^T$ is non-negative and decreases to $I_t^{\infty}(M)$ when $T \to +\infty$ (see Theorem \ref{informationtheorem} and preceding paragraphs). By the monotone convergence theorem it follows that $T \mapsto H_T - H_{T-t}$ decreases to $\E(I_t^\infty(M)) \ge 0$ when $T \to +\infty$.  From this one obtains that 
\[h(M) = \lim_{T \to +\infty} \E(\frac{1}{T}h_T(M)) = \E(\frac{1}{t}I_t^{\infty}(M))\]
for all $t \g 0$.

If $N$ is a manifold with bounded geometry then by Lemma \ref{finitetimeentropy} one has that $I_t^T(N)$ is finite.  It follows from Theorem \ref{informationtheorem} that $I_t^{\infty}(N) = 0$ for some $t \g 0$ if and only if $N$ is Liouville.  Hence $h(M) = 0$ if and only if $M$ is almost surely Liouville as claimed.

On the other hand if $h(M) \g 0$ then, because $t \mapsto I_t^\infty(M)$ is increasing, one has by monotone convergence
\[\E(\lim\limits_{t \to +\infty}I_t^\infty(M)) = \lim\limits_{t \to +\infty}th(M) = +\infty.\]
This implies that $L(M) = \lim\limits_{t \to +\infty}I_t^{\infty}(M) = +\infty$ with positive probability.  Since this $L$ is measurable and $M$ is ergodic one must have $L(M) = +\infty$ almost surely.   By Lemma \ref{finitedimensionallemma} this implies that the dimension of the space of bounded harmonic functions on $M$ is infinite almost surely.
\end{proof}

\subsubsection{Linear Drift}

Given a pointed bounded geometry manifold $(M,o,g)$ we define by
\[\ell_t(M) = \int d(o,x)q(t,o,x)\d x\]
the mean displacement of Brownian motion starting at $o$ from its starting point.  We are interested in the asymptotics of $\ell_t(M)$ when $t \to +\infty$ and $M$ is a stationary random manifold.  In particular we introduce the linear drift of a random manifold $M$ via the following theorem.

\begin{theorem}\label{lineardrift}
Let $M$ be a stationary random manifold.  Then the linear drift
\[\ell(M) = \lim\limits_{t \to +\infty}\E\left(\frac{1}{t}\ell_t(M)\right)\]
exists and is finite.
\end{theorem}
\begin{proof}
Suppose that $M$ takes values in a set of manifolds $\manifolds$ with uniformly bounded geometry.  We begin by establishing that $\ell_t(M)$ is continuous with respect to both $M$ and $t$.

If $(M_n,o_n,g_n)$ is a sequence in $\manifolds$ converging to $(M,o,g)$ then by Theorem \ref{uniformlyboundedgeometry} there exists an exhaustion $U_n$ of $M$ by relatively compact open sets and smooth embeddings $\varphi_n:U_n \to M_n$ with $\varphi_n(o) = o_n$ such that $\varphi_n^*g_n$ converges smoothly $g$ on compact subsets of $M$.

It follows that the Jacobian $J_n$ of $\varphi_n$ converges uniformly to $1$ on compact sets and $d_n(x) = d(o_n,\varphi_n(x))$ converges uniformly on compact sets to $d(o,x)$.  Also, by Theorem 2.1 of \cite{lu2012}, one has that $q_n(t,x) = q(t,o_n,\varphi_n(x))$ converges uniformly on compact sets to $q(t,o,x)$.

Combining these facts one obtains that for each $r \g 0$
\[\ell_t^r(M) = \int\limits_{B_r(o)}q(t,o,x)d(o,x)\d x\]
depends continuously on $(M,o,g) \in \manifolds$.

Let $A,B \g 0$ be given by Lemma \ref{radialvariance} below.  By Jensen's inequality one has for all $M \in \manifolds$ that
\begin{align*}\left(\ell(M) - \ell_t^{At}(M)\right)^2 &= \left(\int\limits_{M \setminus B_{At}(o)} q(t,o,x)d(o,x)\d x\right)^2
\\ &\le \int\limits_{M \setminus B_{At}(o)} q(t,o,x)d(o,x)^2\d x \le Be^{-At}t^2
\end{align*}
which establishes that $\ell_t^r$ converges uniformly to $\ell_t$ on $\manifolds$ when $r \to +\infty$ for all $t \ge 1$.  In particular $\ell_t$ is continuous with respect to $M$ on $\manifolds$ for $t \ge 1$ (in fact this is true for all $t$ but we will not need it).

To establish continuity with respect to $t$ assume $t,s \g 1$ and notice that using Lemma \ref{radialvariance} as above one obtains that for all $T \g \max(s,t)$ the integrals of both $d(o,x)q(t,o,x)$ and $d(o,x)q(s,o,x)$ on $B_{AT}(o)$ are bounded by $Be^{-At}t^2$.  Combining this with the fact that if $s \to t$ then $q(s,o,\cdot)$ converges uniformly on $B_{AT}(o)$ to $q(t,o,\cdot)$ yields the desired result.

It follows that if $s \to t \g 1$ then $\ell_s(M)$ converges uniformly on $\manifolds$ to $\ell_t(M)$ and hence
\[L_t = \E(\ell_t(M))\]
is continuous with respect to $t \ge 1$.

We will now establish that $L_t$ is subadditive, i.e. satisfies $L_{t+s} \le L_t + L_s$, from which the existence of the finite limit $\lim L_t/t$ follows .
For this purpose we calculate using the triangle inequality
\begin{align*}\ell_{t+s}(M) &= \int q(t+s,o,x)d(o,x)\d x
\\ &= \int q(t,o,y)q(s,y,x)d(o,x)\d x\d y
\\ &\le \int q(t,o,y)q(s,y,x)(d(o,y)+d(y,x))\d x\d y = \ell_t(M) + P^t\ell_s(M) 
\end{align*}
which taking expectation yields the desired result.
\end{proof}

Define $\ell^+(M)$ for a stationary random manifold $M$ as the infimum of all $L \g 0$ such that
\[\lim\limits_{t \to +\infty}\E\left(\int\limits_{B_{Lt}(o)} q(t,o,x)\d x\right) = 1\]
we wish to guarantee that $\ell(M) = \ell^+(M)$.

However, a counterexample is given by a random manifold $M$ which is equal to the hyperbolic plane with constant curvature $-1$ or $-2$ each with probability $1/2$.  In this example $\ell(M)$ equals $1.5$ but $\ell^+(M)$ equals $2$.  The problem arises because the distribution of $M$ is a convex combination of other harmonic measures (in this case Dirac deltas).  We say a harmonic measure on a compact set $\manifolds$ of manifolds with uniformly bounded geometry is ergodic if it is extremal among all harmonic measures on this set.  A random manifold is said to be ergodic if its distribution is.

\begin{lemma}\label{lplusequalsl}
Let $M$ be a stationary random manifold. Then $\ell(M) \le \ell^+(M)$.  Furthermore, if $M$ is ergodic then $\ell(M) = \ell^+(M)$.
\end{lemma}
\begin{proof}
Let $A$ and $B$ be given by Lemma \ref{radialvariance}.  For any $L \g 0$ one has
\[\frac{1}{t}\ell_t(M) \le L + A\int\limits_{M \setminus B_{Lt}(o)}q(t,o,x)\d x + \int\limits_{M \setminus B_{At}(o)}\frac{d(o,x)}{t}q(t,o,x)\d x.\]

Using Jensen's inequality and Lemma \ref{radialvariance} one obtains that the third term is bounded from above by $\sqrt{B}e^{-At/2}$.  This implies in particular that $\ell^+(M) \le A$.  If $L \g \ell^+(M)$ then the expectation of the second term goes to zero from which one obtains that $\ell(M) \le \ell^+(M)$.

Proof of the converse inequality will be postponed until the next section (see Corollary \ref{citethisfromchapter2}).
\end{proof}

The main result of \cite{ichihara1988} is that one can compare the radial process of Brownian motion on a Riemannian manifold with that of a model space with constant curvature.  Combining this result with upper bounds for the heat kernel yields the following technical lemma which we have used in the proofs above.
\begin{lemma}\label{radialvariance}
Let $\manifolds$ be a compact subset of the Gromov space consisting manifolds with uniformly bounded geometry.  Then there exist constants $A,B \g 0$ such that for each Brownian motions $X_t$ starting at the origin $o$ of a manifold $(M,o,g)$ in $\manifolds$ one has
\[\E\left( \frac{d(o,X_t)^2}{t^2} 1_{\left\lbrace d(o,X_t) \g At\right\rbrace}\right) \le Be^{-At}\]
for all $t \ge 1$.
\end{lemma}
\begin{proof}
Let $n$ denote the dimension of the manifolds in $\manifolds$ and let $-k^2$ be a lower bound for their curvature. 

Letting $Y_t$ be a Brownian motion starting at the origin of $\R^n$ endowed with a complete metric of constant curvature $-k^2$ one has by \cite[Theorem 2.1]{ichihara1988} that if $X_t$ is any Brownian motion starting at the origin of a manifold in $\manifolds$ then
\[\P(d(Y_0,Y_t) \g x) \ge \P(d(X_0,X_t) \g x)\]
for all $x$.

Notice that if $V,W$ are non-negative random variables with $\P(V \g x) \ge \P(W \g x)$ for all $x$ then $\E(V) \ge \E(W)$.  In particular for any non-decreasing non-negative function $f$ one has $\E(f(V)) \ge \E(f(W))$.

Applying this observation one obtains
\[\E\left(\frac{d(Y_0,Y_t)^2}{t^2}1_{\left\lbrace d(Y_0,Y_t) \g At\right\rbrace}\right) \ge \E\left(\frac{d(X_0,X_t)^2}{t^2}1_{\left\lbrace d(X_0,X_t) \g At\right\rbrace}\right)\]
for all $A \g 0$ and all $t \ge 0$.  So that it suffices to bound the expectation on the left hand side.

Letting $K(t,r)$ denote the probability transition density Brownian motion on the hyperbolic plane with constant curvature $-k^2$.  One has explicitly
\[\E\left(d(Y_0,Y_t)1_{\left\lbrace d(Y_0,Y_t) \g At\right\rbrace}\right) = \int\limits_{At}^{+\infty}\frac{r^2}{t^2}\vol(S^{n-1})\frac{1}{k}\sinh(kr)^{n-1}K(t,r)\d r,\]
where $S^{n-1}$ is the standard $n-1$-dimensional sphere

By Theorem \ref{heatkernelupperbound} there exists a constant $C$ such that one has $K(t,r) \le Ce^{-r^2/3t}$ for all $t \ge 1$.  Applying this, and bounding $\sinh(r)$ by $e^r$, one obtains
\[\E\left(d(Y_0,Y_t)1_{\left\lbrace d(Y_0,Y_t) \g At\right\rbrace}\right) \le \frac{C\vol(S^{n-1})}{k}\int\limits_{At}^{+\infty}\frac{r^2}{t^2}e^{(n-1)kr-r^2/3t}\d r\]
which bounding $e^{-r^2/3t}$ by $e^{-Ar/3}$ and choosing $A = 3nk$ yields
\[\cdots \le \frac{C\vol(S^{n-1})}{kt^2}\int\limits_{At}^{+\infty}r^2e^{-r}\d r = \frac{C\vol(S^{n-1})}{kt^2}(A^2t^2 + 2At + 2)e^{-At}\]
which, choosing $B$ appropriately, yields the desired bound for $t \ge 1$.
\end{proof}
\subsubsection{Inequalities}\label{sectioninequalities}

For a manifold $M$ with a compact quotient it can be shown that the limit
\[\lim\limits_{r \to +\infty} \frac{1}{r}\log(\vol(B_r(x)))\]
exists and has the same value for all $x \in M$.  If, for some Riemannian manifold possibly without a compact quotient, this limit above is zero then we say $M$ has subexponential growth.

We define the volume growth of a stationary random manifold $M$ as
\[v(M) = \liminf\limits_{r \to +\infty}\E\left( \frac{1}{r}\log(\vol(B_r(o_M)))\right).\]

By Bishop's inequality one has a uniform exponential upper bound on the volume of the ball of radius $r$ on any set of manifolds with uniformly bounded geometry.  This implies by dominated convergence that $M$ has subexponential growth almost surely then $v(M) = 0$.  On the other hand a uniform lower bound on volume is given by the fact that $M$ takes values in a space of manifolds with uniformly bounded geometry.  Hence, by Fatou's Lemma one has that if $v(M) = 0$ then $M$ satisfies 
\[\liminf\limits_{r \to +\infty} \frac{1}{r}\log(\vol(B_r(o_M))) = 0\]
almost surely.

We will bound entropy of a stationary random manifold from above and below in terms of its linear drift and volume growth. The upper bound was announced by Kaimanovich in the case of manifolds with a compact quotient (see \cite[Theorem 6]{kaimanovich1986}).  A sharper version of the lower bound, also for a manifolds with a compact quotient, was established by Ledrappier (see \cite{ledrappier2010}).  Analogous results for random walks on discrete groups also exist and have been sucessively improved by several authors, some of the first of these can be attributed to Varopoulos, Carne, and Guivarc'h (see \cite{gouezel-matheus-maucourant2012} and the references therein).  An analogous theorem for stationary random graphs is due to Benjamini and Curien (see \cite[Proposition 3.6]{benjamini-curien2012}).

\begin{theorem}\label{inequalities}
For all ergodic stationary random manifolds $M$ the following holds:
\[\frac{1}{2}\ell(M)^2 \le h(M) \le \ell(M)v(M).\]
\end{theorem}
\begin{proof}
We begin with the lower bound.

For this purpose fix $D \g 2$ and let $C$ be given by Theorem \ref{heatkernelupperbound} so that
\[q(t,x,y) \le C\exp\left(-\frac{d(x,y)^2}{Dt}\right)\]
holds on all manifolds in the range of $M$ for all $t \ge 1$.

Using this upper bound and Jensen's inequality we obtain
\begin{align*}\frac{1}{t}h_t(M) &= -\frac{1}{t} \int q(t,o,x)\log(q(t,o,x))\d x
\\& \ge \frac{1}{t}\int q(t,o,x)(\frac{d(o,x)^2}{Dt} - \log(C))\d x
\\& = \int q(t,o,x)(\frac{d(o,x)}{\sqrt{D}t})^2\d x - \frac{\log(C)}{t}
\\& \ge \left(\int q(t,o,x)\frac{d(o,x)}{\sqrt{D}t}\d x\right)^2 - \frac{\log(C)}{t}
\\& =  \frac{1}{D}(\frac{1}{t}\ell_t(M))^2 - \frac{\log(C)}{t}.
\end{align*}

Taking expectation and using Jensen's inequality once more one obtains
\[\E(\frac{1}{t}h_t(M)) \ge \frac{1}{D}\E\left((\frac{1}{t}\ell_t(M))^2\right) -\frac{\log(C)}{t} \ge \frac{1}{D}\E(\frac{1}{t}\ell_t(M))^2 -\frac{\log(C)}{t}\]
which by taking limit with $t \to +\infty$ yields
\[h(M) \ge \frac{1}{D}\ell(M)^2.\]

Letting $D$ decrease to $2$ one obtains $h(M) \ge \frac{1}{2}\ell(M)^2$ as claimed.

For the lower bound let $K(t,r)$ be the transition density of Brownian motion on $n$-dimensional hyperbolic space of constant curvature $-k^2$ where we assume all manifolds in the range of $M$ have curvature greater than or equal to $-k^2$ and dimension $n$.  By \cite[Theorem 2.2]{ichihara1988} one has
\[q(t,o_M,x) \ge K(t,d(o_M,x))\]
which combined with the upper bound given by Theorem \ref{heatkernelupperbound} yields
\[|\log(q(t,o_M,x))| \le \max(\log(C),\log(K(t,d(o_M,x)))\]
for some constant $C \g 0$ depending only on $\manifolds$, and all $t \ge 1$.

The density $K(t,r)$ is obtained by evaluating the heat kernel of hyperbolic space (curvature $-1$) at $t/2$ and $r/k$.  Hence from the lower bounds for the hyperbolic heat kernel given by \cite[Theorem 3.1]{davies-mandouvalos1988} one obtains constants $a,b,c \g 0$ depending only on $\manifolds$ such that
\[|\log(q(t,o_M,x))| \le a + b\log(t) + c(t+ r+ t/r)\]
for all $t \ge 1$.

Let $A,B$ be given by Lemma \ref{radialvariance}. We obtain uniformly over the range of $M$ (setting $r(x) = d(o_M,x)$) that
\begin{align*}
\frac{1}{t}\int\limits_{M \setminus B_{At}}&q(t,o_M,x)|\log(q(t,o_M,x))|\d x
\\ &\le \frac{a}{t} + \frac{b\log(t)}{t} + c\int\limits_{M \setminus B_{At}}q(t,o_M,x)(1+\frac{r(x)}{t} + \frac{r(x)^2}{t^2})\d x 
\\ &\le \frac{a}{t} + \frac{b\log(t)}{t} + c\int\limits_{M \setminus B_{At}}q(t,o_M,x)(\frac{r(x)^2}{A^2t^2}+\frac{r(x)^2}{At^2} + \frac{r(x)^2}{t^2})\d x 
\\ &\le \frac{a}{t} + \frac{b\log(t)}{t} + c(A^{-2} + A^{-1} + A)Be^{-At}
\end{align*}
for all $t \ge 1$.

This immediately implies (taking expectation and limit) that
\[h(M) = \lim\limits_{t \to +\infty}\E\left(-\frac{1}{t}\int\limits_{B_{At}(o_M)}q(t,o_M,x)\log(q(t,o_M,x))\d x\right).\]

The same upper bound on $|\log(q(t,o_M,x)|$ can now be applied in the ball of radius $At$ (where $r(x)/t$ can be bounded by $A$ before integrating).  One obtains that for any $L \le A$ letting $U_t$ be the annulus between radii $Lt$ and $At$ centered at $o_M$ one has
\[\frac{1}{t}\int\limits_{U_t}q(t,o_M,x)|\log(q(t,o_M,x))|\d x \le \frac{a}{t} + \frac{b\log(t)}{t} + c(1+A+A^2)\int\limits_{M \setminus B_{Lt}(o_M)}q(t,o_M,x)\d x\]
for all $t \ge 1$.

Assuming $L \g \ell(M)$ this last inequality implies
\[h(M) = \lim\limits_{t \to +\infty}\E\left(-\frac{1}{t}\int\limits_{B_{L t}(o_M)}q(t,o_M,x)\log(q(t,o_M,x))\d x\right).\]

To conclude notice that $\varphi(z) = -z\log(z)$ is convex on $z \g 0$.  Setting $v_t = \vol(B_{Lt}(o_M))$ and letting $p_t$ be the integral of $q(t,o_M,x)$ over $x$ in $B_{Lt}(o_M)$ one obtains using Jensen's inequality applied to normalized volume on the ball that
\[v_t \frac{1}{v_t}\int\limits_{B_{Lt}(o_M)}\varphi(q(t,o_M,x))\d x \le v_t\varphi(p_t/v_t) = p_t\log(v_t) - p_t\log(p_t) \le \log(v_t).\]

Taking expectation now yields
\[h(M) \le \liminf\limits_{t \to +\infty}\E\left(\frac{1}{t}\log(\vol(B_{Lt}(o_M)))\right) = L v(M)\]
for all $L \g \ell(M)$.  Which letting $L$ decrease to $\ell(M)$ proves the claimed upper bound.
\end{proof}

Recall that $\R^2$ endowed with the metric $\d s^2 = \d x^2 + (1+x^2)^2\d y^2$ has subexponential volume growth but admits the bounded harmonic function $\arctan(x)$ (this example was attributed to O. Chung by Avez).   Avez proved in 1976 that for manifolds with a transitive isometry group such an example is impossible (see \cite{avez1976}).

\begin{corollary}[Avez]
If $(M,g)$ is a connected Riemannian manifold with subexponential volume growth whose isometry group acts transitively then $M$ satisfies the Liouville property. 
\end{corollary}

A generalization of Avez's result to manifolds admitting a compact quotient under isometries was obtained by Varopolous (see \cite[Theorem 3]{varopoulos1986}).

\begin{corollary}[Varopoulos]
If $(M,g)$ is a Riemannian manifold with subexponential volume growth which admits a compact quotient under isometries then $M$ satisfies the Liouville property. 
\end{corollary}

An analogous result to the previous two for generic leaves of compact foliations (with respect to any harmonic measure) was announced by Kaimanovich and also follows from our theorem above (see Theorem 2 of \cite{kaimanovich1988} and the comments on page 307).   A particular interesting case is the horospheric foliation on the unit tangent bundle of a compact negatively curved manifold.

\begin{corollary}[Kaimanovich]
If $M$ is a negatively curved compact Riemannian manifold then almost every horosphere with respect to any harmonic measure for the horospheric foliation on the unit tangent bundle $T^1M$ satisfies the Liouville property.
\end{corollary}

Another interesting consequence of Theorem \ref{inequalities} is that $h(M) = 0$ if and only if $\ell(M) = 0$.  The following case was established by Karlsson and Ledrappier using a discretization procedure to reduce the proof to an analysis of a random walk on a discrete group (see \cite{karlsson-ledrappier2007}).

\begin{corollary}[Karlsson-Ledrappier]
Let $(M,g)$ be a manifold with bounded geometry that admits a compact quotient under isometries.  Then the $M$ satisfies the Liouville property if and only if its Brownian motion is non-ballistic (i.e. linear drift is zero).
\end{corollary}

One might be tempted to conjecture that a stationary random manifold of exponential growth must have non-constant bounded harmonic functions.   A counterexample is provided by Thurston's Sol-geometry (see \cite[pg. 304]{lyons-sullivan1984} or consider the case with $p = q = 1$ and drift parameter $a = 0$ in the central limit theorem of \cite{brofferio-salvatori-woess2012}).

\begin{example}[Lyons-Sullivan]
Let $M = \R^3$ endowed with the Riemannian metric $\d s^2 = e^{2z} \d x^2 + e^{-2z} \d y^2 + \d z^2$.  Then $M$ has a transitive isometry group, exponential volume growth, and satisfies the Liouville property.
\end{example}

\section{Brownian motion on stationary random manifolds.}
\label{sectionbrownianmotion}
In this section we construct Brownian motion on a stationary random manifold.  The technical steps consist of defining the corresponding path space (where both the manifold and path can vary) and proving that the Weiner measures one has on the paths over each manifold vary with sufficient regularity to define a global measure in path space over any harmonic measure on the set of manifolds under consideration (see Lemma \ref{pathspace} and Theorem \ref{liftedharmonic}).

Using this construction one obtains that linear drift can be defined pathwise allowing us to complete the proof of a fact used in the previous section (see Lemma \ref{lplusequalsl} and Corollary \ref{citethisfromchapter2}).

It also follows that a non-compact ergodic stationary random manifold must almost surely contain infinitely many disjoint diffeomorphic copies of any finite radius ball (see Theorem \ref{ghystypetheorem}).  This is in the spirit of the much more detailed result of Ghys which reduces the possible topologies of non-compact generic leaves of a foliation by surfaces to six possible types (see \cite{ghys1995}).

Generalizing Furstenberg's formula for the largest Lyapunov exponent of a product of random matrices (see \cite[Theorem 8.5]{furstenberg1963} and \cite[pg. 358]{ledrappier1984}) Karlsson and Ledrappier have given a formula expressing the rate of escape of random sequences in metric spaces using the expected increment of a Busemann function along the sequence (see \cite[Theorem 18]{karlsson-ledrappier2011}).  In our context we prove (see Theorem \ref{furstenbergformula}) a Furstenberg-type formula for the linear drift of Brownian motion on a stationary random manifold in terms of increments of a random Busemann function similar to \cite[Proposition 1.1]{ledrappier2010}.

In the last subsection we improve the lower bound $\frac{1}{2}\ell(M)^2 \le h(M)$ for entropy obtained in Theorem \ref{inequalities} to $2\ell(M)^2 \le h(M)$ in the case of certain stationary random Hadamard manifold (see Theorem \ref{lsquarelessthanh}).   In the case of a manifold with compact quotient this result was proved by Kaimanovich and Ledrappier, see \cite[Theorem 10]{kaimanovich1986} and \cite[Theorem A]{ledrappier2010}.  Equality implies that the gradient of almost every Busemann function at the origin must be collinear with that of a positive harmonic function, a condition which has strong rigidity consequences in the case of a single negatively curved manifold with compact quotient (see \cite{ledrappier2010} and \cite{ledrappier-shu2013}).  It is unknown to the author whether similar rigidity results can be obtained for stationary random manifolds.

\subsection{Brownian motion on stationary random manifolds}
\subsubsection{Path space}

We will construct a `path space' over a given set $\manifolds = \manifolds\left(d,r,\lbrace C_k\rbrace\right)$ of manifolds with uniformly bounded geometry.  Since later on we will be interested in time-reversal of Brownian motion we chose to consider paths whose domain is the entire real line instead of $[0,+\infty)$ as it was in previous sections.

For this purpose let $\widehat{\manifolds}'$ be the set of pairs $(M,\omega)$ where $M$ is a manifold in $\manifolds$ (denote by $o_M$ its basepoint and by $g_M$ its Riemannian metric) and $\omega:\R \to M$ is a continuous curve with $\omega_0 = o_M$.  And by $\widehat{\manifolds}$ denote the equivalence classes of elements of $\widehat{\manifolds}'$ where $(M,\omega)$ and $(M',\omega')$ are equivalent if there is a pointed isometry $f:M \to M'$ such that $\omega'  = f\circ \omega$.  We will not be careful in distinguishing elements of $\widehat{\manifolds}$ with their representatives $(M,\omega)$ in $\widehat{\manifolds}'$ since all our definitions will be invariant under the defined equivalence relationship.

Recall that a metric on a disjoint union of two metric spaces is admissible if it coincides with the given metrics when restricted to each half.  We mimick the definition of the distance on the Gromov space to turn $\widehat{\manifolds}$ into a metric space.
\begin{definition}\label{pathspacemetric}
Define the distance between two elements $(M,\omega)$ and $(M',\omega')$ in $\widehat{\manifolds}$ as either $1/2$ or, if such an $\epsilon$ exists, the infimum among all $\epsilon \in (0,1/2)$ such that that there exist an admissible metric on the disjoint union $M \sqcup M'$ with the following properties:
\begin{enumerate}
 \item $d(o,o') \l \epsilon$.
 \item $d(B_{1/\epsilon}(o),M') \l \epsilon$ and $d(M,B_{1/\epsilon}(o')) \l \epsilon$.
 \item $d(\omega_t,\omega'_t) \l \epsilon$ for all $t \in [-1/\epsilon, 1/\epsilon]$.
\end{enumerate}
\end{definition}

Define for each $t \in \R$ the shift map $\shift^t:\widehat{\manifolds} \to \widehat{\manifolds}$ by
\[\shift^t (M,\omega) = (M,\shift^t\omega)\]
where $(\shift^t\omega)_s = \omega_{t+s}$ for all $s$ and the basepoint of $M$ on the right hand has been changed to $\omega_t$ (previously it was $\omega_0$).

Also, one has a projection $\pi:\widehat{\manifolds} \to \manifolds$ which associates to each $(M,\omega)$ the unique pointed manifold in $\manifolds$ isometric (with basepoint) to $(M,\omega_0)$.

Recall that a sequence $(M_n,o_n,g_n)$ of pointed manifolds is said to converge smoothly to $(M,o,g)$ if there exists an exaustion $U_n$ of $M$ be increasing relatively compact open sets and a sequence of smooth embeddings $\varphi_n:U_n \to M_n$ with $\varphi_n(o) = o_n$ such that the pullback metrics $\varphi_n^*g_n$ converge smoothly on compact sets to $g$.   

\begin{lemma}\label{pathspace}
Let $\manifolds$ and $\widehat{\manifolds}$ be as above.   Then $\manifolds$ is compact and metrizable when endowed with the topology of smooth convergence, $\widehat{\manifolds}$ is a complete metric space with the distance defined above, and the projection $\pi:\widehat{\manifolds} \to \manifolds$ is continuous and surjective.   Furthermore, for each $t \in \R$ the shift map $\shift^t$ is a self homeomorphism of $\widehat{\manifolds}$.
\end{lemma}
\begin{proof}
By Theorem \ref{uniformlyboundedgeometry} one has that $\manifolds$ is compact with respect to pointed Gromov-Hausdorff convergence and that this convergence is equivalent to smooth convergence on $\manifolds$.  Since Gromov-Hausdorff convergence is metric this establishes the claim on $\manifolds$. 

The continuity of the projection follows because we have defined $\widehat{\manifolds}$ to be larger than the pointed Gromov-Hausdorff distance between the projected manifolds.  Continuity of the shift map can be verified directly from the definition.  Hence it remains to establish that $\widehat{\manifolds}$ is separable and complete.

We begin by establishing completeness.  For this purpose take a sequence $(M_n,\omega^n)$ in $\widehat{\manifolds}$ such that $d((M_n,\omega^n),(M_{n+1},\omega^{n+1})) \l 2^{-n}$ for all $n$.  For each $n$ let $d_n$ be an admissible metric on the disjoint union $M_n \sqcup M_{n+1}$ satisfying the conditions of Definition \ref{pathspacemetric} for $\epsilon = 2^{-n}$.  Using these distances define a distance on the infinite disjoint union $X = \bigsqcup M_n$ such that if $x \in M_n$ and $y \in M_{n+p}$ the distance between them is the infimum over sequences $x = x_0,x_1,x_2,\ldots,x_p = y$ with $x_i \in M_{n+i}$ of $d_n(x_0,x_1)+\cdots + d_{n+p-1}(x_{n+p-1},x_{n+p})$.  Let $\widehat{X}$ be the completion of $X$, $M = \widehat{X} \setminus X$, and let $\omega$ be the local uniform limit of the curves $\omega^n$ in $\widehat{X}$.

One can verify that $\omega_t$ belongs to $M$ for all $t$.  Furthermore we have shown in Theorem \ref{gromovspacetheorem} that the sequence $M_n$ converges in the pointed Gromov-Hausdorff sense to $M$.  It follows that $M$ is isometric to a pointed manifold in $\manifolds$.  Since the distance between $(M_n,\omega^n)$ and $(M,\omega)$ is less than $2^{-(n-1)}$ this establishes that $\widehat{\manifolds}$ is complete as claimed.

Take a countable dense subset $D$ of $\manifolds$ and consider a countable set $\widehat{D}$ consisting of pairs $(M,\omega)$ where $M$ ranges over $D$ and $\omega$ over a countable dense (with respect to local uniform convergence) subset of curves in each $M$.  We claim that $\widehat{D}$ is dense in $\widehat{\manifolds}$.

To establish this fact fix $(M,\omega)$ in $\widehat{\manifolds}$ and $\epsilon \g 0$.  Since $\omega$ is continuous we may partition the interval $[-1/\epsilon,1/\epsilon]$ into a finite number of disjoint intervals $I_i = [t_i,t_{i+1}]$ such that the diameter of $\lbrace \omega_t: t \in I_i\rbrace$ is less than $\epsilon/3$ for all $i$. Take $\delta \in (0,\epsilon/3)$ so that $1/\delta$ is larger than the diameter of $\lbrace \omega_t : t \in [-1/\epsilon,1/\epsilon]\rbrace$ and notice that there exists $(M',o',g')$ in $D$ at distance less than than $\delta$ from $(M,\omega_0,g)$ ($g$ being the metric on $M$).  It follows that there exists an admissible distance on the disjoint union $M \sqcup M'$ such each point $\omega_{t_i}$ is at distance less than $\epsilon/3$ from some point $p_i \in M'$.  Interpolating between the $p_i$ with geodesic segments one obtains a curve $\omega'_t$ in $M'$ at distance less than $\epsilon$ from $\omega_t$ for all $t \in [-1/\epsilon,1/\epsilon]$.  Arbitrarily close to $\omega'$ in the topology of local uniform convergence there are curves $\omega''$ such that $(M',\omega'')$ belongs to $\widehat{D}$.  This proves that $\widehat{D}$ is dense as claimed.
\end{proof}
\subsubsection{Brownian motion}

By the topology of local uniform convergence on an interval $I \subset \R$ on $\widehat{\manifolds}$ we mean the topology generated by replacing $t \in [-1/\epsilon,1/\epsilon]$ in condition 3 of Definition \ref{pathspacemetric} by $t \in [-1/\epsilon,1/\epsilon] \cap I$ .  We denote by $\F_0^\infty$ the $\sigma$-algebra generated by the topology of local uniform convergence on $[0,+\infty)$ on $\widehat{\manifolds}$.

Recall that $C(\R,M)$ denotes the space of continuous curves on a manifold $M$ and denote by $\F_0^\infty(M)$ the $\sigma$-algebra on $C(\R,M)$ generated by the topology of local uniform convergence on $[0,+\infty)$.  Each element $\omega$ of $C(\R,M)$ is naturally associated to the element $(M,\omega)$ of $\widehat{\manifolds}$ though this association is not necessarily injective (e.g. applying any isometry to a curve yields the same element of $\widehat{\manifolds}$).

\begin{lemma}\label{liftedweiner}
The natural projection $\omega \mapsto (M,\omega)$ from the space of continuous curves $C(\R,M)$ on a manifold $M \in \manifolds$ into $\widehat{\manifolds}$ is continuous when both are endowed with the topology of local uniform convergence on the same interval $I \subset \R$.  In particular each probability measure on $\F_0^\infty(M)$ yields a probability on $\left(\widehat{\manifolds},\F_0^\infty\right)$.
\end{lemma}
\begin{proof}
The uniform distance between two curves $\omega,\omega'$ in $M$ on a compact interval $[a,b]$ is larger than or equal to the corresponding distance on $\widehat{\manifolds}$ endowed with the topology of uniform convergence on $[a,b]$ (this follows directly from Definition \ref{pathspacemetric}).
\end{proof}

Given $(M,o,g) \in \manifolds$ we denote by $P_{(M,o,g)}$ the pushforward of Wiener measure corresponding to Brownian motion starting at $o$ to the $\sigma$-algebra $\F_0^\infty$ on $\widehat{\manifolds}$.  The Markov property of Brownian motion on $M$ takes the following form:
\begin{corollary}\label{foliatedmarkov}
Let $(M,o,g) \in \manifolds$, then for any $t \g 0$ and $A \in F_t^\infty$ one has
\[\P_{(M,o,g)}(A) = \int q(t,o,x) P_{(M,x,g)}(A) \d x.\]
\end{corollary}

By a Brownian motion on a stationary random manifold we mean a random element $(M,\omega)$ of $\widehat{\manifolds}$ with distribution $\widehat{\mu}$ given by the following theorem for some harmonic measure $\mu$ on $\manifolds$.

\begin{theorem}\label{liftedharmonic}
For each probability measure $\mu$ on $\manifolds$ there is a unique probability $\widehat{\mu}$ on $\left(\widehat{\manifolds},\F_0^\infty\right)$ defined by
\[\widehat{\mu} = \int \P_{(M,o,g)} \d \mu(M,o,g).\]
Furthermore the following two properties hold
\begin{enumerate}
 \item If a probability $\mu$ on $\manifolds$ is harmonic then its lift $\widehat{\mu}$ can be extended uniquely to a shift invariant Borel probability on $\widehat{\manifolds}$.
 \item If a probability $\mu$ on $\manifolds$ is harmonic and ergodic then the unique shift invariant extension of $\widehat{\mu}$ is ergodic.
\end{enumerate}
\end{theorem}
\begin{proof}
We have shown in Lemma \ref{liftedweiner} that the Weiner measures on each $(M,o,g)$ in $\manifolds$ lifts to a probability $P_{(M,o,g)}$ on $\widehat{\manifolds}$.  If we show that 
\begin{equation}\label{shouldbecontinuous}
(M,o,g) \mapsto \int f(M,\omega) \d P_{(M,o,g)}(N,\omega)
\end{equation}
is continuous for all continuous bounded $f:\widehat{\manifolds} \to \R$ then $\widehat{\mu} = \int \P_{(M,o,g)}\d \mu$ is well defined as an element of the dual of the space of bounded continuous functions.  By the Riesz representation theorem this means that $\widehat{\mu}$ is well defined as a probability on the Stone-Cech compactification of $\widehat{\manifolds}$.  After this it suffices to show that for each $\epsilon \g 0$ there is a compact subset $K$ of $\widehat{\manifolds}$ with $P_{(M,o,g)}(K) \ge 1-\epsilon$ for all  $(M,o,g)$ in $\manifolds$ to obtain that in fact $\widehat{\mu}$ is a probability on $\widehat{\manifolds}$.

To begin we define the space of orthonormal frames $O(\manifolds)$ as the set of tuples 
\[(M,o,g,v_1,\ldots,v_d)\]
where $(M,o,g)$ is in $\manifolds$ and $v_1,\ldots,v_d$ is an orthonormal basis of the tangent space $T_oM$.  This space is considered up to equivalence by isometries which respect the basepoints and orthonormal frames.  A sequence $(M^n,o^n,g^n,v_1^n,\ldots,v_d^n)$ converges to $(M,o,g,v_1,\ldots,v_d)$ if there exists an exhaustion $U_n$ of $M$ and smooth embeddings $\varphi_n:U_n \to M_n$ satisfying the usual properties for smooth convergence plus that $|D\varphi_n v_i - v_i^n| \to 0$ for $i = 1,\ldots, d$.  We omit the verification that $O(\manifolds)$ is separable, compact and the projection is continuous.  Similarly we define the space of frames $F(\manifolds)$ by dropping the condition that $v_1,\ldots,v_d$ be orthonormal (this space is no longer compact).

On the frame bundle $F(M)$ of any manifold $(M,o,g) \in \manifolds$ there are $d$ unique smooth vector fields $H_1,\ldots, H_d$ with the property that the flow defined by $H_i$ corresponds to parallel transport of each orthonormal frame along the geodesic whose initial condition is the $i$-th vector of the frame.  Any solution to the Stratonovich differential equation
\begin{equation}\label{brownianequation} 
d X_t = \sum\limits_{i = 1}^d H_i(X_t)\circ \d W_t^i
\end{equation}
driven by a standard Brownian motion $(W_t^1,\ldots,W_t^d)$ on $\R^d$ and starting at an orthonormal frame, projects to a Brownian motion on $M$.  In particular the distribution of the solution starting at any orthonormal frame $v_1,\ldots, v_d$ over the basepoint $o$ is the Weiner measure $P_o$ and lifts to $P_{(M,o,g)}$ on $\widehat{\manifolds}$.  We will show that $P_{(M,o,g)}$ depends continuously on $(M,o,g)$ in $\manifolds$ by showing that it depends continuously on the points in $O(\manifolds)$.  Since points projecting to the same element of $\manifolds$ have the same associated measure in $\widehat{\manifolds}$ this approach may seem unnecessarily complicated, however it allows us to use standard theorems on continuity of solutions to rough differential equations and simultaneously solve the problem of finding compact sets with large probability for all $P_{(M,o,g)}$.

Choose $p \in (2,3)$ and consider for each $T \g 0$ the Polish space $\Omega$ of geometric rough paths with locally finite $p$-variation (by which we mean that restricted to each interval $I = [0,T]$ the path belongs to $C^{p-var}(I,G^2(\R^d))$ as defined in \cite[Definition 9.15]{friz-victoir2010}).  Let $\nu$ be the probability measure on $\Omega G_p(\R^d)$ which is the distribution of `Stranovich enhanced Brownian motion' (see \cite[Section 13.2]{friz-victoir2010}).  Given $(M,o,g,v_1,\ldots, v_d)$ and a path $\alpha$ in $\Omega$ there is a unique curve $\omega:[0,+\infty) \to M$ which is the projection to $M$ of the solution to Equation \ref{brownianequation} driven by the rough path $\alpha$ (this is the existence and uniqueness theorem for rough differential equations of Terry Lyons, we will use the version given by \cite[Theorem 10.26]{friz-victoir2010}).  We will show that $(M,\omega)$ is continuous as a function from $\Omega \times O(\manifolds)$ to $\widehat{\manifolds}$ (where the later is endowed 
with the topology of local uniform convergence).  

Assuming the continuity claim has been established we notice that, because the pushforward of $\nu$ under the map $\alpha \mapsto \omega$ is the Weiner measure corresponding to $o$, it would follow by dominated convergence that the function defined by Equation \ref{shouldbecontinuous} is continuous for all bounded $f$.  Furthermore, by letting $\alpha$ vary in a compact subset of $\Omega$ with $\nu$-probability greater than $1-\epsilon$ one would obtain that there is a compact subset of $\widehat{\manifolds}$ with probability greater than $1-\epsilon$ for all $P_{(M,o,g)}$.

We now establish the continuity claim.  Consider a sequence in $O(\manifolds)$ such that 
\[(M^n,o^n,g^n,v_1^n,\ldots, v_d^n) \to (M,o,g,v_1,\ldots,v_d)\]
when $n \to +\infty$, and let $\varphi_n:U_n \to M^n$ be the corresponding embeddings of an exaustion of $M$.  Given $\alpha \in \Omega$ and a compact interval $[0,T]$ we may take a compact $d$-dimensional submanifold with boundary $N \subset F(M)$  such that the solution to Equation \ref{brownianequation} driven by $\alpha$ starting at $v_1,\ldots,v_d$ remains in $N$ on $[0,T]$.   We embed $N$ into some $\R^N$ and consider smooth compactly supported extensions of the horizontal vector fields $H_i$ to $\R^N$.  The pullbacks of the horizontal vector fields $H_i^n$ on $F(M_n)$ under $\varphi_n$ restricted to $N$ can be extended to all of $\R^N$ in such a way that they share a compact support and converge smoothly to the corresponding $H_i$ and the pullbacks of the orthonormal frames $v_i^n$ eventually belong to $N$ and converge to $v_1,\ldots,v_d$.  Hence by \cite[Theorem 10.26]{friz-victoir2010} 
the solutions to Equation \ref{brownianequation} driven by $\alpha$ and the vector fields $H_i^n$ starting at the pullbacks of $v_1^n,\ldots,v_d^n$ converge to the corresponding solution driven by the vector fields $H_i$ uniformly on $[0,T]$.  This implies continuity of the solution with respect to $(M,o,g,v_1,\ldots,v_d)$ in $O(\manifolds)$\begin{footnote}{Here we use the fact that if there is an exhaustion $U_n$ of $M$ and embeddings $\varphi_n:U_n \to M_n$ defining the conditions for smooth convergence, and furthermore there are curves $\omega$ on $M$ and $\omega_n$ on $M_n$ such that $\varphi_n^{-1}\circ \omega_n$ converges to $\omega$ uniformly on compact sets, then $(M_n,\omega_n)$ converges to $(M,\omega)$ as elements of $\widehat{\manifolds}$.  The proof amounts to repeating the arguments in Lemma \ref{cristinaslemma}.}\end{footnote}.   Continuity with respect to $\alpha$ for a fixed manifold with frame in $O(\manifolds)$ also follows from \cite[Theorem 10.26]{friz-victoir2010}.

We have thus far established the existence of $\widehat{\mu}$ for each probability $\mu$ on $\manifolds$. It remains to interpret harmonicity and ergodicity of $\mu$ with properties of $\widehat{\mu}$ relative to the shift maps.

If $\mu$ is harmonic then by Corollary \ref{foliatedmarkov} one has for all $A \in \F_t^\infty$ that $\widehat{\mu}(A) = \shift^t_*\widehat{\mu}(A)$ from this it follows that the pushforward measure $\shift^t_*\widehat{\mu}$ (defined on $\F_{-t}^\infty$ coincides with $\mu$ on $\F_0^\infty$.  Hence we may extend $\widehat{\mu}$ uniquely to a shift invariant probability $\F_t^\infty$ for all $t \in \R$. Since all these extensions are compatible Tulcea's extension theorem implies that there is a unique shift invariant Borel extension.

Suppose that $\mu$ is ergodic, harmonic, and that the unique shift invariant extension of $\widehat{\mu}$ to the Borel $\sigma$-algebra is not ergodic.  Then by definition one can obtain two distinct shift invariant probabilities such that
\[\widehat{\mu} = \alpha\widehat{\mu}_1 + (1-\alpha)\widehat{\mu}_2\]
for some $\alpha \in (0,1)$.  Since $\widehat{\mu}$ projects to $\mu$ this allows one to express $\mu$ as a non-trivial convex combination of the projection $\mu_1$ and $\mu_2$ of $\widehat{\mu}_1$ and $\widehat{\mu}_2$ respectively.  Using Corollary \ref{foliatedmarkov} one sees that because the $\widehat{\mu}_i$ are shift invariant $\mu_1$ and $\mu_2$ are harmonic.  But because $\mu$ is ergodic one must have $\mu_1 = \mu_2 = \mu$ contradicting the fact that $\widehat{\mu}_1 \neq \widehat{\mu}_2$.  Hence each ergodic harmonic measure $\mu$ on $\manifolds$ lifts to (after taking the unique shift invariant extension) an ergodic shift invariant probability $\widehat{\mu}$ as claimed.
\end{proof}

\subsubsection{Pathwise linear drift}

Suppose $X_t$ is a Brownian motion on a manifold with bounded geometry. By the results of Ichihara (see \cite[Example 2.1]{ichihara1988}) one has that $\limsup\limits_{t \to +\infty} d(X_0,X_t)/t$ is finite.  However, the limit
\[\lim\limits_{t \to +\infty} \frac{d(X_0,X_t)}{t}\]
need not exists almost surely.  Furthermore even if the above limit exists it might be random (i.e. take different values with positive probability).  For a concrete example consider a metric on $\R^2$ which contains isometrically embedded copies of a half space with constant curvature $-1$ and a half space of constant curvature $-2$.  For any Brownian motion with respect to such a metric the above limit will take two distinct values with positive probability.

The construction of the Brownian motion process on $\widehat{\manifolds}$ implies that manifolds such as the previously discussed cannot be generic with respect to any ergodic harmonic measure.

\begin{theorem}\label{pathwisedrift}
Let $(M,\omega)$ be a Brownian motion an an ergodic stationary random manifold.  Then there one has
\[\lim\limits_{t \to +\infty} \E\left(\frac{d(\omega_0,\omega_t)}{t}\right) = \ell(M)\]
and 
\[\lim\limits_{t \to +\infty} \frac{d(\omega_0,\omega_t)}{t} = \ell(M)\]
almost surely, where $\ell(M)$ is the linear drift of the stationary random manifold $M$.
\end{theorem}
\begin{proof}
Let $\mu$ be the distribution of $M$ and consider the shift invariant lift $\widehat{\mu}$ of $\mu$ to $\widehat{\manifolds}$ given by Theorem \ref{liftedharmonic}.  We define the function $d_{s,t}$ on $\widehat{\manifolds}$ for $s \l t$ by
\[d_{s,t}(M,\omega) = d(\omega_s,\omega_t).\]

The triangle inequality implies that $d_{s,u} \le d_{s,t} + d_{s,u}$, added to the fact that $\widehat{\mu}$ is shift invariant one obtains that the family $\lbrace d_{s,t}\rbrace$ is a stationary subadditive process with respect to the probability $\widehat{\mu}$.  Therefore by Kingman's subadditive ergodic theorem  (see Theorems 1 and 5 of \cite{kingman1968}) the limits
\[\lim\limits_{t \to +\infty}\frac{d_{s,t}}{t}\]
exists $\widehat{\mu}$ almost surely and in $L^1$.  Since the limit is almost surely shift invariant it must be constant if $\widehat{\mu}$ is ergodic.
\end{proof}

As a consequence of the previous theorem we can complete the proof that $\ell^+(M) = \ell(M)$ for all ergodic stationary random manifolds $M$, a fact which was used in the previous section.

\begin{corollary}\label{citethisfromchapter2}
Let $M$ be an ergodic stationary random manifold then one has
\[\lim_{t \to +\infty}\E\left( \int\limits_{B_{Lt}(o_M)}q(t,o_M,x)\d x \right ) = 1\]
for all $L \g \ell(M)$.
\end{corollary}
\begin{proof}
Let $\mu$ be the distribution of $M$, $\manifolds$ the support of $\mu$, and $\widehat{\mu}$ the lift of $\mu$ to $\widehat{\manifolds}$.  The quantity inside the limit is equal to
\[\widehat{\mu}\left( \left\lbrace (M,\omega) \in \widehat{\mu}: d(\omega_0,\omega_t) \l Lt\right\rbrace \right)\]
which converges to $1$ by Theorem \ref{pathwisedrift}.
\end{proof}

\subsubsection{Recurrence and the correspondence principle}

Szemeredi's theorem states that if a set of integers $A$ has positive upper density then it contains arbitrarily long arithmetic progressions.  Furstenberg's proof of this result begins by associating to $A$ a shift invariant probability measure $\mu$ on the space of infinite strings of zeros and ones with the property that if the measure of the set of strings beginning with a particular finite pattern is positive then this pattern appears in $A$.  In particular the set $U$ of strings beginning with $1$ has positive probability and Szermeridi's theorem is shown to be equivalent to the property that for each $k$ there exist $n$ such that $\mu(U \cap T^{-n}U \cap \cdots \cap T^{-nk}U) \g 0$ where $T$ is the shift transformation (see \cite[pg. 178]{einsiedler-thomas2011}).

The above idea of Furstenberg has been abstracted to a general albeit informal `correspondence principle' which might be roughly stated as follows: The combinatorial properties of a concrete mathematical object can sometimes be codified by a measure preserving dynamical system.  Hence ergodic theorems can yield proofs of combinatorial statements and vice-versa.

In this subsection we try to apply the above principle to a fixed bounded geometry manifold $(M,g)$.   The idea is that the closure of the set of manifolds $(M,x,g)$ where $x$ varies over all of $M$, supports at least one harmonic measure `diffused from $M$'.  If for some $o \in M$ the pointed manifold $(M,o,g)$ belongs to the support of this measure then this imposes strong `recurrence' properties on the Riemannian metric and topology of $M$ forcing a certain finite radius `pattern' in the manifold to appear infinitely many times.  Perhaps the strongest result which has been obtained by this type of reasoning is the theorem of Ghys which states that a non-compact leaf of a compact foliation by surfaces which is generic with respect to a harmonic measure can only have one out of six possible topologies (see \cite{ghys1995}).  We will prove a weaker result which illustrates the same principle.  We begin with a definition.

\begin{definition}
Let $(M,o,g)$ be a pointed bounded geometry manifold in some set of manifolds with uniformly bounded geometry $\manifolds$.  A diffusion measure $(M,o,g)$ is any weak limit when $t \to +\infty$ of a subsequence of the measures
\[\frac{1}{t}\int_0^t (\pi \circ \shift^s)_*P_{(M,o,g)}\d s\]
where $\pi$ is the projection from $\widehat{\manifolds}$ to $\manifolds$ and the shift maps on $\widehat{\manifolds}$ are denoted by $\shift^t$.
\end{definition}

We might say a manifold is `recurrent' if it belongs to the support of at least one of its diffusion measures.  We begin by showing that, with this definition, all manifolds which are generic with respect to an ergodic harmonic measure are recurrent and in fact have a unique diffusion measure.  
\begin{lemma}
Let $\mu$ be an ergodic harmonic measure on a space of manifolds with uniformly bounded geometry $\manifolds$.  Then $\mu$ is the unique diffusion measure for $\mu$-almost every $(M,o,g) \in \manifolds$.
\end{lemma}
\begin{proof}
Consider the ergodic shift invariant lift $\widehat{\mu}$ of $\mu$ to $\widehat{\manifolds}$ given by Theorem \ref{liftedharmonic} and let $\lbrace U_n\rbrace$ be a countable basis of the topology on $\manifolds$ by open sets whose boundaries have zero $\mu$ measure.  By Birkhoff's ergodic theorem one has for $\widehat{\mu}$ almost every $(M,\omega)$ in $\widehat{\manifolds}$ that
\[\mu(U_n) = \widehat{\mu}\left(\pi^{-1}(U_n)\right) = \lim\limits_{t \to +\infty}\frac{1}{t} \int_0^t 1_{U_n}\left(\pi\circ \shift^s (M,\omega)\right) \d s\]
for all $n$.  

This implies that 
\[\mu = \lim\limits_{t \to +\infty}\frac{1}{t}\int_0^t (\pi \circ \shift^s)_*P_{(M,\omega_0,g)}\d s\]
so $\mu$ is a diffusion measure for $\widehat{\mu}$ almost every $(M,\omega_0,g)$.  This implies the claim since $\widehat{\mu}$ projects to $\mu$.
\end{proof}

The following theorem shows, for example, that the plane with one handle cannot be the generic leaf of a compact foliation no matter what bounded geometry Riemannian metric we put on it.
\begin{theorem}\label{ghystypetheorem}
Suppose that $(M,o,g)$ is a non-compact bounded geometry manifold which belongs to the support of its only diffusion measure $\mu$.  Then for each $r \g 0$ there are infinitely many disjoint diffeomorphic copies of $B_r(o)$ embedded in $M$.
\end{theorem}
\begin{proof}
Given $r$ we can choose a neighborhood $U$ of $(M,o,g)$ sufficiently small so that for all $(M',o',g')$ in $U$ there is a diffeomorphism $\varphi: B_r(o) \to M'$ with $\varphi(o) = o'$ and such that $\varphi(B_r(o))$ is contained in $B_{2r}(o')$.

Consider the set $A = \lbrace x \in M: (M,x,g) \in U\rbrace$ and notice that one has
\[\lim\limits_{t \to +\infty} \frac{1}{t} \int_0^t 1_A(\omega_s)\d s = \mu(U) \g 0\]
for $P_o$ almost every Brownian path $\omega$.

On the other hand, since $M$ is non-compact and has bounded geometry it has infinite volume.  This implies that the fraction of time spent by Brownian motion in any compact subset of $M$ converges to $0$ as $t \to +\infty$.  It follows that $A$ is unbounded and hence $M$ contains infinitely many disjoint copies of $B_r(o)$.
\end{proof}

The main result of \cite{ghys1995} implies that the non-compact generic leaves with respect to any ergodic harmonic measure on a compact foliation by surfaces have either zero or infinite genus and either one, two, or infinitely many ends.   The first part of this result follows from the theorem above.

\begin{corollary}[Ghys]\label{ghyscorollary}
Any non-compact generic leaf with respect to an ergodic harmonic measure on a compact foliation by surfaces must have either zero or infinite genus.
\end{corollary}

\subsection{Busemann functions and linear drift}
\subsubsection{Busemann functions}

The Busemann function $\xi_x:X \to \R$ associated to a point $x$ in a pointed proper metric space $(X,o)$ is defined by
\[\xi_x(y) = d(x,y) - d(x,o)\]
and the Busemann compactification of $X$ is the closure of all such functions in the topology of local uniform convergence (this is equivalent to pointwise convergence since all Busemann functions are $1$-Lipschitz).  Given a set of pointed manifolds with uniformly bounded geometry $\manifolds$ we will show in this section that the linear drift of a harmonic measure on $\manifolds$ can be expressed in terms of the increment of a `random Busemann function' on a fixed time interval of Brownian motion.

As a first step we must construct a space containing the Brownian paths $\widehat{\manifolds}$ and Busemann functions.  For this purpose we define $\widehat{\manifolds}^1$ as the set of triplets $(M,\omega,\xi)$ where $(M,\omega)$ is in $\widehat{\manifolds}$ and $\xi:M \to \R$ is $1$-Lipschitz and satisfies $\xi(\omega_0) = 0$.  Elements of $\widehat{\manifolds}^1$ are considered up to the usual equivalence by isometries preserving $\omega$ and $\xi$.

The distance on $\widehat{\manifolds}^1$ is defined similarly to that of $\widehat{\manifolds}$ adding only one condition.
\begin{definition}\label{busemannspacemetric}
Define the distance between two elements $(M,\omega,\xi)$ and $(M',\omega',\xi')$ in $\widehat{\manifolds}$ as either $1/2$ or, if such an $\epsilon$ exists, the infimum among all $\epsilon \in (0,1/2)$ such that that there exist an admissible metric on the disjoint union $M \sqcup M'$ with the following properties:
\begin{enumerate}
 \item $d(o,o') \l \epsilon$.
 \item $d(B_{1/\epsilon}(o),M') \l \epsilon$ and $d(M,B_{1/\epsilon}(o')) \l \epsilon$.
 \item $d(\omega_t,\omega'_t) \l \epsilon$ for all $t \in [-1/\epsilon, 1/\epsilon]$.
 \item $|\xi(x)-\xi'(x')| \l \epsilon$ whenever $d(x,x') \l \epsilon$ and either $x \in B_{1/\epsilon}(o)$ or $x' \in B_{1/\epsilon}(o')$.
\end{enumerate}
\end{definition}

Usually the Busemann functions of a manifold $M$ are defined as the locally uniform limits of the functions $\xi_x$ (where $x \in M$).  We define the generalized Busemann space $\widehat{\manifolds}^b$ as the closure of the Busemann functions $(M,\omega,\xi_x)$ (where $x \in M$) in the above metric space.  The elements of this space will be called generalized Busemann functions.

\begin{lemma}\label{busemannspace}
The metric space $\widehat{\manifolds}^b$ is complete and separable, and the projection $\widehat{\pi}:\widehat{\manifolds}^b \to \widehat{\manifolds}$ is continuous surjective and proper (i.e. preimage of any compact set is compact).
\end{lemma}
\begin{proof}
Completeness of $\widehat{\manifolds}^b$ follows from an argument similar to that of Lemma \ref{pathspace}.  We consider a sequence $(M^n,\omega^n,\xi^n)$ with distance between consecutive elements less than $2^{-n}$.  We take an admissible metric $d_n$ on $M_n \sqcup M_{n+1}$ satisfying the conditions of Definition \ref{busemannspacemetric} and use it to define a metric on the countable disjoint union $\bigsqcup M_n$ by setting $d(x,y)$ with $x \in M^n$ and $y \in M^{n+p}$ to be the infimum of $d_n(x,x_1)+\cdots + d_{n+p-1}(x_{p-1},y)$ over all chains $x_1,\ldots, x_{p-1}$ with $x_i \in M_{n+i}$ for $i = 1,\ldots,p-1$.  We take $X$ to be the completion of this disjoint union and $M = X \setminus \bigsqcup M^n$.   As in Lemma \ref{pathspace} one has that $\omega^n$ converge locally uniformly to a curve $\omega$ in $M$.  

Take a point $x \in M$ and assume that $x = \lim x_n$ where $x_n \in M^n$ and $d(x_n,x_{n+1}) \l 2^{-n}$ for each $n$.  Then because $d(o_n,x_n) \to d(o,x)$ we will have for all $n$ large enough that $|\xi^n(x_n)-\xi^{n+1}(x_{n+1})| \l 2^{-n}$.  This implies that $\lim \xi^n(x_n)$ exists.  Furthermore, if $y_n$ is another sequence with the same properties converging to $x$ then $d(y_n,x_n) \to 0$ and since each $\xi^n$ is $1$-Lipschitz one has that $\lim \xi^n(y_n) = \lim \xi^n(x_n)$.  We define $\xi(x) = \lim \xi^n(x_n)$.  One verifies that $\xi$ is 1-Lipschitz and defined on a dense subset of $M$, therefore it extends uniquely to a continuous function on $M$.

By definition the function $F$ on $X$ which coincides with $\xi^n$ on each $M^n$ and with $\xi$ on $M$ is continuous.  Therefore it is uniformly continuous on compact sets.  Hence, given $\epsilon \g 0$ we may find $\delta \g 0$ such that if $x \in B_{\epsilon + 1/\epsilon}(o)$ and $d(x,y) \l \delta$ then $|F(x) - F(y)| \l \epsilon$.  The admissible distance on $M_n \sqcup M$ which equals $d(x,y) + \epsilon - \delta$ whenever $x \in M_n$ and $y \in M$ and coincides with $d$ otherwise, shows that the distance between $(M_n,\omega^n, \xi^n)$ and $(M,\omega,\xi)$ is less than $\epsilon$ for all $n$ large enough.  Hence $\widehat{\manifolds}^b$ is complete.

Surjectivity of $\widehat{\pi}$ is immediate. Separability of $\widehat{\manifolds}^b$ follows from that of $\widehat{\manifolds}$ once we show that $\widehat{\pi}$ is proper (the preimage of each point in a dense countable subset of $\widehat{\manifolds}$ is compact, hence has a countable dense set, the union of these sets is dense in $\widehat{\manifolds}^b$).

To establish that $\widehat{\pi}$ is proper it suffices to show that if $(M^n,\omega^n,\xi_{x_n})$ is a sequence of Busemann functions in $\widehat{\manifolds}^b$ such that $(M^n,\omega^n)$ converges to $(M,\omega)$ in $\widehat{\manifolds}$, then there is a $1$-Lipschitz function $\xi$ on $M$ such that a subsequence of $(M^n,\omega^n,\xi_{x_n})$ converges to $(M,\omega,\xi)$.

Repeating the construction above we may assume without loss of generality that there is an admissible metric on $\bigsqcup M_n \sqcup M$ such that the curves $\omega^n$ converge locally uniformly to $\omega$ and the metric restricted to $M_n \sqcup M_{n+1}$ satisfies properties 1,2, and 3 of Definition \ref{busemannspacemetric} for $\epsilon = 2^{-n}$.  The same conditions are verified by the metric when restricted to $M_n \sqcup M$ for $\epsilon = 2^{-n+1}$.

Notice that each $\xi_{x_n}$ is defined in terms of the distance on $M^n$.  Hence it extends, using the admissible distance, to a $1$-Lipschitz function $F_n$ on all of $M \sqcup \bigsqcup M_n$.  Since the sequence $F_n$ is equicontinuous there is a locally uniformly convergent subsequence $F_{n_k}$ which converges to a limit $F$.  We define $\xi$ as the restriction of $F$ to $M$.  An argument similar to the one above shows that $(M^{n_k},\omega^{n_k},\xi_{x_{n_k}})$ converges to $(M,\omega,\xi)$.
\end{proof}

\subsubsection{A Furstenberg type formula for linear drift}

We recall the idea of Furstenberg that the largest Lyapunov exponent 
\[\chi = \lim \frac{1}{n}\log(|A_0\cdots A_n|)\]
of a sequence $\ldots,A_{-1},A_0,A_1,\ldots$ of independent identically distributed random $2\times 2$ matrices with determinant equal to $1$ can be expressed as 
\[\chi = \E(\log(|A_0v|))\]
where $v$ is a random vector on the unit circle which is independent from $A_0$ (it belongs to the unstable direction of the sequence of matrices and depends only on $A_{-1},A_{-2},\ldots$).

The distribution of $v$ above is, a priori, unknown but must satisfy a stationarity property, and even without full knowledge of $v$ the formula can be used to establish that the exponent is positive under certain assumptions on the sequence of matrices.  This idea has been generalized considerably and there are now several `Furstenberg type formulas' which express asymptotic quantities of a random trajectory as integrals over a finite time segment of the trajectory involving some random element of a `boundary space' whose distribution is unknown (see for example \cite[Theorem 18]{karlsson-ledrappier2011}).

The purpose of this subsection is to establish a Furstenberg type formula for the linear drift of a harmonic measure (compare with \cite[Proposition 1.1]{ledrappier2010}).  The Busemann function space $\widehat{\manifolds}^b$ will play the role of the boundary (the circle in the above example) and our trajectory is the path of Brownian motion.  We extend the shift maps $\shift^t$ to $\widehat{\manifolds}^b$ so that $(M,\omega,\xi)$ goes to $(M,\omega_{t+\cdot}, \xi - \xi(\omega_t))$ and notice that they are continuous.

By a Brownian motion on a stationary random manifold $M$ we mean a random element $(M,\omega)$ of $\widehat{\manifolds}$ whose distribution is the shift invariant lift $\widehat{\mu}$ of the harmonic distribution $\mu$ of $M$ (see Theorem \ref{liftedharmonic}).   In the following theorem the existence of the random element $(M,\omega,\xi)$ extending $(M,\omega)$ may depend on modifying the domain probability space of $(M,\omega)$ somewhat (without changing the distribution of $(M,\omega)$).

\begin{theorem}[Furstenberg type formula for linear drift]\label{furstenbergformula}
Let $(M,\omega)$ be a Brownian motion on an ergodic stationary random manifold $M$.  Then \textnormal{(}possibly modifying the domain of $(M,\omega)$\textnormal{)} there exists a random generalized Busemann function $\xi$ such that the distribution of $(M,\omega,\xi)$ is shift invariant and one has
\[\E\left(\xi(\omega_t)\right) = t\ell(M)\]
for all $t \in \R \setminus \lbrace 0\rbrace$.  Furthermore $\xi$ can be chosen so that its distribution is ergodic for the shift maps and so that the conditional distribution of $\lbrace \omega_t: t \ge 0\rbrace$ given $M,o=\omega_0$ and $\xi$ is $P_o$.
\end{theorem}
\begin{proof}
Let $(M,\omega)$ be a Brownian motion on an ergodic stationary random manifold $M$ and let $u$ be a uniform random variable in $[0,1]$.  Because the projection $\widehat{\pi}$ from $\widehat{\manifolds}^b$ to $\widehat{\manifolds}$ is proper the distributions of the random variables $(M,\omega,\xi_{\omega_{-uT}})$ (where $T$ ranges over $\R$) are tight (i.e. for each $\epsilon \g 0$ there is a compact subset of $\widehat{\manifolds}^b$ with probability greater than $1-\epsilon$ for all the distributions).  Hence there exists a sequence $T_n \to +\infty$ such that $(M,\omega,\xi_{\omega_{-uT_n}})$ converges in distribution to some random element $(M,\omega,\xi)$ in $\widehat{\manifolds}^b$ (see \cite[Theorem 5.1]{billingsley1999}, notice that this element projects to our original Brownian motion $(M,\omega)$).

Since all functions in $\widehat{\manifolds}^b$ are $1$-Lipschitz we have that $|\xi(\omega_t)| \le d(\omega_0,\omega_t)$.  Notice also that, as shown in the previous section, the expectation of $d(\omega_0,\omega_t)$ is finite for any Brownian motion $(M,\omega)$ on a stationary random manifold.   Also, by Skorohod's representation theorem (see \cite[Theorem 6.7]{billingsley1999}) there exist random elements with the distribution of $(M,\omega,\xi_{\omega_{-uT_n}})$ which converge pointwise to a random element with the distribution of $(M,\omega,\xi)$.  Hence we may use dominated convergence to establish the first in the following chain of equalities (we assume for simplicity that $t \g 0$)
\begin{align*}\E\left(\xi(\omega_t)\right) &= \lim\limits_{n \to +\infty}\E\left(\xi_{\omega_{-uT_n}}(\omega_t)\right) = \lim\limits_{n \to +\infty}\E\left(d(\omega_{-uT_n},\omega_t) - d(\omega_{-uT_n},\omega_0)\right)
\\ &= \lim\limits_{n \to +\infty}\E\left(d(\omega_0,\omega_{uT_n+t}) - d(\omega_0,\omega_{uT_n})\right).
\end{align*}

Notice that $u+t/T_n$ is uniformly distributed on $[t/T_n,t/T_n+1]$. Therefore, setting $v = u$ if $u \g t/T_n$ and $v = 1+u$ otherwise one obtains that $v$ has the same distribution as $u+t/T_n$ from which it follows that
\[\E\left(\xi(\omega_t)\right) = \lim\limits_{n \to +\infty}\E\left(d(\omega_0,\omega_{vT_n}) - d(\omega_0,\omega_{uT_n})\right)\]
where the inner terms cancel except on a set of probability $t/T_n$ where $u \l t/T_n$.

Hence we have obtained
\[\E\left(\xi(\omega_t)\right) = t\lim\limits_{n \to +\infty}\E\left(\frac{d(\omega_0,\omega_{s+T_n})}{T_n}\right) - t\lim\limits_{n \to +\infty}\E\left(\frac{d(\omega_0,\omega_{s})}{T_n}\right)\]
where $s$ is independent from $(M,\omega)$ and uniformly distributed on $[0,t]$.

In second term we may bound the expected value of $d(\omega_0,\omega_s)$, using Ichihara's comparison result (see \cite{ichihara1988}), by the expected diameter of a segment of length $t$ of Brownian motion on a space of constant curvature lower than that of all possible values of $M$.  Since $T_n \to +\infty$ it follows that the second limit is $0$.  For the first term we can bound the expected value of $|d(\omega_0,\omega_{s+T_n}) - d(\omega_0,\omega_{T_n})|$ by the same.  This implies that one has
\[\E\left(\xi(\omega_t)\right) = t\lim\limits_{n \to +\infty}\E\left(\frac{d(\omega_0,\omega_{T_n})}{T_n}\right) = t\ell(M)\]
where $\ell(M)$ is the linear drift of the stationary random manifold $M$.  We omit the proof for negative $t$ which is very similar.

To see that the distribution of $(M,\omega,\xi)$ is shift invariant take any $s \g 0$ and notice that
\[\shift^{-s}(M,\omega,\xi_{\omega_{uT_n}}) = (M,\omega_{\cdot - s},\xi_{\omega_{uT_n}} - \xi_{\omega_{uT_n}}(\omega_{-s}))\]
has the same distribution as $(M,\omega,\xi_{\omega_{uT_n + s}})$.  We may define $v$ so that $uT_n + s$ has the same distribution as the $vT_n$ by setting $v = u$ if $uT_n \l s$ and $v = u+1$ otherwise.  Hence the $\shift^{-s}(M,\omega,\xi_{\omega_{uT_n}})$ has the same distribution as $(M,\omega,\xi_{vT_n})$ which coincides with $(M,\omega,\xi_{tT_n})$ outside of a set of probability $s/T_n$.  Taking limits when $n \to +\infty$ one obtains that $\shift^{-s}(M,\omega,\xi)$ has the same distribution as $(M,\omega,\xi)$ for all $s \g 0$, so the distribution of $(M,\omega,\xi)$ is shift invariant.

Since $u$ is independent from $(M,\omega)$ and the conditional distribution of $\lbrace \omega_t: t \ge 0\rbrace$ relative to $M,o = \omega_0$ is $P_o$ one obtains the same property for conditioning relative to $\xi_{\omega_{-uT_n}}$ and by taking limits also for $\xi$.

The possible shift invariant distributions of $(M,\omega,\xi)$ satisfying $\E(\xi(\omega_t)) = t\ell(M)$ for all $t \neq 0$ and such that the conditional distribution of $\lbrace \omega_t: t \ge 0\rbrace$ relative to $M,o = \omega_0, \xi$ are $P_{o}$, form a convex and weakly compact set.  Any extremal element of this set (and such an element exists by the Krein-Milman theorem) is ergodic with respect to the shift maps.  This shows that the distribution of $\xi$ can be chosen to be ergodic.
\end{proof}

\subsection{Entropy of reversed Brownian motion}
\subsubsection{Reversibility}

The purpose of this section is to improve the inequality $\frac{1}{2}\ell(M)^2 \le h(M)$ obtained in Theorem \ref{inequalities} in the case of Brownian motions on a stationary random Hadamard manifold to $2\ell(M)^2 \le h(M)$ (recall that a Hadamard manifold is a manifold isometric to $\R^d$ endowed with a complete Riemannian metric of non-positive sectional curvature).  This improvement was established by Kaimanovich and Ledrappier in the case of a single manifold with a compact quotient (no curvature assumption) and has strong rigidity consequences if for manifolds with negative curvature in this case (see \cite{ledrappier2010} and the references therein\begin{footnote}{Since we use $q(t,x,y) = p(t/2,x,y)$ our definitions of $\ell$ and $h$ differ from Ledrappier's by a factor of two.  Hence then inequality $h \le \ell v$ remains the same with both conventions but our claim that $2\ell^2 \le h$ corresponds to $\ell^2 \le h$ in Ledrappier's notation.}\end{footnote}).  We are able to prove the inequality 
under the 
assumption that Brownian motion is reversible; a technical hypothesis which is automatically satisfied in the case of a single manifold with compact quotient.  We will discuss this assumption briefly in this subsection.

Consider a Brownian motion $(M,\omega)$ on a stationary random manifold $M$.  The reverse process is defined by $(M,\omega')$ where $\omega'_t = \omega_{-t}$.  We say the Brownian motion is reversible if $(M,\omega)$ and $(M,\omega')$ have the same distribution.

As an example consider a compact manifold $M$.  The unique stationary measure for Brownian motion is the normalized volume measure hence the unique shift invariant measure on $C([0,+\infty),M)$ is $\int P_x \d x/\vol(M)$.  We use this to define a unique shift invariant measure on the continuous paths $C(\R,M)$ defined on all of $\R$.  Giving rise to a Brownian motion on a stationary random manifold $(M,\omega)$ where $M$ is fixed and the basepoint $\omega_0$ is uniformly distributed.

We claim that $(M,\omega)$ thus defined is reversible.  For this purpose notice that given Borel sets $A_0,A_1$ in $M$ and $t \g 0$ one has
\begin{align*}
\P\left(\omega_0 \in A_0, \omega_t \in A_1\right) = \int\limits_{A_0 \times A_1}q(t,x,y)\d x \d y/\vol(M) =  \P\left(\omega_0 \in A_1, \omega_t \in A_0\right)
\end{align*}
because $q(t,x,y) = q(t,y,x)$.  The claim follows from repeating this calculation for an arbitrary finite number of sets and times.

Given a compact foliation $X$ and a harmonic measure $\mu$ one can define a measure on the space of paths $C([0,+\infty),X)$ corresponding to `leafwise Brownian motion' with initial distribution $\mu$.  Similarly to the case discussed above where $X = M$ was a single compact leaf, the fact that $\mu$ is harmonic allows one to uniquely extend this probability in a shift invariant way to all of $C(\R,X)$.  The results of Deroin and Klepsyn (see \cite[Theorem B]{deroin-kleptsyn2007}) imply that for a minimal codimension one foliation without any transverse invariant measure the corresponding Brownian motion indexed on $\R$ is not reversible.  If, furthermore, the foliation is constructed so that no two leaves are isometric then one obtains a stationary random manifold with non-reversible Brownian motion.
\subsubsection{Furstenberg type formula for Hadamard manifolds with pinched negative curvature}

Recall that the Busemann functions of a manifold $(M,o,g)$ are defined as local uniform limits of the function of the form $\xi_x(y) = d(x,y) - d(x,o)$.  The Busemann functions which are not of the form $\xi_x$ form the so-called Busemann boundary of $M$.

\begin{lemma}\label{busemannlimit}
Let $(M,o,g)$ be a Hadamard manifold with curvature bounded between two negative constants.  Then for $P_o$ almost every Brownian path $\omega$ the following limit exists and is an element of the Busemann boundary
\[\xi = \lim\limits_{t \to +\infty}\xi_{\omega_t}.\]
\end{lemma}
\begin{proof}
Convergence of $\xi_{\omega_t}$ to a boundary Busemann function follows if one shows that $d(\omega_0,\omega_t) \to +\infty$ and that the geodesic segment joining $o$ to $\omega_t$ converges to a geodesic ray (see \cite[Proposition 2]{wang2011} and the references therein).  Both of these properties of Brownian motion were established by Prat in the mid 70s (see \cite[Theorem 3.2]{arnaudon-thalmaier2011}, \cite{prat1971} and \cite{prat1975}).
\end{proof}

The Martin boundary of a manifold $M$ consists of the limits of functions of the form $y \mapsto G(x,y)/G(x,o)$ where $G$ is a minimal Green's function (see \cite{wang2011}). The elements of the Martin boundary are positive harmonic functions which we choose to consider up to multiplication by positive constants.  In the context of Hadamard manifolds with pinched negative curvature the equivalence between Martin and Busemann boundaries was established by Anderson and Schoen. 

\begin{lemma}\label{reversedistribution}
Let $(M,o,g)$ be a Hadamard manifold with curvature bounded between two negative constants. Then there is a natural homeomorphism $\xi \mapsto k_{\xi}$ between the Busemann and Martin boundaries.  Furthermore the probability transition density of Brownian motion conditioned on the value of $\xi = \lim_{t\to +\infty} \xi_{\omega_t}$ is given by
\[\frac{k_{\xi}(y)}{k_\xi(x)}q(t,x,y).\]
\end{lemma}
\begin{proof}
For the homeomorphism between the two boundaries see \cite[Theorem 6.3]{anderson-schoen1985}.  The conditional distribution for Brownian motion is verified on page 36 of \cite{ancona1990} (this is a special case of the so-called $h$-transform due to Doob, see \cite{doob2001}).
\end{proof}

We can now established a refined version of Theorem \ref{furstenbergformula} for random Hadamard manifolds with pinched negative curvature.

\begin{lemma}\label{hadamardfurstenbergformula}
Let $(M,\omega)$ be a reversible Brownian motion on a stationary random Hadamard manifold with sectional curvatures bounded between two negative constants.  Then letting $\xi = \lim_{t \to +\infty}\xi_{\omega_{-t}}$ one has
\[\E(\xi(\omega_t)) = t \ell(M)\]
for all $t \neq 0$.
\end{lemma}
\begin{proof}
By Lemma \ref{busemannlimit} the distribution of $(M,\omega,\xi_{\omega_{-uT}})$ where $u$ is uniformly distributed in $[0,1]$ and independent from $(M,\omega)$, converges to that of $(M,\omega,\xi)$.  Hence the result follows exactly as in the proof of Theorem \ref{furstenbergformula}.
\end{proof}

Boundary Busemann functions are at least two times continuously differentiable on any Hadamard manifold (see \cite{heintze-imhof1977}) this allows us to pass to the limit when $t \to 0$ in the formulas above (this idea for obtaining infinitesimal formulas for the linear drift goes back to \cite{kaimanovich1986}).  Doing so along positive and negative $t$ yields different formulas, in order to use Lemma \ref{reversedistribution} to obtain the distribution of reversed Brownian motion conditioned to $\xi$ we must impose the hypothesis that our Brownian motion is reversible.

\begin{lemma}\label{theformula}
Let $(M,\omega)$ be a reversible Brownian motion on a stationary random Hadamard manifold with sectional curvatures bounded between two negative constants and let $\xi = \lim_{t \to +\infty}\xi_{\omega_{-t}}$.  Then one has
\[\ell(M) = \E\left(\frac{1}{2}\Delta\xi(o_M)\right) = -\E\left(\frac{1}{2}\langle \nabla \log k_\xi(o_M), \nabla \xi(o_M)\rangle\right).\]
\end{lemma}
\begin{proof}
By Lemma \ref{hadamardfurstenbergformula} on has
\[\ell(M) = \E\left( \frac{1}{t}\int q(t,o_M,x)\xi(x)\d x\right)\]
for all $t \g 0$.  Since the Busemann functions $\xi$ are $1$-Lipschitz and $\int q(t,o_M,x)d(o_M,x)\d x$ can be uniformly bounded on the support of the distribution of $M$, one can take limit when $t \to 0^+$ from which one obtains that
\[\ell(M) = \E\left(\frac{1}{2}\Delta\xi(o_M)\right).\]

Similarly, by our assumption of reversibility, $\omega_{-t}$ has distribution $q(t,o_M,x)$ on $M$.  The conditional distribution with respect to $\xi$ is given by Lemma \ref{reversedistribution} and one has
\[\ell(M) = -\E\left( \frac{1}{t}\int \frac{k_\xi(x)}{k_\xi(o_M)}q(t,o_M,x)\xi(x)\d x\right)\]
for all $t \g 0$.  Passing to the limit with $t \to 0^-$ one obtains\begin{footnote}{Notice that $k_\xi(x)/k_\xi(o_M)$ is bounded by $C\exp(Cd(o_M,x))$ for some $C$ (see \cite[Corollary 4.5]{arnaudon-driver-thalmaier2007}) and $\xi$ is Lipschitz so one has a uniform bound for the inner integral on all manifolds in the support of the distribution of $M$ by virtue of the uniform upper heat kernel bounds.}\end{footnote}
\[\ell(M) = -\E\left(\frac{1}{2}\Delta\xi(o_M) + \langle \nabla \log k_\xi(o_M) ,\nabla \xi(o_M)\rangle\right)\]
from which the second claimed formula for $\ell(M)$ follows using the first.
\end{proof}

As a toy example of the previous formulas for drift consider the hyperbolic half plane $M = \lbrace (x,y) \in \R^2: y \g 0\rbrace$ with the metric $\d s^2 = y^{-2}(\d x^2 + \d y^2)$ and base point $o_M = (0,1)$.  Notice that $(M,o_M)$ is a stationary random manifold so one may indeed apply Lemma \ref{theformula}.  All Busemann functions are obtained by applying isometries to $\xi(x,y) = -\log(y)$.  Calculating the Laplacian of $\xi$ at $o_M$ one obtains $1$ so $\ell(M) = 1/2$.  On the other hand the positive harmonic function associated to $\xi$ (the Poisson kernel function associated to the boundary point at infinity) is $k_\xi(x,y) = y$ so that in the second formula $\nabla \log k_\xi = -\nabla \xi$ and one obtains again $\ell(M) = 1/2$.

As stated before our objective in this section is to obtain the inequality $2\ell(M)^2 \le h(M)$ for stationary random Hadamard manifolds with curvature bounded between two negative constants and reversible Brownian motion.  So far we have the following.

\begin{corollary}\label{protoinequality}
Let $(M,\omega)$ be a reversible Brownian motion on a stationary random Hadamard manifold with sectional curvatures bounded between two negative constants and let $\xi = \lim_{t \to +\infty}\xi_{\omega_{-t}}$.  Then one has
\[2\ell(M)^2 \le \E\left(\frac{1}{2}|\nabla \log k_\xi(o_M)|^2\right)\]
\end{corollary}
\begin{proof}
Using Lemma \ref{theformula} followed by Jensen's inequality and the Cauchy-Schwarz inequality one obtains
\[2\ell(M)^2 \le \E\left(\frac{1}{2}|\langle \nabla \log k_\xi(o_M),\nabla \xi(o_M)\rangle|^2\right) \le \E\left(\frac{1}{2}|\nabla \log k_\xi(o_M)|^2\right)\]
where we have also used the fact that $\xi$ is $1$-Lipschitz.
\end{proof}
\subsubsection{Reverse entropy and entropy difference}

The purpose of this section is to understand the right-hand side of the inequality in Corollary \ref{protoinequality}.  We want to show that it is smaller than $h(M)$ in order to obtain the inequality $2\ell(M)^2 \le h(M)$.  In short the proof consists in establishing that the right hand side $k(M)$ is the difference between the entropy $h(M)$ and an entropy associated to the reversed process (conditioned on $\xi = \lim\limits_{t \to +\infty}\xi_{\omega_{-t}}$), since this `reversed entropy' is non-negative one obtains $k(M) \le h(M)$.  As a first step we obtain an alternate formula for $k(M)$ as the expected increment of $\log k_\xi$ along a Brownian path.

\begin{lemma}
Let $(M,\omega)$ be a reversible Brownian motion on a stationary random Hadamard manifold with curvature bounded between two negative constants, and let $\xi = \lim\limits_{t \to +\infty}\xi_{\omega_{-t}}$.  Then setting 
\[k(M) = \E\left(\frac{1}{2}|\nabla \log k_\xi(o_M)|^2\right)\]
one has
\[k(M) = -\E\left( \frac{1}{t}\log\left(\frac{k_\xi(\omega_t)}{k_\xi(o_M)}\right)\right) = \E\left( \frac{1}{t}\log\left(\frac{k_\xi(\omega_{-t})}{k_\xi(o_M)}\right)\right)\]
for all $t \g 0$.
\end{lemma}
\begin{proof}
Using that $k_\xi$ is shift invariant one obtains for
\[K_t = -\E\left(\log\left(\frac{k_\xi(\omega_{t})}{k_\xi(o_M)}\right)\right)\]
that $K_{t+s} = K_t + K_s$.

Furthermore since there is a uniform bound for $|\grad \log k_\xi|$ over the entire support of $(M,\omega,\xi)$ (see \cite[Corollary 4.4]{arnaudon-driver-thalmaier2007}) one obtains that $K_t$ is continuous with respect to $t$ by dominated convergence.  This implies $K_t = tK_1$ for all $t \g 0$.

In particular one has for all $t \g 0$ that
\[K_1 = -\E\left(\frac{1}{t}\int q(t,o_M,x)\log\left(\frac{k_\xi(x)}{k_\xi(o_M)}\right)\d x\right)\]
and taking limit when $t \to 0^+$  (which is justified by dominated convergence, again by the uniform bounds of \cite[Corollary 4.4]{arnaudon-driver-thalmaier2007}) one obtains
\[K_1 = \E\left(\frac{1}{2}\Delta \log k_\xi(o_M)\right) = \E\left(\frac{1}{2}|\nabla \log k_\xi(o_M)|^2\right)\]
as claimed.

The equality
\[-\E\left( \frac{1}{t}\log\left(\frac{k_\xi(\omega_t)}{k_\xi(o_M)}\right)\right) = \E\left( \frac{1}{t}\log\left(\frac{k_\xi(\omega_{-t})}{k_\xi(o_M)}\right)\right)\]
follows by shift invariance of $k_\xi$.
\end{proof}

Given a Hadamard manifold $(M,o,g)$ with curvature bounded by two negative constants and a boundary Busemann function $\xi$ let $P_{x,\xi}$ be the probability measure on $C([0,+\infty),M)$ which is the distribution of Brownian motion conditioned to exit at $\xi$ (the transition probability densities are given by Lemma \ref{reversedistribution}).  Recall that we denote by $\F_t$ and $\F^T$ the $\sigma$-algebras generated by $\omega_s$ with $s \le t$ and $s \ge T$ respectively and by $\F^\infty$ the tail $\sigma$-algebra on $C([0,\infty),M)$.

In this context let $I_t^T(M,\xi)$ (where $0 \l t \l T \le \infty$) be the mutual information between $\F_t$ and $\F^T$ with respect to the probability $P_{o,\xi}$.

\begin{lemma}\label{reverseinformation}
Let $(M,o_M,g)$ be a Hadamard manifold with curvature strictly bounded between two negative constants and $\xi$ be a boundary Busemann function.  Then the following properties hold for all $0 \l t \l T \l \infty$:
\begin{enumerate}
 \item $I_t^T(M,\xi)= \int \log\left( \frac{k_\xi(o_M)}{k_\xi(x)}\frac{q(T-t,x,y)}{q(T,o_M,y)}\right) \frac{k_\xi(y)}{k_\xi(o_M)}q(t,o,x)q(T-T,x,y)\d x \d y$.
 \item The function $T \mapsto I_t^T(M,\xi)$ is non-negative and non-increasing.
\end{enumerate}
\end{lemma}
\begin{proof}
The formula for $I_t^T(M,\xi)$ follows from Lemma \ref{reversedistribution} and the Gelfand-Yaglom-Peres theorem (see the proof of Theorem \ref{informationtheorem}).  Non-negativity and monotonicity follow directly form the definition of mutual information.
\end{proof}

The following result identifies $k(M)$ as the difference between $h(M)$ and an entropy for the reverse process of $(M,\omega)$ conditioned to its limit Busemann function $\xi$.

\begin{lemma}\label{reverseentropy}
Let $(M,\omega)$ be a reversible Brownian motion on a stationary random Hadamard manifold with curvature bounded between two negative constants, and let $\xi = \lim\limits_{t \to +\infty}\xi_{\omega_{-t}}$.  Then for all $t \g 0$ one has
\[0 \le \lim\limits_{T \to +\infty}\E\left(I_t^T(M,\xi)\right) = t(h(M) - k(M)).\]
In particular $k(M) \le h(M)$.
\end{lemma}
\begin{proof}
We calculate using Lemma \ref{reverseinformation}, Lemma \ref{reversedistribution}, and the fact that $(M,\omega,\xi)$ is shift invariant to obtain
\begin{align*}
0 \le \E\left(I_t^T(M,\xi)\right) &= \E\left(\log\left(\frac{k_\xi(\omega_0)}{k_\xi(\omega_{-t})}\right) + \log(q(T-t,\omega_{-t},\omega_{-T}) - \log(q(T,\omega_0,\omega_{-T})\right)
\\ &= \E\left(\log\left(\frac{k_\xi(\omega_t)}{k_\xi(\omega_{0})}\right) + \log(q(T-t,\omega_{0},\omega_{T-t})) - \log(q(T,\omega_0,\omega_{T}))\right)
\\ &= tk(M) - H_{T-t} + H_T
\end{align*}
where $H_t = \E(\int \log(q(t,o_M,x))q(t,o_M,x)\d x)$.

When $T \to +\infty$ one has that $H_T-H_{T-t}$ converges to $th(M)$ (see the proof of Theorem \ref{entropytheorem}) and therefore one has
\[0 \le \E\left(I_t^T(M,\xi)\right) = t(h(M) - k(M))\]
as claimed.
\end{proof}

To conclude we combine the previous results to obtain a sharp lower bound for the entropy $h(M)$ of a stationary random Hadamard manifold in terms of its linear drift $\ell(M)$.

\begin{theorem}\label{lsquarelessthanh}
Let $(M,\omega)$ be a reversible Brownian motion on a stationary random Hadamard manifold with curvature bounded between two negative constants, and let $\xi = \lim\limits_{t \to +\infty}\xi_{\omega_{-t}}$.  Then one has $2\ell(M)^2 \le h(M)$ with equality if and only if 
\[\nabla \log k_\xi(o_M) = -2\ell(M)\nabla\xi(o_M)\]
almost surely.
\end{theorem}
\begin{proof}
By Lemma \ref{theformula} one has
\[\ell(M) = -\E\left(\frac{1}{2}\langle \nabla \log k_\xi(o_M), \nabla \xi(o_M)\rangle \right).\]

Squaring and applying Jensen's inequality one obtains
\[2\ell(M)^2 \le \E\left(\frac{1}{2}|\langle \nabla \log k_\xi(o_M), \nabla \xi(o_M)\rangle|^2 \right).\]

Notice at this point that if the equality $\E(X)^2 = \E(X^2)$ holds for some random variable $X$ then $X$ is almost surely constant.  Hence if equality holds in the last inequality above one obtains that $\ell(M) = -\frac{1}{2}\langle \nabla \log k_\xi(o_M), \nabla \xi(o_M)\rangle$ almost surely.

Next we apply the Cauchy-Schwartz inequality and the fact that $\xi$ is $1$-Lipschitz to obtain
\[\E\left(\frac{1}{2}|\langle \nabla \log k_\xi(o_M), \nabla \xi(o_M)\rangle|^2 \right) \le k(M)\]
and by Lemma \ref{reverseentropy} conclude that $2\ell(M)^2 \le k(M) \le h(M)$.

At this step we notice that equality would imply that $\nabla\log k_\xi(o_M)$ and $\nabla \xi(o_M)$ are collinear almost surely.  Combining this with the previous observation about equality in the Jensen's inequality one sees that if $2\ell(M)^2 = h(M)$ then one would have
\[\nabla\log k_\xi(o_M) = -2\ell(M)\nabla\xi(o_M)\]
almost surely as claimed.
\end{proof}

\def\cprime{$'$}


\begin{thebibliography}{{Gra}11}

\bibitem[{\'A}C03]{alvarez-candel2003}
J.~A. {{\'A}lvarez L{\'o}pez} and A.~{Candel}.
\newblock Generic aspects of the geometry of leaves of foliations.
\newblock \url{http://www.csun.edu/~ac53971/research/kyoto_GF2003.pdf}, 2003.

\bibitem[ADT07]{arnaudon-driver-thalmaier2007}
Marc Arnaudon, Bruce~K. Driver, and Anton Thalmaier.
\newblock Gradient estimates for positive harmonic functions by stochastic
  analysis.
\newblock {\em Stochastic Process. Appl.}, 117(2):202--220, 2007.

\bibitem[Anc90]{ancona1990}
A.~Ancona.
\newblock Th\'eorie du potentiel sur les graphes et les vari\'et\'es.
\newblock In {\em \'{E}cole d'\'et\'e de {P}robabilit\'es de {S}aint-{F}lour
  {XVIII}---1988}, volume 1427 of {\em Lecture Notes in Math.}, pages 1--112.
  Springer, Berlin, 1990.

\bibitem[AS85]{anderson-schoen1985}
Michael~T. Anderson and Richard Schoen.
\newblock Positive harmonic functions on complete manifolds of negative
  curvature.
\newblock {\em Ann. of Math. (2)}, 121(3):429--461, 1985.

\bibitem[AT11]{arnaudon-thalmaier2011}
Marc Arnaudon and Anton Thalmaier.
\newblock Brownian motion and negative curvature.
\newblock In {\em Random walks, boundaries and spectra}, volume~64 of {\em
  Progr. Probab.}, pages 143--161. Birkh\"auser/Springer Basel AG, Basel, 2011.

\bibitem[Ave76]{avez1976}
A.~Avez.
\newblock Harmonic functions on groups.
\newblock In {\em Differential geometry and relativity}, pages 27--32.
  Mathematical Phys. and Appl. Math., Vol. 3. Reidel, Dordrecht, 1976.

\bibitem[AW06]{aizenman-warzel2006}
Michael Aizenman and Simone Warzel.
\newblock The canopy graph and level statistics for random operators on trees.
\newblock {\em Math. Phys. Anal. Geom.}, 9(4):291--333 (2007), 2006.

\bibitem[BBI01]{burago-burago-ivanov2001}
Dmitri Burago, Yuri Burago, and Sergei Ivanov.
\newblock {\em A course in metric geometry}, volume~33 of {\em Graduate Studies
  in Mathematics}.
\newblock American Mathematical Society, Providence, RI, 2001.

\bibitem[BC12]{benjamini-curien2012}
Itai Benjamini and Nicolas Curien.
\newblock Ergodic theory on stationary random graphs.
\newblock {\em Electron. J. Probab.}, 17:no. 93, 1--20, 2012.

\bibitem[Bil99]{billingsley1999}
Patrick Billingsley.
\newblock {\em Convergence of probability measures}.
\newblock Wiley Series in Probability and Statistics: Probability and
  Statistics. John Wiley \& Sons Inc., New York, second edition, 1999.
\newblock A Wiley-Interscience Publication.

\bibitem[BSW12]{brofferio-salvatori-woess2012}
Sara Brofferio, Maura Salvatori, and Wolfgang Woess.
\newblock Brownian motion and harmonic functions on {${\rm Sol}(p,q)$}.
\newblock {\em Int. Math. Res. Not. IMRN}, (22):5182--5218, 2012.

\bibitem[Can03]{candel2003}
Alberto Candel.
\newblock The harmonic measures of {L}ucy {G}arnett.
\newblock {\em Adv. Math.}, 176(2):187--247, 2003.

\bibitem[Cha84]{chavel1984}
Isaac Chavel.
\newblock {\em Eigenvalues in {R}iemannian geometry}, volume 115 of {\em Pure
  and Applied Mathematics}.
\newblock Academic Press Inc., Orlando, FL, 1984.
\newblock Including a chapter by Burton Randol, With an appendix by Jozef
  Dodziuk.

\bibitem[Cri08]{cristina2008}
Jan Cristina.
\newblock Gromov-hausdorff convergence of metric spaces.
\newblock \url{http://www.helsinki.fi/~cristina/pdfs/gromovHausdorff.pdf},
  2008.

\bibitem[Der76]{derriennic1976}
Yves Derriennic.
\newblock Lois ``z\'ero ou deux'' pour les processus de {M}arkov.
  {A}pplications aux marches al\'eatoires.
\newblock {\em Ann. Inst. H. Poincar\'e Sect. B (N.S.)}, 12(2):111--129, 1976.

\bibitem[Der85]{derriennic1985}
Y.~Derriennic.
\newblock {\em Entropie, th\'eor\`emes limite et marches al\'eatoires}.
\newblock Publications de l'Institut de Recherche Math\'ematique de Rennes.
  [Publications of the Rennes Institute of Mathematical Research]. Universit\'e
  de Rennes I Institut de Recherche Math\'ematique de Rennes, Rennes, 1985.

\bibitem[DGM77]{debiard-gaveau-mazet1976}
A.~Debiard, B.~Gaveau, and E.~Mazet.
\newblock Th\'eor\`emes de comparaison en g\'eom\'etrie riemannienne.
\newblock {\em Publ. Res. Inst. Math. Sci.}, 12(2):391--425, 1976/77.

\bibitem[DK07]{deroin-kleptsyn2007}
Bertrand Deroin and Victor Kleptsyn.
\newblock Random conformal dynamical systems.
\newblock {\em Geom. Funct. Anal.}, 17(4):1043--1105, 2007.

\bibitem[dlR93]{delarue1993}
Thierry de~la Rue.
\newblock Espaces de {L}ebesgue.
\newblock In {\em S\'eminaire de {P}robabilit\'es, {XXVII}}, volume 1557 of
  {\em Lecture Notes in Math.}, pages 15--21. Springer, Berlin, 1993.

\bibitem[DM88]{davies-mandouvalos1988}
E.~B. Davies and N.~Mandouvalos.
\newblock Heat kernel bounds on hyperbolic space and {K}leinian groups.
\newblock {\em Proc. London Math. Soc. (3)}, 57(1):182--208, 1988.

\bibitem[Doo01]{doob2001}
Joseph~L. Doob.
\newblock {\em Classical potential theory and its probabilistic counterpart}.
\newblock Classics in Mathematics. Springer-Verlag, Berlin, 2001.
\newblock Reprint of the 1984 edition.

\bibitem[EW11]{einsiedler-thomas2011}
Manfred Einsiedler and Thomas Ward.
\newblock {\em Ergodic theory with a view towards number theory}, volume 259 of
  {\em Graduate Texts in Mathematics}.
\newblock Springer-Verlag London Ltd., London, 2011.

\bibitem[Fur63]{furstenberg1963}
Harry Furstenberg.
\newblock Noncommuting random products.
\newblock {\em Trans. Amer. Math. Soc.}, 108:377--428, 1963.

\bibitem[FV10]{friz-victoir2010}
Peter~K. Friz and Nicolas~B. Victoir.
\newblock {\em Multidimensional stochastic processes as rough paths}, volume
  120 of {\em Cambridge Studies in Advanced Mathematics}.
\newblock Cambridge University Press, Cambridge, 2010.
\newblock Theory and applications.

\bibitem[Gar81]{garnett1981}
Lucy~Jane Garnett.
\newblock {\em Functions and measures harmonic along the leaves of a foliation
  and the ergodic theorem}.
\newblock ProQuest LLC, Ann Arbor, MI, 1981.
\newblock Thesis (Ph.D.)--Dartmouth College.

\bibitem[Gar83]{garnett1983}
Lucy Garnett.
\newblock Foliations, the ergodic theorem and {B}rownian motion.
\newblock {\em J. Funct. Anal.}, 51(3):285--311, 1983.

\bibitem[Ghy95]{ghys1995}
{\'E}tienne Ghys.
\newblock Topologie des feuilles g\'en\'eriques.
\newblock {\em Ann. of Math. (2)}, 141(2):387--422, 1995.

\bibitem[GMM12]{gouezel-matheus-maucourant2012}
S.~{Gou{\"e}zel}, F.~{Math{\'e}us}, and F.~{Maucourant}.
\newblock {Sharp lower bounds for the asymptotic entropy of symmetric random
  walks}.
\newblock {\em ArXiv e-prints}, September 2012.

\bibitem[{Gra}11]{gray2011}
Robert~M. {Gray}.
\newblock {\em {Entropy and information theory. 2nd ed.}}
\newblock New York, NY: Springer, 2nd ed. edition, 2011.

\bibitem[Gri09]{grigoryan2009}
Alexander Grigor'yan.
\newblock {\em Heat kernel and analysis on manifolds}, volume~47 of {\em AMS/IP
  Studies in Advanced Mathematics}.
\newblock American Mathematical Society, Providence, RI, 2009.

\bibitem[Gro81]{gromov1981}
Mikhael Gromov.
\newblock Groups of polynomial growth and expanding maps.
\newblock {\em Inst. Hautes \'Etudes Sci. Publ. Math.}, (53):53--73, 1981.

\bibitem[HIH77]{heintze-imhof1977}
Ernst Heintze and Hans-Christoph Im~Hof.
\newblock Geometry of horospheres.
\newblock {\em J. Differential Geom.}, 12(4):481--491 (1978), 1977.

\bibitem[Hsu02]{hsu2002}
Elton~P. Hsu.
\newblock {\em Stochastic analysis on manifolds}, volume~38 of {\em Graduate
  Studies in Mathematics}.
\newblock American Mathematical Society, Providence, RI, 2002.

\bibitem[Ich88]{ichihara1988}
Kanji Ichihara.
\newblock Comparison theorems for {B}rownian motions on {R}iemannian manifolds
  and their applications.
\newblock {\em J. Multivariate Anal.}, 24(2):177--188, 1988.

\bibitem[Ka{\u\i}86]{kaimanovich1986}
V.~A. Ka{\u\i}manovich.
\newblock Brownian motion and harmonic functions on covering manifolds. {A}n
  entropic approach.
\newblock {\em Dokl. Akad. Nauk SSSR}, 288(5):1045--1049, 1986.

\bibitem[Ka{\u\i}88]{kaimanovich1988}
V.~A. Ka{\u\i}manovich.
\newblock Brownian motion on foliations: entropy, invariant measures, mixing.
\newblock {\em Funktsional. Anal. i Prilozhen.}, 22(4):82--83, 1988.

\bibitem[Kai92]{kaimanovich1992}
Vadim~A. Kaimanovich.
\newblock Measure-theoretic boundaries of {M}arkov chains, {$0$}-{$2$} laws and
  entropy.
\newblock In {\em Harmonic analysis and discrete potential theory ({F}rascati,
  1991)}, pages 145--180. Plenum, New York, 1992.

\bibitem[Kec95]{kechris1995}
Alexander~S. Kechris.
\newblock {\em Classical descriptive set theory}, volume 156 of {\em Graduate
  Texts in Mathematics}.
\newblock Springer-Verlag, New York, 1995.

\bibitem[Kin68]{kingman1968}
J.~F.~C. Kingman.
\newblock The ergodic theory of subadditive stochastic processes.
\newblock {\em J. Roy. Statist. Soc. Ser. B}, 30:499--510, 1968.

\bibitem[KL07]{karlsson-ledrappier2007}
Anders Karlsson and Fran{\c{c}}ois Ledrappier.
\newblock Propri\'et\'e de {L}iouville et vitesse de fuite du mouvement
  brownien.
\newblock {\em C. R. Math. Acad. Sci. Paris}, 344(11):685--690, 2007.

\bibitem[KL11]{karlsson-ledrappier2011}
Anders Karlsson and Fran{\c{c}}ois Ledrappier.
\newblock Noncommutative ergodic theorems.
\newblock In {\em Geometry, rigidity, and group actions}, Chicago Lectures in
  Math., pages 396--418. Univ. Chicago Press, Chicago, IL, 2011.

\bibitem[KV83]{kaimanovich-vershik1983}
V.~A. Ka{\u\i}manovich and A.~M. Vershik.
\newblock Random walks on discrete groups: boundary and entropy.
\newblock {\em Ann. Probab.}, 11(3):457--490, 1983.

\bibitem[KW02]{kaimanovich-woess2002}
Vadim~A. Kaimanovich and Wolfgang Woess.
\newblock Boundary and entropy of space homogeneous {M}arkov chains.
\newblock {\em Ann. Probab.}, 30(1):323--363, 2002.

\bibitem[Led84]{ledrappier1984}
F.~Ledrappier.
\newblock Quelques propri\'et\'es des exposants caract\'eristiques.
\newblock In {\em \'{E}cole d'\'et\'e de probabilit\'es de {S}aint-{F}lour,
  {XII}---1982}, volume 1097 of {\em Lecture Notes in Math.}, pages 305--396.
  Springer, Berlin, 1984.

\bibitem[Led96]{ledrappier1996}
F.~Ledrappier.
\newblock Profil d'entropie dans le cas continu.
\newblock {\em Ast\'erisque}, (236):189--198, 1996.
\newblock Hommage {\`a} P. A. Meyer et J. Neveu.

\bibitem[Led10]{ledrappier2010}
Fran\c{c}ois Ledrappier.
\newblock Linear drift and entropy for regular covers.
\newblock {\em Geom. Funct. Anal.}, 20(3):710--725, 2010.

\bibitem[{Les}14]{lessa2013}
P.~{Lessa}.
\newblock {Reeb stability and the Gromov-Hausdorff limits of leaves in compact
  foliations}.
\newblock {\em To appear in Asian Journal of Mathematics}, 2014.

\bibitem[LS84]{lyons-sullivan1984}
Terry Lyons and Dennis Sullivan.
\newblock Function theory, random paths and covering spaces.
\newblock {\em J. Differential Geom.}, 19(2):299--323, 1984.

\bibitem[LS12]{ledrappier-shu2013}
F.~{Ledrappier} and L.~{Shu}.
\newblock {Entropy rigidity of symmetric spaces without focal points}.
\newblock {\em ArXiv e-prints}, May 2012.

\bibitem[Lu12]{lu2012}
Peng Lu.
\newblock Convergence of fundamental solutions of linear parabolic equations
  under {C}heeger-{G}romov convergence.
\newblock {\em Math. Ann.}, 353(1):193--217, 2012.

\bibitem[Pet06]{petersen2006}
Peter Petersen.
\newblock {\em Riemannian geometry}, volume 171 of {\em Graduate Texts in
  Mathematics}.
\newblock Springer, New York, second edition, 2006.

\bibitem[Pin64]{pinsker1964}
M.S. Pinsker.
\newblock {\em {Information and information stability of random variables and
  processes.}}
\newblock {Holden-Day Series in Time Series Analysis. San
  Francisco-London-Amsterdam: Holden-Day, Inc. XII. 243 p. }, 1964.

\bibitem[Pra71]{prat1971}
Jean-Jacques Prat.
\newblock \'{E}tude asymptotique du mouvement brownien sur une vari\'et\'e
  riemannienne \`a courbure n\'egative.
\newblock {\em C. R. Acad. Sci. Paris S\'er. A-B}, 272:A1586--A1589, 1971.

\bibitem[Pra75]{prat1975}
Jean-Jacques Prat.
\newblock \'{E}tude asymptotique et convergence angulaire du mouvement brownien
  sur une vari\'et\'e \`a courbure n\'egative.
\newblock {\em C. R. Acad. Sci. Paris S\'er. A-B}, 280(22):Aiii, A1539--A1542,
  1975.

\bibitem[Roh52]{rohlin1952}
V.~A. Rohlin.
\newblock On the fundamental ideas of measure theory.
\newblock {\em Amer. Math. Soc. Translation}, 1952(71):55, 1952.

\bibitem[Sol70]{solovay1970}
Robert~M. Solovay.
\newblock A model of set-theory in which every set of reals is {L}ebesgue
  measurable.
\newblock {\em Ann. of Math. (2)}, 92:1--56, 1970.

\bibitem[Var86]{varopoulos1986}
Nicholas~Th. Varopoulos.
\newblock Information theory and harmonic functions.
\newblock {\em Bull. Sci. Math. (2)}, 110(4):347--389, 1986.

\bibitem[Wan11]{wang2011}
Xiaodong Wang.
\newblock Compactifications of complete {R}iemannian manifolds and their
  applications.
\newblock In {\em Surveys in geometric analysis and relativity}, volume~20 of
  {\em Adv. Lect. Math. (ALM)}, pages 517--529. Int. Press, Somerville, MA,
  2011.

\end{thebibliography}
\end{document}